\documentclass[11pt]{amsart}
\usepackage[margin=1.1in]{geometry}
\usepackage{amsmath,amssymb,amsthm,mathtools}
\usepackage{enumitem}
\usepackage{tikz}
\usetikzlibrary{positioning,arrows.meta,calc}
\usepackage{graphicx}
\usepackage{xcolor}
\usepackage[hidelinks]{hyperref}
\usepackage[normalem]{ulem}

\newtheorem{theorem}{Theorem}[section]
\newtheorem{proposition}{Proposition}[section]
\newtheorem{lemma}{Lemma}[section]
\newtheorem{corollary}{Corollary}[section]
\theoremstyle{definition}
\newtheorem{definition}{Definition}[section]
\newtheorem{example}{Example}[section]
\newtheorem{assumption}{Assumption}[section]
\theoremstyle{remark}
\newtheorem{remark}{Remark}[section]
\numberwithin{equation}{section}
\numberwithin{figure}{section}

\newcommand{\R}{\mathbb{R}}
\newcommand{\E}{\mathbb{E}}
\newcommand{\Prob}{\mathbb{P}}
\newcommand{\Q}{\mathbb{Q}}
\newcommand{\Qt}{\widetilde{\mathbb{Q}}}
\newcommand{\Lt}{\widetilde{L}}
\newcommand{\mut}{\widetilde{\mu}}
\newcommand{\Ft}{\widetilde{F}}
\newcommand{\cE}{\mathcal{E}}
\newcommand{\ip}[2]{\langle #1,#2\rangle}
\newcommand{\tr}{\operatorname{tr}}
\newcommand{\dv}{\operatorname{div}}
\newcommand{\Sig}{\Sigma}
\DeclareMathOperator{\dist}{dist}

\title[Martingale property of multi-dimensional stochastic exponentials]{Necessary and sufficient conditions for the martingale property of stochastic exponentials driven by multi-dimensional diffusions}
\author{Erhan Bayraktar}
\thanks{E. Bayraktar is supported in part by the National Science Foundation under grant DMS-2602036.}
\address{Department of Mathematics, University of Michigan, Ann Arbor, MI 48109, USA}
\email{erhan@umich.edu}
\author{Hao Xing}
\address{Questrom School of Business, Boston University, 595 Commonwealth Avenue, Boston, MA 02215, USA}
\email{haoxing@bu.edu}
\date{November 2026}

\begin{document}
\begin{abstract}
We give necessary and sufficient coefficient conditions for stochastic exponentials driven by multidimensional diffusions to be true martingales. The classical reduction to non-exit of the diffusion with Girsanov-shifted drift separates the analysis into smooth boundary points, infinity, and nonsmooth boundary pieces. On smooth domains, under regularity and non-explosion hypotheses, martingality is equivalent to nonpositivity of the shifted Fichera drift, provided an explicit local condition holds where that drift vanishes and the stochastic integrand is unbounded. At infinity, Feller's test gives exact criteria when a scalar function of the diffusion reduces to a one-dimensional diffusion up to a time change with rate comparable to a positive function of the scalar coordinate. General radial and projection tests supply sufficient conditions. Applications give explicit necessary and sufficient conditions for inverse-CIR and $4/2$ models with correlated factors and singular stochastic integrands, including critical equality cases. The coupled non-affine examples go beyond direct scalar Feller tests and affine Riccati criteria, and exhibit boundary exits not decided by tests at infinity. Further exact criteria cover logarithmic boundary singularities, weighted distances, and an absorbed factor. The martingale conditions also yield uniqueness criteria for associated valuation equations.
\end{abstract}
\maketitle

\section{Introduction}\label{sec:intro}

Let $D\subseteq\R^d$ be a domain. We work on the canonical path space carrying the coordinate pair $(Y,W)$, with possible finite lifetime and a cemetery state, as specified in Section~\ref{sec:canonical}. Under $\Prob_y$, $W$ is an $m$-dimensional Brownian motion and $Y$ solves
\begin{equation}\label{eq:Y}
dY_t=\mu(Y_t)\,dt+\sigma(Y_t)\,dW_t,\qquad Y_0=y\in D,
\end{equation}
with $\mu:D\to\R^d$, $\sigma:D\to\R^{d\times m}$ and $a:=\sigma\sigma^{\top}$. Given $\theta:D\to\R^m$, we consider the positive local martingale
\begin{equation}\label{eq:S}
S_t=\cE\Big(\int_0^{\cdot}\theta(Y_s)\cdot dW_s\Big)_t=\exp\Big(\int_0^t\theta(Y_s)\cdot dW_s-\tfrac12\int_0^t|\theta(Y_s)|^2ds\Big)
\end{equation}
and ask when it is a true martingale. Equivalently, $dS_t/S_t=\theta(Y_t)\cdot dW_t$ with $S_0=1$. When $S$ is an asset price, $(S,Y)$ is a multifactor stochastic volatility model; the coefficients $\sigma$ and $\theta$ determine the correlations between the asset return and the factors. This framework contains the one-factor models of \cite{Sin98,AP07,LM07,BKX12,MU12} and the vector models of Dandapani--Protter \cite{DP22}. When $S$ represents a discounted asset price under an equivalent local martingale measure, its strict local martingality corresponds to an asset-price bubble \cite{CH05}. The martingale property also enters the study of optimal arbitrage in the sense of Fernholz--Karatzas \cite{FK10}, the valuation of American options \cite{BKX12b}, and uniqueness for the valuation equation in the natural growth class \cite{ET09,BX10,BKX12}.

In one dimension the question is settled: for a scalar factor, $S$ is a true martingale if and only if an auxiliary scalar diffusion does not reach the endpoints of its state interval, which Feller's test decides in terms of the coefficients \cite{DS02,MU12,BKX12}. In the affine class it is settled as well, through the global existence of solutions of generalized Riccati equations \cite{KMK10,MMS11,KRM15}. For multidimensional radial models, Criens--Glau \cite[Cor.~5.2]{CG18} give necessary and sufficient martingality conditions through one-dimensional diffusion tests. For general multi-factor models, useful sufficient conditions for martingality or strict local martingality include Novikov-type conditions, the Lyapunov conditions of Klebaner--Liptser \cite{KL14}, Criens \cite{Cri18} and Dandapani--Protter \cite{DP22}, the growth conditions of Lions--Musiela \cite{LM07}, and, on the orthant, the Fichera-drift condition of Fernholz--Karatzas \cite[Prop.~3]{FK10}. The purpose of this paper is to give necessary and sufficient coefficient conditions under stated boundary hypotheses and to obtain exact criteria for coupled non-affine models, including correlated factors with boundary singularities in $\theta$. We identify the classical ingredients and compare the scope of the resulting criteria with these exact and sufficient tests.

\subsection{The strategy}\label{sec:strategy}
Since $S$ is a nonnegative local martingale it is a supermartingale, and it is a true martingale on $[0,T]$ if and only if $\E[S_T]=1$. We test this equality using the exit probability of an auxiliary diffusion. Let
\begin{equation}\label{eq:Ltilde}
\Lt=\frac12\sum_{i,j=1}^d a_{ij}(y)\,\partial_{ij}+\sum_{i=1}^d\mut_i(y)\,\partial_i,\qquad \mut:=\mu+\sigma\theta,
\end{equation}
be the generator of the factor coordinate under the auxiliary law $\Qt_y$ on the same canonical space. The exit time $\tau_D$ is the same path functional under $\Prob_y$ and $\Qt_y$; its distribution depends on the law. Under $\Qt_y$, the factor drift is $\mut$ and the diffusion may leave $D$. Then, provided $Y$ itself never leaves $D$ under $\Prob_y$, Lemma~\ref{thm:main} shows that
\begin{equation}\label{eq:mass-intro}
\E_y[S_T]=\Qt_y(\tau_D>T).
\end{equation}
Hence, \emph{$S$ is a true martingale on $[0,T]$ if and only if the diffusion with Girsanov-shifted drift $\mu+\sigma\theta$ does not exit $D$ by time $T$ under $\Qt_y$.}

The assumption that the original diffusion stays in $D$ under $\Prob$ is essential here. If $S$ is a true martingale on $[0,T]$, the measure $d\Q^T=S_T\,d\Prob$ is equivalent to $\Prob$ on $\mathcal F_T$, so $Y$ still stays in $D$ up to $T$; Girsanov's theorem and uniqueness in law identify its law with that of the shifted diffusion. To see the converse and the role of the assumption, exhaust $D$ by relatively compact domains $D_n$ and let $\tau_n:=\inf\{t\ge0:Y_t\notin D_n\}$ on the canonical space. Stopped Girsanov gives
\[
\E_y[S_T\mathbf 1_{\{\tau_n>T\}}]=\Qt_y(\tau_n>T).
\]
Original non-exit means $\tau_n\uparrow\infty$ $\Prob$-a.s., so passage to the limit gives \eqref{eq:mass-intro}. The shifted diffusion can therefore exit even though the original one cannot. When $\Qt_y(\tau_D\le T)>0$, the stopped terminal densities $S_{T\wedge\tau_n}$ fail to be uniformly integrable, and $\E_y[S_T]<1$. If the original diffusion can reach the boundary, its behaviour there must be included in the model before applying such a criterion; see Remark~\ref{rem:reflect}.

Versions of formula \eqref{eq:mass-intro} appear in \cite{CFY05,CFR14,Ruf15}, and it is also the starting point of \cite{DP22}. By itself it does not answer the question, because the non-exit of a potentially degenerate multi-dimensional diffusion is not a condition on the coefficients. What this paper does is to turn it into one for different boundary geometries.

The $\Lt$-diffusion can leave $D$ in three ways: through a smooth boundary point, at infinity, or through a singular part of the boundary, such as a corner or an axis of the orthant. We treat the three separately, and at each we ask for a condition that is necessary as well as sufficient. Two principles organize the analysis.

First, at a boundary point $x$ of the $C^2$ part of $\partial D$, with outward normal $\nu(x)$, the sign of the \emph{Fichera drift}
\begin{equation}\label{eq:fichera-intro}
\Ft(x):=\Big\langle\mut(x)-\tfrac12\dv a(x),\,\nu(x)\Big\rangle,\qquad(\dv a)_i:=\sum_j\partial_ja_{ij},
\end{equation}
of $\Lt$ determines local boundary accessibility under the regularity and characteristic-boundary hypotheses of Section~\ref{sec:boundary}: positive values imply accessibility and negative values imply non-attainability, whereas the zero case requires additional information. This is Fichera's classification \cite{Fic56,OR73}, whose probabilistic content is due to Friedman \cite{Fri75} and Stroock--Varadhan \cite{SV72b}. Since $Y$ does not leave $D$ under $\Prob$, the boundary is characteristic, $\ip{a\nu}{\nu}=0$, and the Fichera drift $F$ of $L$ is $\le0$; the coefficient $\theta$ enters through
\begin{equation}\label{eq:shift-intro}
\Ft=F+\ip{\sigma\theta}{\nu},
\end{equation}
and the question is whether the Girsanov term $\ip{\sigma\theta}{\nu}$ can turn a repelling boundary point into an attracting one.

Second, \emph{the Girsanov drift can create an exit only where the coefficient $\theta$ is unbounded.} Near a boundary point where $\theta$ is bounded the measures $\Prob$ and $\Qt$ are locally equivalent, so the boundary behaviour under $\Qt$ is the one under $\Prob$, with no smoothness or sign condition whatsoever (Proposition~\ref{prop:inherit}); under the boundary regularity assumptions, $\Ft=F$ there (Corollary~\ref{cor:inherit}). The boundary analysis is therefore needed only where $\theta$ blows up, which in stochastic volatility models means the points where the volatility of the asset becomes infinite, and at infinity.

\subsection{Main results}\label{sec:mainresults}
We write $F$ and $\Ft$ for the Fichera drifts of $L$ and $\Lt$ on the $C^2$ part of $\partial D$, and, for $y\in D\setminus\{0\}$, we use the scalar functions
\begin{equation}\label{eq:AB-intro}
\mathcal A(y):=\ip{y}{a(y)y},\qquad\mathcal B(y):=\frac{\tr a(y)+2\ip{y}{\mut(y)}}{\ip{y}{a(y)y}},
\end{equation}
in terms of which $\Lt f(\tfrac12|y|^2)=\tfrac12\mathcal A\,[f''+\mathcal Bf']$. The main results are the following; the precise statements, with all side conditions, are the results quoted.
\begin{itemize}[leftmargin=2em]
\item \emph{Smooth boundary points} (Theorem~\ref{thm:fichera}). Let $\partial D\in C^2$, $a\in C^1(\bar D)$, $\mu,\mut\in C(\bar D)$, and suppose $\Qt$ does not explode at infinity, for instance because $\tr a+2\ip{y}{\mut}\le C(1+|y|^2)$. Then
\begin{equation}\label{eq:iff-intro}
S\ \text{is a true martingale}\iff \Ft\le0\ \text{on }\partial D\iff \ip{\sigma\theta}{\nu}\le-\Big\langle\mu-\tfrac12\dv a,\,\nu\Big\rangle\ \text{on }\partial D,
\end{equation}
provided that, at the characteristic points with vanishing Fichera drift $\{\Ft=0\}$ near which $\theta$ is unbounded, a modulus-of-continuity condition (T), stated in Definition~\ref{def:T}, holds. Condition (T) is automatic where $\Ft<0$, and at characteristic points with vanishing Fichera drift it holds when the coefficients are Lipschitz up to a $C^3$ boundary or, as for CIR-type factors, when $\ip{a\nabla\rho}{\nabla\rho}/\rho$ and $\Lt\rho$ are Lipschitz in the normal variable $\rho$ (Lemma~\ref{lem:T}). The implication ``$\Rightarrow$'' needs neither (T) nor the condition at infinity.
\item \emph{Infinity} (Propositions~\ref{cor:osgood} and \ref{thm:radial-iff}). Let $D$ be unbounded and suppose exits through its boundary $\partial D$ are excluded. Here $\partial D$ consists of finite boundary points. Explosion at infinity is treated separately. When $\mathcal A(y)$ and $\mathcal B(y)$ are bounded above or below by functions $A(\tfrac12|y|^2)$ and $B(\tfrac12|y|^2)$, Khasminskii's radial tests \cite{SV06} give sufficient conditions for non-explosion and explosion under $\Qt$. These include the Osgood-growth condition of Lions--Musiela \cite{LM07}. If, outside a ball, $\mathcal A(y)$ and $\mathcal B(y)$ depend only on the radius $|y|$ (the \emph{radial case}), so that $\mathcal A(y)=A(\tfrac12|y|^2)$ and $\mathcal B(y)=B(\tfrac12|y|^2)$, then $S$ is a true martingale iff the integral in Feller's test at $+\infty$ for the scalar diffusion with generator $\tfrac12A(r)[\tfrac{d^2}{dr^2}+B(r)\tfrac{d}{dr}]$ is infinite. Radiality here concerns these two scalar combinations of the coefficients; it does not require rotational symmetry of the full diffusion.
\item \emph{Orthant models} (Theorem~\ref{thm:SV}). On $D=(0,\infty)^d$, under the regularity and growth conditions above and a projection condition at intersections of boundary faces (boundary strata) $\{y_{i_1}=\dots=y_{i_k}=0\}$, $k\ge2$, $S$ is a true martingale iff $\Ft\le0$ on every face.
\end{itemize}

The boundary criterion combines classical local accessibility, the sufficient condition (T) at boundary points with zero Fichera drift, and preservation of boundary non-attainment where $\theta$ is bounded (Proposition~\ref{prop:inherit}). Under the boundary regularity hypotheses, $\Ft(x^0)>0$ at a smooth boundary point $x^0$ implies a positive probability of exit from $D$ through a smooth boundary patch on which $\Ft>0$, and hence rules out martingality on $[0,\infty)$. When $\Ft(x^0)\leq 0$, condition (T) provides a sufficient condition for non-attainability. When $\Ft(x_0)<0$, (T) is automatically satisfied. When $\Ft(x_0) =0$, the traditional sign test of \cite{Fri75} and \cite{SV72b} is inconclusive. Condition (T) fills this gap. A critical logarithmic family shows that condition (T) is not necessary (Example~\ref{ex:critical}). Theorem~\ref{thm:coupled-log} gives an exact threshold for coefficients with logarithmic corrections, modulated by correlated factors, including critical equality. Theorem~\ref{thm:correlated} gives a separate exact criterion for correlated square-root factors with boundary singularities in $\theta$. Their mixed covariance and shifted drift do not satisfy the boundary regularity of Theorem~\ref{thm:fichera}.

At infinity, Criens--Glau \cite[Cor.~5.2]{CG18} give exact martingality criteria under global radial assumptions on the shifted radial drift, radial variance, and quadratic-variation rate of the exponential. Proposition~\ref{thm:radial-iff} localizes the radial assumptions to the exterior of a ball, allows the quadratic-variation rate to vary along spheres, and combines the infinity test with conditions excluding finite-boundary exits. Proposition~\ref{prop:adapted-infinity} gives an exact test through more general functions with compact sublevel sets, and Corollary~\ref{cor:variable-clock} allows time-change rates that grow or decay with the scalar coordinate. Example~\ref{ex:anisotropic} applies this formulation to anisotropic dynamics. The sufficient conditions of Lions--Musiela \cite{LM07} and Dandapani--Protter \cite{DP22} are identified as the one-sided versions of the criterion at infinity, and Theorem~10 of \cite{DP22} becomes an equivalence in the radial case (Proposition~\ref{prop:LMDP}). At nonsmooth boundary pieces, Proposition~\ref{prop:sing} is a classical sufficient test; Proposition~\ref{thm:weighted-stratum} gives an exact criterion for a family with inverse-distance singularities and correlated modulation, including critical equality, and Example~\ref{ex:weighted-stratum} demonstrates the limitation of Euclidean projection.

The applications in Section~\ref{sec:SV} are what the criteria are designed for. The $3/2$ model \cite{CS07} and Grasselli's $4/2$ model \cite{Gra17} are recovered from the Fichera drift at the point where the volatility is infinite; for a two-factor inverse-CIR model and for a $4/2$ model modulated by a second factor $Y^2$ with modulating factor $g(Y^2)$ we obtain the closed-form conditions
\[
2\alpha_1-\eta_1^2\ \ge\ 2\eta_1\,\rho^-\sup g\qquad\text{and}\qquad 2\alpha_1-\eta_1^2\ \ge\ 2\eta_1\,(\rho b)^-\sup g ,
\]
For nonconstant $g$ and a nonzero Girsanov correction, the shifted first factor depends on the second and the shifted system is non-affine. Theorem~\ref{thm:correlated} extends these inequalities to any number of square-root factors with a positive definite correlation matrix and a bounded, locally H\"older modulation $g(Y^2,\ldots,Y^d)$. The threshold depends on the asset correlation with $Y^1$ and on $\sup g$, and includes equality; Example~\ref{ex:correlated} checks it when $\sup g$ is not attained. Section~\ref{sec:correlated} compares these coupled models with the scalar, affine and existing multifactor $4/2$ criteria. For affine diffusions the criteria reproduce the Riccati answer, and for Wishart processes the boundary condition $\Omega\Omega^{\top}\succeq(n+1)Q^{\top}Q$ appears as the sign of the Fichera drift on the cone, the curvature of the cone accounting for the $n+1$ (Section~\ref{sec:affine}). 
Proposition~\ref{prop:closed-state} and Theorem~\ref{thm:absorbed} treat closed state spaces and a singular density that vanishes at absorption. Section~\ref{sec:open} distinguishes these exact classes from the general problems that remain open.

For models satisfying the open-state-space Assumption~\ref{ass:standing}, the martingale criteria also yield uniqueness theorems for linear parabolic equations, as detailed in Appendix~\ref{sec:pde}. The loss of mass $1-\E_y[S_t]$ solves the Cauchy problem $\partial_tm=\Lt m$, $m(0,\cdot)=0$, with no boundary condition on $\partial D$ and none at infinity, so $S$ is a strict local martingale iff this problem has a nontrivial bounded solution (Remark~\ref{rem:pde}, Lemma~\ref{prop:CP}); and through the Doob transform with the harmonic weight $1+s$, $S$ is a true martingale iff the valuation equation of $(S,Y)$ has a unique solution in the class of linear growth in $s$ (Theorem~\ref{thm:VE}). The growth class therefore shifts the Fichera drift of the pricing operator on $\partial D$ by $\ip{\sigma\theta}{\nu}$, and the criteria above become explicit uniqueness classes (Proposition~\ref{thm:pde-unique}): for the $4/2$, inverse-CIR and Wishart valuation equations, uniqueness in the linear-growth class holds iff the parameter conditions above hold, with no boundary condition imposed where the volatility of the asset is infinite. In one dimension this is the mechanism of \cite{ET09,ET10,BX10,BKX12}; the multi-dimensional statements are new. 

\subsection{Relation to the literature}
For the underlying absolute-continuity theory, see Kabanov, Liptser and Shiryaev \cite{KLS80} and Liptser and Shiryaev \cite[Ch.~7]{LS01}. Criens--Glau \cite{CG18} develop this theory for possibly explosive semimartingales under a weaker uniqueness assumption than local uniqueness. For multidimensional diffusions with a conservative original law, Criens \cite[Cor.~2.10]{Cri21} states the equivalence between martingality and non-explosion of the drift-changed diffusion. The loss-of-mass formula \eqref{eq:mass-intro} appears explicitly in \cite[Thm.~2.4 and Rem.~2.6]{CFY05}; see also \cite{CFR14,Ruf15}. Khasminskii's criterion and its parabolic form are in \cite{Kha60,Aze74}. For Brownian motion on Riemannian manifolds, the equivalence between non-explosion (stochastic completeness) and uniqueness of bounded solutions of the heat equation is reviewed by Grigor'yan \cite[Thm.~6.2 and Cor.~6.3]{Gri99}. In one dimension the link between strict local martingales and non-uniqueness of the valuation equation in the linear-growth class is due to Ekstr\"om--Tysk \cite{ET09,ET10} and \cite{BX10,BKX12}.

Criens--Glau \cite[Sec.~5.2, Cor.~5.2]{CG18} obtain exact radial criteria for both martingality and uniform integrability; their partially radial comparison in \cite[Sec.~5.3, Prop.~5.5]{CG18} concerns global absolute continuity. Section~\ref{sec:LMDP} compares their radial hypotheses with those used here. Criens \cite[Thm.~2.2, Cor.~2.8]{Cri21} relates uniform integrability to almost-sure explosion of an auxiliary diffusion obtained by a time change, when the quadratic-variation rate is locally bounded away from zero, and derives Khasminskii-type tests in \cite[Sec.~3.1]{Cri21}. Uniform integrability over the infinite horizon is stronger than the martingale property on every finite horizon; Corollary~\ref{cor:coupled-log-ui} gives the exact distinction for the coupled logarithmic model.

The boundary classification of Section~\ref{sec:boundary} is Fichera's \cite{Fic56,OR73}; its probabilistic content, the attainability of the boundary points at which the noise is non-degenerate or the Fichera drift is positive, and the non-attainability of those at which it is negative (the sets $\Sig_2\cup\Sig_3$ and $\widetilde\Sig_1$ of Section~\ref{sec:boundary}), is due to Friedman \cite[Ch.~9 and 11]{Fri75} and Stroock--Varadhan \cite{SV72b}, and Lemma~\ref{thm:local} and Proposition~\ref{thm:attain} are versions of their results with a single Lyapunov function per case. At characteristic boundary points with zero Fichera drift, the sign test is inconclusive; condition (T) gives a sufficient condition for non-attainment. Related theories of boundary regularity for Kimura- and Heston-type operators are developed by Epstein--Mazzeo \cite{EM13} and Feehan--Pop \cite{FP15}. 

The radial test, the comparison theorem and the projection argument near intersections of boundary faces are due to Khasminskii \cite{Kha60,Kha12}, Ikeda--Watanabe \cite{IW89} and Ramasubramanian \cite{Ram88}. In stochastic portfolio theory, on the orthant, the Fichera-drift criterion for the deflator of a market is Proposition~3 of Fernholz--Karatzas \cite{FK10} (see also \cite{FK09}), and their Theorem~1 characterizes the optimal arbitrage as the smallest nonnegative supersolution of a Cauchy problem, which is the non-uniqueness phenomenon of Remark~\ref{rem:pde} in that setting. The comparison with their Proposition~3 is made precisely in Remark~\ref{rem:FK}; the local boundary mechanism is classical. Theorem~\ref{thm:fichera} treats general coefficients $\theta$ in the stochastic exponential. Under its regularity hypotheses, $\Ft(x^0)>0$ at even one smooth boundary point implies a positive probability of exit through a nearby boundary patch, and therefore rules out martingality on $[0,\infty)$. By continuity, positivity at $x^0$ persists on that patch; the argument does not require hitting the individual point $x^0$. 

For affine models the exact Riccati-based characterizations of Kallsen--Muhle-Karbe \cite{KMK10}, Mayerhofer--Muhle-Karbe--Smirnov \cite{MMS11} and Keller-Ressel--Mayerhofer \cite{KRM15} also give the exponential moments and their explosion times \cite{AP07,KR11}, which are finiteness questions for exponential functionals rather than martingale questions and are not addressed here.

\subsection{Organization}
Section~\ref{sec:exact} contains the reduction to a non-exit problem, together with the two tools by which non-exit is decided: Khasminskii's criterion and the Lyapunov lemma. Sections~\ref{sec:boundary}, \ref{sec:infty} and \ref{sec:sing} treat the three kinds of exit: smooth boundary points, infinity (including the comparison with \cite{CG18,Cri21,LM07,DP22}), and singular strata. Section~\ref{sec:SV} collects the criteria in the form in which they apply to stochastic volatility models and works out the examples, Section~\ref{sec:affine} checks them against the affine theory, and Section~\ref{sec:open} concludes. The closed-state-space formulation is in Section~\ref{sec:extended-state}; the logarithmic and weighted-distance criteria are in Sections~\ref{sec:log-extension} and \ref{sec:weighted-stratum}, and the criterion using a function with compact sublevel sets is in Section~\ref{sec:adapted-infinity}. Appendix~\ref{app:reg} contains the regularity lemma on which the finite-horizon and PDE statements rest. Appendix~\ref{sec:pde} reports uniqueness results for linear parabolic equations. Appendix~\ref{app:prelim-proofs} contains the proofs of Section~\ref{sec:exact}. Figure~\ref{fig:hierarchy} in Appendix~\ref{app:hierarchy} shows how the statements fit together.

\section{The reduction to a non-exit problem}\label{sec:exact}

This section fixes the canonical formulation and records the loss-of-mass identity, its finite-horizon consequences, and the Khasminskii and Lyapunov criteria used below. It also specifies how the identity changes when the state space includes an accessible boundary. Proofs are collected in Appendix~\ref{app:prelim-proofs}.

\subsection{Canonical formulation and the loss-of-mass formula}\label{sec:canonical}

\begin{assumption}\label{ass:standing}\leavevmode
\begin{enumerate}[label=(\roman*),leftmargin=2em]
\item $D\subseteq\R^d$ is a domain, exhausted by bounded domains $D_n\subset\subset D_{n+1}$, $\bigcup_n D_n=D$.
\item $a=\sigma\sigma^{\top}$ is continuous on $D$ and $\ip{a(y)v}{v}>0$ for $y\in D$, $v\ne0$; $\mu$, $\sigma$, $\theta$ are locally bounded and measurable on $D$.
\item \emph{($\Prob$-non-exit)} The generalized martingale problem for $L=\frac12\sum_{i,j}a_{ij}\partial_{ij}+\sum_i\mu_i\partial_i$ on $D$ is well posed and conservative.
\end{enumerate}
\end{assumption}
Let $(\Omega,\mathcal F,(\mathcal F_t))$ be the canonical space of maximal paths of $(Y,W)$ in $D\times\R^m$, continuous up to their lifetime and sent to a cemetery state thereafter, with the coordinate filtration before completion under any law. The lifetime is the limit of exit times from a relatively compact exhaustion of $D\times\R^m$. We take $Y_0=y$ and $W_0=0$. Under $\Prob_y$, the coordinate pair solves the generalized martingale problem with generator
\begin{equation}\label{eq:canonical-generator}
\mathcal G f
=\frac12\sum_{i,j=1}^d a_{ij}(y)\partial_{y_i y_j}f
+\sum_{i=1}^d\mu_i(y)\partial_{y_i}f
+\frac12\Delta_w f
+\sum_{i=1}^d\sum_{\alpha=1}^m\sigma_{i\alpha}(y)\partial_{y_i w_\alpha}f.
\end{equation}
Here the test functions are in $C_c^2(D\times\R^m)$ and the martingale problem is understood up to the lifetime. Its local characteristics give $d\langle W^\alpha,W^\gamma\rangle_t=\delta_{\alpha\gamma}dt$ and $d\langle Y^i,W^\alpha\rangle_t=\sigma_{i\alpha}(Y_t)dt$. This martingale problem is associated with $(Y,W)$ where $Y$ follows \eqref{eq:Y}. 

The stochastic exponential \eqref{eq:S} is defined from these coordinates under $\Prob_y$. The resulting pair $(S,Y)$, with $S_0=1$, solves the martingale problem for the joint generator $\mathcal L$ in \eqref{eq:pricing-op}. Including $W$ retains the Brownian components of the integrand that may be orthogonal to the factor noise when $m>d$.

Choose the exhaustion in Assumption~\ref{ass:standing}(i) with $\partial D_n\in C^2$ and define, as coordinate path functionals,
\begin{equation}\label{eq:canonical-exit}
\tau_n:=\inf\{t\ge0:Y_t\notin D_n\},\qquad
\tau_D:=\lim_{n\to\infty}\tau_n.
\end{equation}
Under each of the two laws considered here, before each $\tau_n$, the coordinate $W$ has unit covariance and a bounded drift while $Y$ remains in a compact subset of $D$. Thus $W$ cannot explode before $\tau_D$, and the joint lifetime is $\tau_D$. Assumption~\ref{ass:standing}(iii) says $\Prob_y(\tau_D=\infty)=1$.

The auxiliary law $\Qt_y$ is a solution on the same canonical space of the generalized martingale problem with generator
\begin{equation}\label{eq:canonical-shift}
\widetilde{\mathcal G}f
=\mathcal G f+(\sigma\theta)(y)\cdot\nabla_yf+\theta(y)\cdot\nabla_wf.
\end{equation}
Thus $\widehat W_t:=W_t-\int_0^t\theta(Y_s)ds$ is Brownian up to $\tau_D$ under $\Qt_y$, and the factor coordinate has generator $\Lt$. Under Assumption~\ref{ass:standing}(ii), the generalized martingale problem for $\Lt$ is well posed \cite[Thm.~1.13.1]{Pin95}; a weak realization with diffusion coefficient $\sigma$ gives the stated joint law, and all probabilities involving only $Y$ are determined by this unique factor law. 

\begin{lemma}\label{thm:main}
Under Assumption~\ref{ass:standing}, for every $T>0$ and $y\in D$,
\begin{equation}\label{eq:mass}
\E_y[S_T]=\Qt_y(\tau_D>T).
\end{equation}
Consequently:
\begin{enumerate}[label=(\roman*),leftmargin=2em]
\item $S$ is a true martingale on $[0,T]$ iff $\Qt_y(\tau_D\le T)=0$;
\item $S$ is a true martingale on $[0,\infty)$ iff $\Qt_y(\tau_D<\infty)=0$; by \cite[Thm.~1.15.1]{Pin95} this holds for one $y\in D$ iff it holds for all $y\in D$.
\end{enumerate}
\end{lemma}

The martingality--non-explosion equivalence is part of the general absolute-continuity theory of \cite{CG18}; see also \cite[Cor.~2.10]{Cri21} for the conservative diffusion setting. Formula \eqref{eq:mass} is the ``loss of mass'' formula of Cheridito--Filipovi\'c--Yor \cite[Thm.~2.4, Rem.~2.6]{CFY05} (there for $\theta=\sigma^{\top}\varphi$), of Carr--Fisher--Ruf \cite{CFR14} and Ruf \cite{Ruf15} (general $\theta$), of \cite{DP22} in the vector models, and, in one dimension, of Delbaen--Shirakawa \cite{DS02}, Mijatovi\'c--Urusov \cite{MU12} and \cite{BKX12}.

\begin{remark}[reflecting factors]\label{rem:reflect}
Assumption~\ref{ass:standing}(iii) requires the original process to stay in the open set $D$. A reflected factor, such as a CIR factor that can reach zero and return, instead lives on a state space containing that boundary. Proposition~\ref{prop:closed-state} proves the loss-of-mass identity on a closed state space when $\theta$ is bounded on its compact subsets and both drift specifications are well posed with the same boundary rule. Absorbed factors require the density itself to be specified at absorption: it may retain positive mass, or a singular coefficient $\theta$ may force it to zero. Section~\ref{sec:extended-state} gives the resulting identities, and Theorem~\ref{thm:absorbed} gives an exact criterion for an absorbed factor coupled to correlated factors. For the one-dimensional boundary classification see \cite{MU12,MNU10,DLR21}.
\end{remark}

\subsection{Finite horizons}
For time-homogeneous coefficients the finite-horizon and infinite-horizon statements are almost equivalent, so that Lemma~\ref{thm:main}(ii), the form in which the criteria below are stated, loses nothing.

\begin{proposition}\label{prop:finite}
Under Assumption~\ref{ass:standing}:
\begin{enumerate}[label=(\roman*),leftmargin=2em]
\item If $\E_y[S_T]=1$ for some $T>0$ and all $y\in D$, then $S$ is a true martingale on $[0,\infty)$ for every starting point.
\item Assume moreover that the coefficients $a,\mut$ are locally H\"older continuous on $D$. Then either $S$ is a true martingale on $[0,\infty)$ for every starting point, or $S$ is a strict local martingale on $[0,T]$ for every $T>0$ and every starting point.
\end{enumerate}
\end{proposition}

The regularity of the loss of mass used in part (ii) is proved by localization on the auxiliary domains $D_n$ in Appendix~\ref{app:reg}.

\subsection{Khasminskii's criterion and Lyapunov functions}
Lemma~\ref{thm:main} reduces the martingale property to non-exit under $\Qt$. Khasminskii's criterion \cite{Kha60} expresses non-exit as the absence of bounded positive solutions of an elliptic equation.

\begin{lemma}\label{thm:khas}
Let Assumption~\ref{ass:standing} hold and let $a,\mut$ be locally H\"older continuous on $D$. Then $S$ is a true martingale on $[0,\infty)$ if and only if, for some (equivalently, every) $\lambda>0$, the equation $\Lt u=\lambda u$ has no bounded positive solution in $D$.
\end{lemma}

\begin{remark}[PDE formulation]\label{rem:pde}
The loss of mass $m(t,y):=1-\E_y[S_t]=\Qt_y(\tau_D\le t)$ is a bounded nonnegative solution of the Cauchy problem $\partial_tm=\Lt m$ on $(0,\infty)\times D$ with $m(0,\cdot)=0$ (classically under local H\"older continuity, Lemma~\ref{lem:regularity}). Hence $S$ is a strict local martingale if and only if this Cauchy problem has a nontrivial bounded solution, i.e.\ iff uniqueness fails in the bounded class; and, through the Doob transform with the harmonic weight $1+s$, iff uniqueness fails for the valuation equation of $(S,Y)$ in the class of linear growth in $s$. Every criterion of Sections~\ref{sec:boundary}--\ref{sec:sing} is therefore also a uniqueness criterion for a linear parabolic equation. In one dimension this is the mechanism of \cite{ET09,ET10,BX10,BKX12}; on the whole space it is the stochastic completeness theory of \cite{Aze74,Gri99}. The precise statements are collected in Appendix~\ref{sec:pde}.
\end{remark}

Lemma~\ref{thm:khas} is sharp but not ``on the coefficients''. Its two practical one-sided versions are the Lyapunov criteria of Khasminskii \cite{Kha60,Kha12}, see \cite[Thm.~6.7.1]{Pin95}, \cite[Thms.~10.2.1--10.2.2]{SV06}; they contain Theorems~3 and~4 of \cite{DP22} and underlie the deterministic criteria of \cite{KL14,Cri18,Li17}. The forms used below and their proofs in Appendix~\ref{app:prelim-proofs} specify the boundary and positivity conditions needed for the applications. In the following result, $B\subset D$ is relatively closed in $D$ (typically $B=\bar B_0\cap D$ for a ball $B_0$), and $\partial_DB$ denotes its boundary relative to $D$.

\begin{lemma}\label{thm:lyap}
Let Assumption~\ref{ass:standing}(ii) hold for $\Lt$ on $D$.
\begin{enumerate}[label=(\roman*),leftmargin=2em]
\item (\emph{non-exit}) Suppose there are $\lambda\ge0$ and $V\in C^2(D)$, $V\ge0$, with $\Lt V\le\lambda V$ on $D$ and $V(y)\to\infty$ as $y\to\partial D\cup\{\infty\}$, i.e.\ $\{V\le k\}\subset\subset D$ for every $k$. Then $\Qt_y(\tau_D<\infty)=0$ for all $y\in D$. The same conclusion holds if $\Lt V\le\lambda V$ only on $D\setminus B$ with $B\subset\subset D$.
\item (\emph{exit}) Suppose there are $\lambda>0$, a set $B\subset D$ relatively closed in $D$ (possibly empty), $y_0\in D\setminus B$ and a bounded $V\in C^2(D\setminus B)\cap C\big((D\setminus B)\cup\partial_DB\big)$ with $\Lt V\ge\lambda V$ on $D\setminus B$ and $V(y_0)>\max\{0,\sup_{\partial_DB}V\}$ (with $\sup\emptyset=-\infty$). Then $\Qt_{y_0}(\tau_D<\infty)>0$, hence $\Qt_y(\tau_D<\infty)>0$ for all $y\in D$.
\end{enumerate}
\end{lemma}

The coefficient tests of Sections~\ref{sec:boundary}--\ref{sec:sing} use specific Lyapunov functions: a power or a logarithm of the distance to the boundary in Section~\ref{sec:boundary}, a function of $\tfrac12|y|^2$ in Section~\ref{sec:infty}, and a function of the distance to a stratum in Section~\ref{sec:sing}. The necessity statements come from part (ii) and the exit lemma of the next section, the sufficiency statements from part (i).

\subsection{Boundary behaviour on closed state spaces}\label{sec:extended-state}

Boundary behaviour is part of the specification of the model. A reflected factor remains in its state space when it touches the boundary; an absorbed factor may leave a positive terminal density there, or may drive a singular density to zero. These cases give different martingale criteria. The one-dimensional theory treats absorption in \cite{MU12,MNU10} and reflection in \cite{DLR21}. On a closed state space, localization must preserve the prescribed boundary rule.

\begin{proposition}\label{prop:closed-state}
Let $E\subset\R^d$ be closed and let $Y$ be a conservative continuous solution in $E$ of
\[
dY_t=\mu(Y_t)\,dt+\sigma(Y_t)\,dW_t+dK_t,
\]
where $K$ is a continuous process of finite variation specified by a reflection rule at $\partial E$. Suppose the systems with drifts $\mu$ and $\mu+\sigma\theta$, the same diffusion coefficient and the same boundary rule are well posed in law up to their lifetimes. Assume that $\theta$ is bounded on every compact subset of $E$, let $\tau_n = \inf\{t\geq 0\,:\, Y_t \notin E \text{ or } |Y_t| \geq n\}$, and $\zeta = \lim_{n\rightarrow \infty}\tau_n$.
Then, for $S=\cE(\int\theta(Y)\cdot dW)$ and every $T>0$,
\begin{equation}\label{eq:closed-mass}
\E_y[S_T]=\Qt_y(\zeta>T).
\end{equation}
Consequently $S$ is a true martingale on every finite horizon if and only if the shifted system is conservative.
\end{proposition}

For example, let $X$ be Brownian motion reflected at $0$ and take $S=\cE(\int\lambda X^\gamma\,dW)$, with $\gamma>0$ and $\lambda\in\R$. The shifted reflected diffusion has drift $\lambda x^\gamma$. Feller's test at infinity gives
\[
S\text{ is a true martingale on every finite horizon}
\quad\Longleftrightarrow\quad \lambda\le0\ \text{or}\ \gamma\le1.
\]
Hitting $0$ does not remove mass in this example. Singularities of $\theta$ at an attainable boundary require a specification of the density at absorption; see Section~\ref{sec:absorbed-model}.

There is a different identity when the density has a positive limit at absorption and is then frozen. Writing $\tau$ for that absorption time, the decomposition is
\begin{equation}\label{eq:retained-mass}
\E_y[S_T]=\E_y[S_T;\,T<\tau]
+\E_y[S_\tau;\,\tau\le T].
\end{equation}
When stopped Girsanov identifies the first term with the survival probability of the killed shifted system, the second term is the mass retained at the boundary. It is zero in Theorem~\ref{thm:absorbed}; for $\theta\equiv0$ the two terms are $\Prob_y(T<\tau)$ and $\Prob_y(\tau\le T)$. Thus absorption alone does not determine which boundary events cause loss of mass. A general classification with both types of boundary remains beyond this result.

\section{Smooth boundary points: the Fichera criterion}\label{sec:boundary}

This section treats the first kind of exit, through the smooth part of the boundary, and proves the first main result, Theorem~\ref{thm:fichera}: on the $C^2$ part of $\partial D$ the sign of the Fichera drift of $\Lt$ decides the martingale property, and it does so through the value of the Girsanov term $\ip{\sigma\theta}{\nu}$ at the boundary points where $\theta$ is unbounded only.

\subsection{The Fichera drift and the statement of the theorem}
Let $\partial D^{(2)}$ denote the $C^2$ part of $\partial D$ and, for $x\in\partial D^{(2)}$, let $\nu(x)$ be the outward unit normal and $\rho(x):=\dist(x,\partial D)$, which is $C^2$ in $D$ near $\partial D^{(2)}$ with $\nabla\rho=-\nu$ on $\partial D^{(2)}$. For a generator $L_b=\frac12\sum a_{ij}\partial_{ij}+\sum b_i\partial_i$ whose coefficients satisfy, near a point $x^0\in\partial D^{(2)}$,
\begin{equation}\label{eq:bdryreg}
a\in C^1(\bar D\cap N),\qquad b\in C(\bar D\cap N)\qquad\text{for a neighbourhood $N$ of $x^0$},
\end{equation}
define the \emph{Fichera drift} (Fichera \cite{Fic56}, Ole\u{\i}nik--Radkevi\v{c} \cite{OR73}, Friedman \cite[Ch.~9, 11]{Fri75})
\begin{equation}\label{eq:fichera}
F_b(x):=\Big\langle b(x)-\tfrac12\dv a(x),\,\nu(x)\Big\rangle,\qquad x\in\partial D^{(2)}\cap N ,
\end{equation}
and split $\partial D^{(2)}$ into
\begin{align*}
\Sig_3&:=\{\ip{a\nu}{\nu}>0\},& \Sig_2&:=\{\ip{a\nu}{\nu}=0,\ F_b>0\},\\
\widetilde\Sig_1&:=\{\ip{a\nu}{\nu}=0,\ F_b<0\},& \Sig_1^0&:=\{\ip{a\nu}{\nu}=0,\ F_b=0\}.
\end{align*}
On $\Sig_3$ the boundary is non-characteristic and is reached by the noise. On the characteristic part the Fichera drift decides whether the boundary attracts ($\Sig_2$, drift pointing outward) or repels ($\widetilde\Sig_1$). On $\Sig_1^0$, the boundary is characteristic and $F_b=0$. Thus the corrected drift $b-\tfrac12\dv a$ has zero normal component there, although $b$ itself need not be tangent to the boundary. On this set the sign of $F_b$ alone does not decide attainability; the sufficient modulus-of-continuity condition below covers both the Lipschitz case of \cite{Fri75} and the CIR-type square-root degeneracy.

\begin{definition}[condition (T)]\label{def:T}
$L_b$ satisfies (T) at $x^0\in\partial D^{(2)}$ if there exist an open neighbourhood $N$ of $x^0$, as in \eqref{eq:bdryreg} and small enough that $\rho\in C^2(D\cap N)$, and a constant $C\ge0$ such that $\ip{a\nu}{\nu}=0$ on $\partial D\cap N$ and
\begin{equation}\label{eq:T}
\frac{\ip{a\nabla\rho}{\nabla\rho}}{\rho}-2L_b\rho\ \le\ C\rho\qquad\text{in }D\cap N.
\end{equation}
This is a local condition: once it holds on $N$, it holds with the same $C$ on every smaller neighbourhood of $x^0$.
\end{definition}

Condition (T) is a local criterion for non-attainment of the boundary. The requirement $\ip{a\nu}{\nu}=0$ means that the noise has no normal component at the boundary. In the interior, put $h:=\ip{a\nabla\rho}{\nabla\rho}$ and $\beta:=L_b\rho$. By It\^o's formula, $\beta(X_t)$ is the instantaneous drift of the distance process $\rho(X_t)$ and $h(X_t)$ is its instantaneous variance rate. The inequality in (T) is equivalent to
\[
\beta\ \ge\ \frac{h}{2\rho}-\frac C2\rho,
\qquad\text{or, equivalently,}\qquad
L_b(-\log\rho)=\frac{h}{2\rho^2}-\frac{\beta}{\rho}\ \le\ \frac C2.
\]
Thus the drift that increases the distance compensates for the fluctuation term $h/(2\rho)$, up to an error of order $\rho$. The function $-\log\rho$ tends to $+\infty$ at the boundary but has drift bounded above under (T); after localizing to a sufficiently small neighbourhood, the stopped It\^o argument rules out reaching that boundary patch in finite time (Lemma~\ref{thm:local}(ii)). Lemma~\ref{lem:T} shows that $h/\rho-2\beta$ tends to $2F_b$ at the boundary. Consequently (T) implies $F_b\le0$ on the patch and controls the positive part of $h/\rho-2\beta$ by $C\rho$ in the interior.

For a CIR factor $dY=(\alpha-\kappa Y)dt+\eta\sqrt Y\,dZ$ near $0$, one has $\rho(y)=y$, $h(y)=\eta^2y$ and $\beta(y)=\alpha-\kappa y$, so (T) reads $\eta^2-2\alpha+2\kappa y\le Cy$. It holds exactly when $2\alpha\ge\eta^2$, with $C=2\max\{\kappa,0\}$: this is Feller's non-attainment condition, including the critical equality. Thus (T) covers square-root diffusion coefficients, which are not Lipschitz at $0$. In general it is sufficient, but not necessary, for non-attainment (Example~\ref{ex:critical}).

Throughout the paper $F:=F_\mu$ and $\Ft:=F_{\mut}$ denote the Fichera drifts of $L$ and of $\Lt$, so that $\Ft=F+\ip{\sigma\theta}{\nu}$ wherever both are defined. The main result of this section reads as follows; its hypotheses are stated globally for simplicity, and the local statements from which it is assembled are Propositions~\ref{thm:attain} and \ref{prop:inherit} below.

\begin{theorem}\label{thm:fichera}
Let Assumption~\ref{ass:standing} hold and assume $\partial D\in C^2$, $a\in C^1(\bar D)$ and $\mu,\mut\in C(\bar D)$. Then:
\begin{enumerate}[label=(\roman*),leftmargin=2em]
\item $\ip{a\nu}{\nu}=0$ and $F\le0$ on $\partial D$.
\item \emph{(necessity)} If $S$ is a true martingale on $[0,\infty)$ from some (equivalently, every) $y\in D$, or if for some $T>0$ it is a true martingale on $[0,T]$ from every $y\in D$, then $\Ft\le0$ on $\partial D$. If the coefficients are locally H\"older continuous on $D$, it suffices that $S$ be a true martingale on some $[0,T]$ from some $y\in D$ (Proposition~\ref{prop:finite}(ii)).
\item \emph{(sufficiency)} Suppose
\begin{equation}\label{eq:lingrowth}
\tr a(y)+2\ip{y}{\mut(y)}\ \le\ C\,(1+|y|^2)\qquad(y\in D),
\end{equation}
or, more generally, $\Qt_y\big(\tau_D<\infty,\ \sup_{t<\tau_D}|Y_t|=\infty\big)=0$ for all $y\in D$ (e.g.\ under the radial conditions of Lemma~\ref{thm:radial-sandwich}(i)). If $\Ft\le0$ on $\partial D$ and $\Lt$ satisfies (T) at every point of $\{\Ft=0\}$ near which $\theta$ is unbounded (in particular if $\Ft<0$ on $\partial D$, or if $\theta$ is locally bounded on $\bar D$, or if $\partial D\in C^3$ and $\sigma,\mut$ are Lipschitz on $\bar D$), then $S$ is a true martingale on $[0,\infty)$ for every $y\in D$.
\end{enumerate}
Consequently, suppose the non-explosion condition in (iii) holds and either $\{\Ft=0\}=\emptyset$ or (T) holds at the points of $\{\Ft=0\}$ near which $\theta$ is unbounded. Then
\begin{equation}\label{eq:iff}
\begin{aligned}
S\ \text{is a true martingale}&\iff \Big\langle \mut-\tfrac12\dv a,\,\nu\Big\rangle\le0\ \text{ on }\partial D\\
&\iff \ip{\sigma\theta}{\nu}\le -\Big\langle \mu-\tfrac12\dv a,\,\nu\Big\rangle\ \text{ on }\partial D .
\end{aligned}
\end{equation}
\end{theorem}

The proof occupies the rest of the section and has three steps. The first (Section~\ref{sec:local}) is local: with one Lyapunov function per case we show that the diffusion could leave $D$ from the points of $\Sig_2\cup\Sig_3$ and cannot leave $D$ from the points of $\widetilde\Sig_1$, nor from the points of $\Sig_1^0$ at which (T) holds. The second (Section~\ref{sec:global}) passes from local to global by a covering argument. The third (Section~\ref{sec:inherit}) shows that, where $\theta$ is bounded near a boundary portion, non-attainment under $\Prob$ implies non-attainment under $\Qt$. The proof of the theorem is completed in Section~\ref{sec:proof-fichera}, and Section~\ref{sec:tangential} discusses the characteristic boundary where the Fichera drift vanishes but (T) fails, on which the theorem only provides a necessary condition.

\subsection{Local analysis}\label{sec:local}
We begin with an identity that expresses the Fichera drift through the distance function. It is behind both Lyapunov functions used below, it explains where condition (T) comes from, and it contains a curvature term, $\tr(a\nabla^2\rho)$, which cannot be neglected on curved boundaries (it produces the correct Wishart threshold in Proposition~\ref{prop:wishart}).

\begin{lemma}\label{lem:T}
Assume \eqref{eq:bdryreg} and $\ip{a\nu}{\nu}=0$ on $\partial D\cap N$. Then, as $D\ni x\to\bar x\in\partial D\cap N$,
\begin{equation}\label{eq:limitT}
\frac{\ip{a\nabla\rho}{\nabla\rho}}{\rho}(x)-2L_b\rho(x)\ \longrightarrow\ 2F_b(\bar x).
\end{equation}
In particular (T) holds at every point of $\widetilde\Sig_1\cap N$. It also holds at a point $x^0\in\Sig_1^0$ provided $F_b\le0$ on $\partial D\cap N'$ for some neighbourhood $N'\subset N$ of $x^0$ and either $\partial D\in C^3$ near $x^0$ and $\sigma$, $b$ are Lipschitz up to the boundary, or $\ip{a\nabla\rho}{\nabla\rho}/\rho$ and $L_b\rho$ are Lipschitz in the normal variable, uniformly near $x^0$. The sign condition on the neighbouring boundary patch cannot be dropped (Example~\ref{ex:isolated}).
\end{lemma}

\begin{proof}
Put $n:=\nabla\rho$ (inward normal, extended as a $C^1$ field near $\partial D$) and $h:=\ip{an}{n}$. Since $a\ge0$ and $h=0$ on $\partial D\cap N$, we have $an=0$ there, so $\partial_n h=\sum_{ijk}n_in_jn_k\partial_ka_{ij}$ on the boundary. Differentiating $\sum_ia_{ij}n_i=0$ along tangential directions and using $an=0$ once more gives $\ip{\dv a}{n}=\partial_nh-\tr(a\nabla^2\rho)$ on $\partial D\cap N$. Hence, on the boundary,
\[
\partial_nh-2L_b\rho=\partial_nh-\tr(a\nabla^2\rho)-2\ip{b}{n}=\ip{\dv a}{n}-2\ip{b}{n}=2F_b ,
\]
because $\nu=-n$. Finally $h(x)/\rho(x)=\partial_nh(\xi)$ for some $\xi$ on the normal segment from $\bar x$ to $x$ (mean value theorem, $h(\bar x)=0$), and $\partial_nh$, $L_b\rho$ are continuous up to the boundary by \eqref{eq:bdryreg}; this gives \eqref{eq:limitT}. If $F_b(x^0)<0$ the left side of \eqref{eq:T} is negative near $x^0$. In the Lipschitz case, $\sigma^{\top}n$ vanishes on the boundary and is Lipschitz, so $h\le C\rho^2$; moreover $L_b\rho(x)\ge L_b\rho(\bar x)-C\rho$ (this uses $\nabla^2\rho$ Lipschitz, i.e.\ $\partial D\in C^3$), where $\bar x$ is the nearest boundary point of $x$, and $-2L_b\rho(\bar x)=2F_b(\bar x)-\partial_nh(\bar x)\le0$ because $F_b(\bar x)\le0$ by the sign hypothesis on the patch and $\partial_nh(\bar x)\ge0$ ($h\ge0$ vanishes on the boundary). Combining, $h/\rho-2L_b\rho\le C\rho+0+2C\rho$. In the last case, $\varphi:=h/\rho-2L_b\rho$ satisfies $\varphi(x)\le\lim_{x'\to\bar x}\varphi(x')+C\rho(x)=2F_b(\bar x)+C\rho(x)\le C\rho(x)$ by \eqref{eq:limitT} and $F_b(\bar x)\le0$. Both arguments use the sign of $F_b$ at the nearest boundary point, not only at $x^0$.
\end{proof}

For the rest of this section $\Prob^b_x$ denotes the solution of the generalized martingale problem for $L_b$ on $D$ (Assumption~\ref{ass:standing}(ii) for $a$ and $b$), $X$ the coordinate process and $\tau_D$ its exit time. For an open $W\subset D$ we write $\tau_W$ for the exit time from $W$; whenever the coefficients are bounded on $W$, the limit $X_{\tau_W-}:=\lim_{t\uparrow\tau_W}X_t$ exists on $\{\tau_W<\infty\}$ (the drift and the stochastic integral in the representation of $X$ extend continuously to $\tau_W$). The next lemma is the mechanism by which a Lyapunov function forces an exit; it is a probabilistic barrier argument.

\begin{lemma}\label{lem:exit}
Let $x^0\in\partial D$, let $W=D\cap B(x^0,\varepsilon_0)$, and suppose that $a$ and $b$ are bounded on $W$ and that there are $\gamma>0$ and $U\in C^2(W)\cap C(\bar W)$ with $U(x^0)=0$, $U>0$ on $\bar W\setminus\{x^0\}$ and $L_bU\le-\gamma$ in $W$. Then, for $x\in W$,
\[
\E^b_x[\tau_W]\le U(x)/\gamma ,\qquad\text{and}\qquad \lim_{W\ni x\to x^0}\Prob^b_x\big(\tau_W<\infty,\ |X_{\tau_W-}-x^0|<\delta\big)=1\quad\text{for every }\delta>0 .
\]
In particular $\Prob^b_x(\tau_D<\infty)>0$ for $x\in W$ close to $x^0$, hence for all $x\in D$.
\end{lemma}

\begin{proof}
Let $W_n:=W\cap D_n$ and $\tau^n:=\tau_{W_n}$. By It\^o's formula and optional stopping (the coefficients and $U$ are bounded on $W_n$), $\E^b_x[U(X_{t\wedge\tau^n})]-U(x)=\E^b_x\big[\int_0^{t\wedge\tau^n}L_bU(X_s)ds\big]\le-\gamma\E^b_x[t\wedge\tau^n]$, so $\gamma\E^b_x[t\wedge\tau^n]\le U(x)-\E^b_x[U(X_{t\wedge\tau^n})]\le U(x)$; letting $n\to\infty$ and $t\to\infty$ gives the first claim. For the second, fix $\delta>0$, $t>0$ and $x\in W$ with $|x-x^0|<\delta/2$. Since $a,b$ are bounded by some $K$ on $W$, $\E^b_x[\sup_{s\le t\wedge\tau_W}|X_s-x|^2]\le c(t):=2K^2t^2+8dKt$ (Doob's inequality), so
\[
\Prob^b_x\big(\tau_W<\infty,\ |X_{\tau_W-}-x^0|<\delta\big)\ \ge\ \Prob^b_x\Big(\tau_W\le t,\ \sup_{s\le t\wedge\tau_W}|X_s-x|<\delta/2\Big)\ \ge\ 1-\frac{U(x)}{\gamma t}-\frac{4c(t)}{\delta^2}.
\]
Let $x\to x^0$ and then $t\downarrow0$. Finally, if $\delta<\varepsilon_0$, the points of $\partial W$ within distance $\delta$ of $x^0$ lie on $\partial D$, so $\{\tau_W<\infty,|X_{\tau_W-}-x^0|<\delta\}\subset\{\tau_D<\infty\}$; the extension to all $x\in D$ is \cite[Thm.~1.15.1]{Pin95}.
\end{proof}

The next lemma supplies the local boundary arguments for Theorem~\ref{thm:fichera}. Under its hypotheses, part (i) gives positive probability of exit near a non-characteristic boundary point or a characteristic point with positive Fichera drift. Applied with $b=\mut$, this rules out martingality on $[0,\infty)$ by Lemma~\ref{thm:main}. Part (ii) excludes exit through a sufficiently small characteristic boundary patch near a point where the Fichera drift is negative, or vanishes and (T) holds. The covering argument of Proposition~\ref{thm:attain} combines these local exclusions; together with control of the remaining exits, they yield the sufficiency part of Theorem~\ref{thm:fichera}.

\begin{lemma}\label{thm:local}
Let Assumption~\ref{ass:standing}(ii) and \eqref{eq:bdryreg} hold at $x^0\in\partial D^{(2)}$.
\begin{enumerate}[label=(\roman*),leftmargin=2em]
\item If $x^0\in\Sig_3$, or if $x^0\in\Sig_2$ and $\ip{a\nu}{\nu}=0$ on $\partial D\cap N$, then $\lim_{D\ni x\to x^0}\Prob^b_x(\tau_D<\infty,\ |X_{\tau_D-}-x^0|<\delta)=1$ for every $\delta>0$; in particular $\Prob^b_x(\tau_D<\infty)>0$ for every $x\in D$.
\item Suppose $\ip{a\nu}{\nu}=0$ on $\partial D\cap N$ (the boundary is characteristic near $x^0$). If $x^0\in\widetilde\Sig_1$, or if $x^0\in\Sig_1^0$ and (T) holds at $x^0$, then there is $\varepsilon_0>0$ such that, with $W:=D\cap B(x^0,\varepsilon_0)$,
\[
\Prob^b_x\big(\tau_W<\infty,\ X_{\tau_W-}\in\partial D\big)=0\qquad\text{for all }x\in W .
\]
\end{enumerate}
\end{lemma}

\begin{proof}
(i) We construct $U$ as in Lemma~\ref{lem:exit} on $W=D\cap B(x^0,\varepsilon_0)$ with $\varepsilon_0$ small; recall $h=\ip{a\nabla\rho}{\nabla\rho}$, and note that $\rho\le|x-x^0|\le\varepsilon_0$ on $W$ and that $L_b\rho$, $\tr a$, $b$ are bounded on $W$ by \eqref{eq:bdryreg}.

\emph{Case $x^0\in\Sig_3$.} Let $U:=\rho(2\varepsilon_0-\rho)+2\varepsilon_0|x-x^0|^2$. Then $U\ge0$ with equality only at $x^0$, and
\[
L_bU=(2\varepsilon_0-2\rho)L_b\rho-h+2\varepsilon_0\big(\tr a+2\ip{b}{x-x^0}\big)\le-h(x^0)/2+C\varepsilon_0\le-h(x^0)/4
\]
on $W$ for $\varepsilon_0$ small, using $h\ge h(x^0)/2$ near $x^0$.

\emph{Case $x^0\in\Sig_2$, boundary characteristic near $x^0$.} Let $U:=\rho^\alpha+\varepsilon|x-x^0|^2$ with $\alpha\in(0,1)$ and $\varepsilon>0$ to be chosen. Then $U\in C^2(W)\cap C(\bar W)$, $U\ge0$ with equality only at $x^0$, and
\[
L_b\rho^\alpha=\alpha\rho^{\alpha-1}\Big[L_b\rho-\tfrac12(1-\alpha)\frac{h}{\rho}\Big].
\]
By Lemma~\ref{lem:T} (and its proof; this is where the characteristic patch is used), as $x\to x^0$ the bracket tends to $L_b\rho(x^0)-\tfrac12(1-\alpha)\partial_nh(x^0)=-F_b(x^0)+\tfrac\alpha2\partial_nh(x^0)$. Choose $\alpha$ so small that this limit is $\le-F_b(x^0)/2<0$; then the bracket is $\le-F_b(x^0)/4$ on $W$ for $\varepsilon_0$ small, and since $\rho^{\alpha-1}\ge\varepsilon_0^{\alpha-1}$ on $W$ we get $L_b\rho^\alpha\le-2\gamma$ with $\gamma:=\alpha\varepsilon_0^{\alpha-1}F_b(x^0)/8$. Finally $L_b[\varepsilon|x-x^0|^2]=\varepsilon(\tr a+2\ip{b}{x-x^0})\le\gamma$ for $\varepsilon$ small, so $L_bU\le-\gamma$ on $W$. In both cases Lemma~\ref{lem:exit} gives (i).

(ii) Let $U:=-\log\rho$ on $W=D\cap B(x^0,\varepsilon_0)$ with $\varepsilon_0$ so small that \eqref{eq:T} holds on $W$ and $\rho\le e^{-1}$ there (so $U\ge1$); by Lemma~\ref{lem:T}, \eqref{eq:T} holds automatically near points of $\widetilde\Sig_1$ when the boundary is characteristic near $x^0$. Then
\[
L_bU=-\frac{L_b\rho}{\rho}+\frac{h}{2\rho^2}=\frac1{2\rho}\Big(\frac h\rho-2L_b\rho\Big)\le\frac C2\le\lambda U ,\qquad\lambda:=\tfrac C2\ (\text{any }\lambda>0\text{ if }C=0).
\]
Hence $e^{-\lambda t}U(X_t)$ is a nonnegative local supermartingale on $[0,\tau_W)$, and for $k\ge1$, $t>0$, optional sampling at $\inf\{s:U(X_s)\ge k\}\wedge\tau_W\wedge t$ gives $\Prob^b_x(\sup_{s<\tau_W\wedge t}U(X_s)\ge k)\le e^{\lambda t}U(x)/k$. Since $U(X_s)\to\infty$ as $s\uparrow\tau_W$ on $\{\tau_W\le t,\ X_{\tau_W-}\in\partial D\}$, this event has probability $0$ for every $t$.
\end{proof}

The Lyapunov functions in (i) are those of Friedman \cite[Ch.~11]{Fri75}; see also Stroock--Varadhan \cite{SV72b}, who obtain the boundary behaviour of degenerate diffusions under Fichera-type conditions by different means. The use of $-\log\rho$ in (ii) replaces the functions $\rho^{-\alpha}$ (for $\widetilde\Sig_1$) and the Lipschitz-based argument of \cite[Ch.~9, Thm.~4.1]{Fri75} (for $\Sig_1^0$) by a single computation, and isolates in (T) exactly what is needed at characteristic points with vanishing Fichera drift. The patch hypotheses in Lemmas~\ref{lem:T} and \ref{thm:local}, which require the boundary to be characteristic, and the Fichera drift to be nonpositive, on a whole neighbourhood of $x^0$ rather than at $x^0$ alone, cannot be dropped:

\begin{example}[isolated characteristic points]\label{ex:isolated}
In the half-plane $D=\{(x,r):r>0\}$ let $a=\operatorname{diag}(1,r+x^2)$ and $b=(0,1)$. At $x^0=(0,0)$ the boundary is characteristic and $F_b(x^0)=-(1-\tfrac12)<0$, but every boundary neighbourhood of $x^0$ consists otherwise of non-characteristic points $(x,0)$, $x\neq0$, which are attainable by Lemma~\ref{thm:local}(i); so no neighbourhood of $x^0$ is free of exits, and the conclusion of Lemma~\ref{thm:local}(ii) fails without the patch hypothesis. Similarly, with $a=\operatorname{diag}(1,r)$ and $b=(0,\tfrac12-x^2)$ the boundary is characteristic, $F_b(x,0)=x^2$ vanishes at $x^0$ only, the coefficients are Lipschitz in $r$, and yet $a_{rr}/r-2L_br=2x^2$ is not $O(r)$ near $x^0$: (T) fails at the point $x^0$ where $F_b=0$ because $F_b>0$ on the neighbouring patch. In Theorem~\ref{thm:fichera} both patch hypotheses hold on the whole boundary. We are not aware of these examples in the literature; the classical statements \cite{Fri75,SV72b} assume the boundary to be characteristic on a whole component and do not address isolated points.
\end{example}

\subsection{From local to global}\label{sec:global}
The local statements combine into a global one by a covering argument and the Markov property.

\begin{proposition}\label{thm:attain}
Let Assumption~\ref{ass:standing}(ii) hold for $L_b$ on $D$ and let \eqref{eq:bdryreg} hold at every point of $\partial D^{(2)}$.
\begin{enumerate}[label=(\roman*),leftmargin=2em]
\item If $\Sig_2\cup\Sig_3\neq\emptyset$, then $\Prob^b_x(\tau_D<\infty)>0$ for all $x\in D$.
\item Suppose $\Sig_3=\emptyset$, so that $\partial D^{(2)}$ is characteristic, and let $\Gamma\subset\partial D^{(2)}$ consist of points of $\widetilde\Sig_1$ and of points of $\Sig_1^0$ at which (T) holds. Then, for all $x\in D$,
\[
\Prob^b_x\big(\tau_D<\infty,\ X_{\tau_D-}\text{ exists and lies in }\Gamma\big)=0 .
\]
In particular, if $\Sig_2\cup\Sig_3=\emptyset$ and (T) holds on $\Sig_1^0$, the diffusion cannot leave $D$ through $\partial D^{(2)}$.
\end{enumerate}
\end{proposition}

\begin{proof}
(i) If $\Sig_3\neq\emptyset$, apply Lemma~\ref{thm:local}(i) at a point of $\Sig_3$; if $\Sig_3=\emptyset$, then $\partial D^{(2)}$ is characteristic and Lemma~\ref{thm:local}(i) applies at a point of $\Sig_2$. For (ii), for each $x^0\in\Gamma$ let $W(x^0)=D\cap B(x^0,\varepsilon_0(x^0))$ be as in Lemma~\ref{thm:local}(ii), and extract a countable family $W_k=D\cap B(x_k,\varepsilon_k)$ with $\Gamma\subset\bigcup_kB(x_k,\varepsilon_k)$. On the event in question there are $k$ and a rational $q<\tau_D$ with $X_t\in W_k$ for $t\in[q,\tau_D)$ and $X_{\tau_D-}\in\partial D\cap B(x_k,\varepsilon_k)$; then $\tau_D=q+\tau_{W_k}\circ\vartheta_q$ and $X_{\tau_{W_k}-}\circ\vartheta_q\in\partial D$. By the Markov property at time $q$, the probability of the event is at most
\[
\sum_{k,q}\E^b_x\Big[\Prob^b_{X_q}\big(\tau_{W_k}<\infty,\ X_{\tau_{W_k}-}\in\partial D\big);\ q<\tau_D,\ X_q\in W_k\Big]=0
\]
by Lemma~\ref{thm:local}(ii).
\end{proof}

Proposition~\ref{thm:attain} is a statement about an arbitrary generator $L_b$. Applied to $L$ it says what Assumption~\ref{ass:standing}(iii) forces on the coefficients of the model; applied to $\Lt$ it decides the exits of the drift-changed diffusion. The next proposition transfers boundary non-attainment from $L$ to $\Lt$ wherever $\theta$ is locally bounded.

\subsection{Bounded integrands}\label{sec:inherit}
Recall that, under Assumption~\ref{ass:standing}(iii), $Y$ does not leave $D$ under $\Prob$. The Girsanov density $S$, stopped upon leaving a region on which $\theta$ is bounded, is a true martingale by Novikov's criterion; so on such a region $\Qt$ is equivalent to $\Prob$, and whatever $Y$ cannot do under $\Prob$ it cannot do under $\Qt$ either.

\begin{proposition}\label{prop:inherit}
Let Assumption~\ref{ass:standing} hold, let $\Gamma\subset\partial D$, and let $N$ be an open set containing $\Gamma$ such that $\theta$ is bounded on $D\cap N$. If $\Prob_y\big(\tau_D<\infty,\ Y_{\tau_D-}\text{ exists and lies in }\Gamma\big)=0$ for all $y\in D$ (for instance by Assumption~\ref{ass:standing}(iii)), then
\[
\Qt_y\big(\tau_D<\infty,\ Y_{\tau_D-}\text{ exists and lies in }\Gamma\big)=0\qquad\text{for all }y\in D .
\]
\end{proposition}

\begin{proof}
Cover $\Gamma$ by countably many open balls $N_k\subset N$ and put $W_k:=D\cap N_k$, with exit time $\tau_k$. Fix $k$, $T>0$ and $y\in W_k$. Since $\theta$ is bounded on $W_k$, the stopped exponential $Z:=S^{\tau_k}$ is a true $\Prob_y$-martingale (Novikov) and $\Prob^k:=Z_T\cdot\Prob_y$ is a probability measure equivalent to $\Prob_y$ on $\mathcal F_T$. By Girsanov's theorem, under $\Prob^k$ the process $Y$ solves the SDE with drift $\mut$ up to $\tau_k\wedge T$, so by uniqueness in law for the drift-changed diffusion stopped at the exit from $W_k$ (Assumption~\ref{ass:standing}(ii) and \cite[Thm.~1.13.1]{Pin95}), the $\Prob^k$-law of $(Y_{t\wedge\tau_k\wedge T})_{t\ge0}$ is its $\Qt_y$-law. The event $E_k:=\{\tau_k\le T,\ Y_{\tau_k-}\text{ exists and lies in }\Gamma\}$ is determined by this stopped path, and on $E_k$ one has $\tau_k=\tau_D$; hence $\Qt_y(E_k)=\Prob^k(E_k)=0$, because $\Prob_y(E_k)=0$ by hypothesis and $\Prob^k\sim\Prob_y$. Finally, on $\{\tau_D<\infty,\ Y_{\tau_D-}\in\Gamma\}$ there are $k$ and a rational $q<\tau_D$ with $Y_t\in N_k$ for $t\in[q,\tau_D)$, and the Markov property at time $q$, as in the proof of Proposition~\ref{thm:attain}, shows that this event is $\Qt_y$-null.
\end{proof}

\subsection{Proof of Theorem~\ref{thm:fichera}}\label{sec:proof-fichera}
All the ingredients are now in place: Proposition~\ref{thm:attain} applied to $L$ gives (i), applied to $\Lt$ it gives the necessity (ii), and for the sufficiency (iii) the boundary is split into the part where $\theta$ is bounded, handled by Proposition~\ref{prop:inherit}, and the part where it is not, handled by Proposition~\ref{thm:attain}(ii).

\begin{proof}
(i) If $\ip{a\nu}{\nu}(x^0)>0$ or $F(x^0)>0$ at some $x^0\in\partial D$, Proposition~\ref{thm:attain}(i) applied to $L$ gives $\Prob_y(\tau_D<\infty)>0$, contradicting Assumption~\ref{ass:standing}(iii).

(ii) If $\Ft(x^0)>0$ for some $x^0$, Proposition~\ref{thm:attain}(i) applied to $\Lt$ gives $\Qt_y(\tau_D<\infty)>0$ for all $y$. By Lemma~\ref{thm:main}(ii), $S$ is not a martingale on $[0,\infty)$ from any starting point, and by Proposition~\ref{prop:finite}(i), for every $T>0$ it fails to be a martingale on $[0,T]$ from some starting point. Under local H\"older continuity it is a strict local martingale on every $[0,T]$ from every $y$ by Proposition~\ref{prop:finite}(ii).

(iii) By Lemma~\ref{thm:main} we must show $\Qt_y(\tau_D<\infty)=0$. Under \eqref{eq:lingrowth} there is no explosion at infinity by Proposition~\ref{cor:osgood}(i); so on $\{\tau_D<\infty\}$ the path is bounded, and since $\mut$ and $\sigma$ are bounded on bounded subsets of $\bar D$, the limit $Y_{\tau_D-}$ exists and lies in $\partial D$. Let $\Gamma_b$ be the set of boundary points near which $\theta$ is bounded (relatively open in $\partial D$) and $\Gamma_u:=\partial D\setminus\Gamma_b$. Choose a countable cover of $\Gamma_b$. Applying Proposition~\ref{prop:inherit} to each patch in this cover and taking a countable union shows that $Y_{\tau_D-}\in\Gamma_b$ has probability zero. By (i) every boundary point is characteristic ($\Sig_3=\emptyset$), and by hypothesis every point of $\Gamma_u$ lies in $\widetilde\Sig_1$ or in $\Sig_1^0$ with (T); Proposition~\ref{thm:attain}(ii) with $\Gamma=\Gamma_u$ shows that $Y_{\tau_D-}\in\Gamma_u$ has probability zero as well. Hence $\Qt_y(\tau_D<\infty)=0$. The equivalence \eqref{eq:iff} follows from (ii), (iii) and $\mut-\mu=\sigma\theta$.
\end{proof}

The following corollary shows that, where $\theta$ is bounded, the change of drift leaves the boundary Fichera drift unchanged; it uses only Theorem~\ref{thm:fichera}(i) and the Cauchy--Schwarz inequality.

\begin{corollary}\label{cor:inherit}
In the setting of Theorem~\ref{thm:fichera}, let $x^0\in\partial D$ and suppose $\theta$ is bounded on $D\cap N$ for a neighbourhood $N$ of $x^0$. Then $\Ft(x^0)=F(x^0)\le0$. Consequently $\Ft(x^0)>0$, hence loss of the martingale property due to a boundary exit, is possible only at boundary points near which $\theta$ is unbounded; and near boundary points where $\theta$ is bounded, non-attainment follows from Proposition~\ref{prop:inherit}.
\end{corollary}

\begin{proof}
By the Cauchy--Schwarz inequality, for $x\in D\cap N$ and $\nu=\nu(x^0)$, $|\ip{\sigma\theta}{\nu}(x)|=|\ip{\theta(x)}{\sigma(x)^{\top}\nu}|\le|\theta(x)|\sqrt{\ip{a(x)\nu}{\nu}}\to0$ as $x\to x^0$, since $a$ is continuous on $\bar D$ and $\ip{a(x^0)\nu}{\nu}=0$ by Theorem~\ref{thm:fichera}(i). As $\mut-\mu=\sigma\theta$ is continuous on $\bar D$, $\Ft(x^0)-F(x^0)=\lim_{x\to x^0}\ip{\sigma\theta}{\nu}(x)=0$.
\end{proof}

\begin{remark}[Fernholz--Karatzas]\label{rem:FK}
In \cite[Thm.~2]{FK10}, $U(T,y)$ is the survival probability of the auxiliary diffusion, so $1-U$ is the loss of mass. Their Proposition~3 treats the orthant and the specific auxiliary drift $\hat b_i(y)=\sum_j a_{ij}(y)/(\sum_j y_j)$. It gives survival when the inward Fichera drifts are nonnegative on every face, and exit when they are negative on every face. Our Theorem~\ref{thm:fichera} treats a general shift $\sigma\theta$ and implies a positive probability of exit if the outward Fichera drift $\Ft$ is positive at any smooth boundary point. Continuity makes $\Ft$ positive on a nearby boundary patch, to which the local accessibility argument applies. The theorem also states the critical condition (T) explicitly. Both use classical boundary accessibility; square-root degeneracy already occurs in \cite[Sec.~12]{FK10}.
\end{remark}

\subsection{Boundary points with vanishing Fichera drift}\label{sec:tangential}
Theorem~\ref{thm:fichera} is an ``iff'' at boundary points where $\Ft\ne0$, and at points where $\Ft=0$ under (T). The following examples explain the role and limitations of (T).

(a) Mere continuity of $\mut$ up to the boundary is not enough at characteristic points with vanishing Fichera drift: in $d=1$, $a(y)=y$, $b(y)=\tfrac12-1/\log(1/y)$ has $F_b(0)=0$ but $0$ is attainable (the scale function is finite at $0$), whereas $b(y)=\tfrac12-y^\gamma$, $\gamma>0$, is not. This is why Theorem~\ref{thm:fichera}(iii) needs (T) on $\{\Ft=0\}$.

(b) The difference between (T) for $L$ and for $\Lt$ is the term $-2\ip{\sigma\theta}{\nabla\rho}$, which is $O(\rho)$ if $\sigma$ is Lipschitz and $\theta$ bounded, and only $O(\sqrt\rho)$ under the standing assumptions. The one-dimensional Feller test shows that an $O(\sqrt\rho)$ perturbation of the drift does not change attainability at a CIR-type boundary point where the Fichera drift vanishes, so (T) is stronger than necessary there; a sharper local condition would need a Lyapunov function of the form $-\log\rho+c\sqrt\rho$, or the analytic machinery developed for Kimura- and Heston-type operators at such boundaries by Epstein--Mazzeo \cite{EM13} and Feehan--Pop \cite{FP15}.

(c) The following example shows that a zero Fichera drift at a boundary point does not determine whether the shifted diffusion can reach that boundary point. Since all other exits are excluded in the example, this attainability question determines whether $S$ is a true martingale on every finite horizon. The example also shows that (T) is sufficient but not necessary for boundary non-attainment and gives the exact criterion for this family.

\begin{example}[the critical boundary]\label{ex:critical}
Let $d=1$, $D=(0,\infty)$, $a(y)=y$, $\mu\equiv\tfrac12$ (a CIR process at the critical Feller parameter, so $0$ is inaccessible under $\Prob$ and Assumption~\ref{ass:standing}(iii) holds), and for $c\ge0$ let
\[
\theta(y)=-\frac{c}{\sqrt y\,\log(1/y)}\quad\text{for }y\le e^{-2},
\]
extended smoothly on $(0,\infty)$ and cut off for large $y$; $\theta$ is locally bounded on $D$ and, for $c>0$, unbounded at $0$. Then $\sigma\theta=-c/\log(1/y)$ near $0$, so $\mut(y)=\tfrac12-c/\log(1/y)$ there, $\mut$ is continuous on $\bar D$, and the Fichera drifts of both $L$ and $\Lt$ vanish at $0$ for every $c$. The scale density of $\Lt$ near $0$ is
\[
p(y)=\exp\Big(-\int^y\frac{2\mut(z)}{a(z)}dz\Big)=\frac{1}{y\,[\log(1/y)]^{2c}}
\]
up to a constant, and $\int_{0}p(y)\,dy=\infty$ iff $2c\le1$; for $c>\tfrac12$ the scale integral is finite and the speed-scale integral $\int_0(s(y)-s(0))m(y)\,dy$, whose integrand is of order $\log(1/y)$, is finite too, so $0$ is accessible. By Feller's test and Lemma~\ref{thm:main} (there is no explosion at infinity since $\mut$ is bounded and $a$ has linear growth),
\[
S\ \text{is a true martingale on every finite horizon}\iff c\le\tfrac12 .
\]
On the other hand, $a(y)/y-2\mut(y)=2c/\log(1/y)$ is not $O(y)$ for any $c>0$, so (T) fails for every $c>0$, including the martingale cases $0<c\le\tfrac12$. Adding independent tangential Brownian coordinates embeds the example in any dimension with a characteristic hyperplane boundary. For this family, the exact condition comes from Feller's test in the normal variable and is not captured by the sign of the Fichera drift or by (T). Theorem~\ref{thm:coupled-log} allows the logarithmic coefficient to depend on correlated tangential factors and proves both directions of the resulting threshold.
\end{example}

\subsection{A coupled logarithmic threshold}\label{sec:log-extension}

Example~\ref{ex:critical} shows that the sign of the Fichera drift and condition (T) do not characterize the critical boundary. The following result allows the coefficient of the logarithmic correction to depend on correlated factors.

\begin{theorem}\label{thm:coupled-log}
Let $(W^0,W^1,\ldots,W^k)$, $k\ge1$, be a Brownian motion with a positive definite constant correlation matrix. On $(0,\infty)\times\R^k$, consider
\[
dX_t=(\tfrac12-\kappa X_t)\,dt+\sqrt{X_t}\,dW^0_t,
\qquad dZ^i_t=-\lambda_iZ^i_t\,dt+\eta_i\,dW^i_t,
\]
where $\kappa,\lambda_i,\eta_i>0$. Let $c:\R^k\to[0,\infty)$ be bounded and locally H\"older continuous. Choose a nonnegative smooth function $h$ on $(0,\infty)$ such that
\[
h(x)=\frac1{\log(1/x)}\quad(0<x\le e^{-2}),
\qquad h(x)=0\quad(x\ge e^{-1}),
\]
and put $S=\cE(-\int c(Z)h(X)X^{-1/2}\,dW^0)$. For every $T>0$ and every starting point,
\begin{equation}\label{eq:coupled-log-iff}
S\text{ is a true martingale on }[0,T]
\quad\Longleftrightarrow\quad \sup_{z\in\R^k}c(z)\le\tfrac12.
\end{equation}
If the inequality fails, $\E_{x,z}[S_T]<1$ for every $T>0$ and every starting point.
\end{theorem}

\begin{proof}
The original critical CIR factor is strictly positive and conservative, and the Ornstein--Uhlenbeck factors are conservative. Set $\bar c=\sup c$ and $\rho_i=d\langle W^0,W^i\rangle_t/dt$. The shifted drifts are
\[
\mut_X=\tfrac12-\kappa x-c(z)h(x),\qquad
\mut_{Z^i}=-\lambda_i z_i-\eta_i\rho_i\frac{c(z)h(x)}{\sqrt x}.
\]
The joint covariance is locally uniformly positive definite on the state space, and stopped Girsanov gives local well-posedness of the shifted system.

Suppose $\bar c\le1/2$. Stop $X$ in $(1/n,n)$, so the stopped integrand is bounded. Under the probability measure defined by the stopped density, compare $X$ with scalar diffusions having the same square-root diffusion coefficient and respective drifts
\[
\tfrac12-\kappa x-\bar c h(x)
\quad\text{and}\quad \tfrac12-\kappa x.
\]
The lower comparison has scale density proportional, near $0$, to
\[
\frac{e^{2\kappa x}}{x[\log(1/x)]^{2\bar c}}.
\]
Its scale integral at $0$ diverges for $2\bar c\le1$, making $0$ inaccessible. Both comparisons are non-explosive at infinity, since $h$ vanishes for large $x$. Their probabilities of leaving $(1/n,n)$ before $T$ tend to zero. Writing $\tau_n$ for the stopping time and $\Q^{n,T}$ for this measure, original non-exit and monotone convergence give $\E_{x,z}[S_T]=\lim_n\Q^{n,T}(\tau_n>T)=1$.

If $\bar c>1/2$, choose a closed box $K\subset\R^k$, $z^0\in\operatorname{int}K$, and $c_*>1/2$ such that $c(z)\ge c_*$ on $K$. Choose $0<\beta<2c_*-1$, write $\ell=\log(1/x)$, and set $f(x)=\ell^{-\beta}$. On $(0,\varepsilon)\times K$, for sufficiently small $\varepsilon$,
\begin{equation}\label{eq:log-barrier}
\Lt f(x)=\frac{\beta((\beta+1)/2-c(z))}{x\ell^{\beta+2}}
-\frac{\kappa\beta}{\ell^{\beta+1}}
\le-\frac{\delta}{x\ell^{\beta+2}}
\end{equation}
for some $\delta>0$. The weight $w(x)=1/(x\ell^{\beta+2})$ tends to infinity and dominates $1+x^{-1/2}$ near $0$. For the first exit $\sigma_n$ from $(1/n,\varepsilon)\times\operatorname{int}K$, stopped It\^o yields
\[
\E^{\Qt}_{x,z^0}[f(X_{\sigma_n})]
+\delta\E^{\Qt}_{x,z^0}\int_0^{\sigma_n}w(X_s)ds\le f(x).
\]
On this cylinder the tangential drift is bounded by $C(1+x^{-1/2})$, and its covariance is bounded. Thus the estimate controls the expected stopping time, the expected total variation of that drift, and the expected squared maximum of the tangential martingale by $Cf(x)$. In particular the probability of a tangential displacement reaching $\dist(z^0,\partial K)$ is at most $Cf(x)$, and the probability of reaching $x=\varepsilon$ is at most $f(x)/f(\varepsilon)$. Since $f(x)\to0$, with positive probability $\sigma=\lim_n\sigma_n$ is finite and the process approaches $x=0$ while $Z$ stays in $\operatorname{int}K$. Integrability of the drifts and finiteness of the quadratic variations ensure that the coordinates have limits at $\sigma$. Lemma~\ref{thm:main} and Proposition~\ref{prop:finite}(ii) now give strictness for every starting point and finite horizon.
\end{proof}

For independent normal and tangential drivers, the boundary coefficients satisfy the regularity hypotheses of Section~\ref{sec:boundary}. At $x=0$ the shifted Fichera drift is zero, whereas
\[
\frac{a_{XX}}x-2\mut_X=2\kappa x+\frac{2c(z)}{\log(1/x)}
\qquad(x\le e^{-2}).
\]
Consequently (T) fails wherever $c(z)>0$, including martingale cases in \eqref{eq:coupled-log-iff}. The scalar threshold is Feller's test from Example~\ref{ex:critical}. With correlated drivers, the occupation bound \eqref{eq:log-barrier} controls the singular tangential shift and yields necessity on any neighbourhood where $c>1/2$.

For a concrete coupled family, take $k=1$, any correlation $\rho_1\in(-1,1)$, and $c(z)=a+bz^2/(1+z^2)$ with $a,b\ge0$. The exact condition is $a+b\le1/2$. When $b>0$ and $a+b=1/2$, the critical supremum is approached only as $|z|\to\infty$, and the density is still a true martingale.

Uniform integrability over the infinite horizon is characterized through global absolute continuity in \cite[Cor.~5.1]{CG18} and \cite[Prop.~2.6]{Cri21}. For the present mean-reverting model, recurrence gives a complete classification.

\begin{corollary}\label{cor:coupled-log-ui}
Under the hypotheses of Theorem~\ref{thm:coupled-log}, for every starting point, $S$ is a uniformly integrable martingale on $[0,\infty)$ if and only if $c\equiv0$. If $c\not\equiv0$, then
\[
\lim_{t\to\infty}S_t=0\qquad\Prob_{x,z}\text{-almost surely}.
\]
In particular, if $0<\sup c\le1/2$, then $S$ is a true martingale on every finite horizon and is not uniformly integrable on $[0,\infty)$.
\end{corollary}

\begin{proof}
We first establish recurrence of the original factor process $(X,Z)$ on $E=(0,\infty)\times\R^k$. Its coefficients are smooth in $E$, and its covariance matrix is locally uniformly positive definite because the Brownian correlation matrix is positive definite. The process is conservative. On a smooth, connected, relatively compact exhaustion $E_n$ of $E$, the killed transition kernels have positive continuous densities in the interior, by local elliptic diffusion theory \cite[Ch.~1]{Pin95}. Hence the process is irreducible with respect to Lebesgue measure and is a $T$-process in the sense of \cite[Sec.~3.2]{MT93III}. Indeed, the sum of the killed time-one kernels with weights $2^{-n}$ is a lower semicontinuous kernel dominated by the full time-one kernel, with positive total mass at every starting point.

Write $\ell(x)=\log(1/x)$. Choose a nonnegative $C^2$ function $\psi$ on $(0,\infty)$ equal to $\log\ell(x)$ for $x\le e^{-3}$ and to zero for $x\ge e^{-2}$, and put
\[
V(x,z)=1+x+|z|^2+\psi(x).
\]
The sublevel sets of $V$ are compact in $E$. For the original generator $L$, direct differentiation gives, when $x\le e^{-3}$,
\begin{align*}
L\psi&=-\frac{1}{2x\ell(x)^2}+\frac{\kappa}{\ell(x)},\\
LV&=\frac12-\kappa x+\sum_{i=1}^k(\eta_i^2-2\lambda_i z_i^2)
       -\frac{1}{2x\ell(x)^2}+\frac{\kappa}{\ell(x)}.
\end{align*}
The mixed derivatives of $V$ vanish, so the correlations contribute no terms to this calculation. The $x$-dependent part of the global expression for $LV$ is $1/2-\kappa x+L\psi(x)$. It tends to $-\infty$ as $x\downarrow0$ by the displayed calculation, and as $x\to\infty$ because $\psi=0$ for $x\ge e^{-2}$. It is bounded above on the intervening compact interval. Since every $\lambda_i$ is positive, there are a compact set $K\subset E$ and a finite constant $b$ such that
\[
LV\le-1+b\mathbf1_K\qquad\text{on }E.
\]
It\^o's formula on the exhaustion supplies this inequality for the stopped extended generators. The Foster--Lyapunov criterion \cite[Thm.~4.4]{MT93III} therefore gives positive Harris recurrence. In particular, by the occupation-time characterization of Harris recurrence \cite[Sec.~3.2]{MT93III}, every open set $U\subset E$ of positive Lebesgue measure satisfies
\[
\int_0^\infty\mathbf1_U(X_t,Z_t)\,dt=\infty
\qquad\Prob_{x,z}\text{-almost surely for every }(x,z)\in E.
\]

Suppose $c\not\equiv0$. Continuity gives a bounded open ball $U_z$ with compact closure and a constant $c_0>0$ such that $c\ge c_0$ on $\overline U_z$. Choose an interval $U_x$ with $\overline U_x\subset(0,e^{-3})$ and set $U=U_x\times U_z$. On $\overline U$, the function $c(z)^2h(x)^2/x$ is bounded below by some $\varepsilon>0$. Thus the continuous local martingale
\[
M_t=-\int_0^t\frac{c(Z_s)h(X_s)}{\sqrt{X_s}}\,dW_s^0
\]
has quadratic variation satisfying
\[
\langle M\rangle_\infty
=\int_0^\infty\frac{c(Z_t)^2h(X_t)^2}{X_t}\,dt
\ge\varepsilon\int_0^\infty\mathbf1_U(X_t,Z_t)\,dt=\infty
\quad\Prob_{x,z}\text{-almost surely}.
\]
The Dambis--Dubins--Schwarz theorem \cite[Thm.~3.4.6]{KS91} gives $M_t=B_{\langle M\rangle_t}$ for a Brownian motion $B$. Since $B_u/u\to0$ almost surely,
\[
\frac{\log S_t}{\langle M\rangle_t}
=\frac{B_{\langle M\rangle_t}}{\langle M\rangle_t}-\frac12
\longrightarrow-\frac12,
\]
so $S_t\to0$ almost surely. A uniformly integrable martingale with $S_0=1$ has $\E[S_\infty]=1$, which excludes that possibility. If $c\equiv0$, then $S\equiv1$. The final assertion follows from Theorem~\ref{thm:coupled-log}.
\end{proof}

\section{Explosion at infinity}\label{sec:infty}

After exits through finite boundary points have been excluded, martingality is equivalent to non-explosion at infinity. Exact radial martingality criteria for multidimensional diffusions were established by Criens--Glau \cite[Cor.~5.2]{CG18}. Khasminskii's radial tests bound the squared radius by scalar diffusions (Lemma~\ref{thm:radial-sandwich}); when the bounds coincide, Proposition~\ref{thm:radial-iff} gives a necessary and sufficient condition. Proposition~\ref{cor:osgood} identifies the linear-growth and Osgood regimes.

The coefficient structure that makes comparison exact is the focus of Proposition~\ref{prop:adapted-infinity}: a scalar function with compact sublevel sets has drift and variance depending on its value, up to a common time-change rate bounded above and away from zero. Corollary~\ref{cor:variable-clock} also allows rates comparable to a positive function of the scalar coordinate. This gives an exact martingality criterion for the anisotropic models of Example~\ref{ex:anisotropic}, including nonlinear motion tangent to the level sets and a time-change rate that varies along them. Section~\ref{sec:LMDP} compares this formulation with the radial criteria of Criens--Glau, the time-change criterion of Criens, and the sufficient criteria of Lions--Musiela and Dandapani--Protter.

\subsection{The two-sided test}
For $y\in D$, $y\neq0$, set
\begin{equation}\label{eq:AB}
\mathcal A(y):=\ip{y}{a(y)y}>0,\qquad \mathcal B(y):=\frac{\tr a(y)+2\ip{y}{\mut(y)}}{\ip{y}{a(y)y}},\qquad r:=\tfrac12|y|^2 .
\end{equation}
For $V(y)=f(\tfrac12|y|^2)$ with $f\in C^2$, a direct computation gives
\begin{equation}\label{eq:radialgen}
\Lt V(y)=\tfrac12\mathcal A(y)\big[f''(r)+\mathcal B(y)f'(r)\big].
\end{equation}
Given continuous $A:[r_0,\infty)\to(0,\infty)$ and $B:[r_0,\infty)\to\R$, put $C(\rho):=\exp\int_{r_0}^\rho B(s)ds$ and
\begin{equation}\label{eq:feller}
v_{A,B}(\infty):=\int_{r_0}^\infty\frac1{C(\rho)}\int_{r_0}^\rho\frac{C(s)}{A(s)}\,ds\,d\rho\ \in(0,\infty] .
\end{equation}
$v_{A,B}(\infty)$ is the integral appearing in Feller's test for the one-dimensional diffusion with generator $\tfrac12A(r)\big[\tfrac{d^2}{dr^2}+B(r)\tfrac{d}{dr}\big]$: its scale density is $1/C$, its speed density $2C/A$, and $+\infty$ is reached in finite time with positive probability iff $v_{A,B}(\infty)<\infty$ \cite[Thm.~5.5.29]{KS91}. The radial test compares $\tfrac12|Y|^2$ with this diffusion, from above and from below. The following form of the test is given in \cite[Thms.~10.2.3--10.2.4]{SV06} and \cite[Thms.~6, 8]{DP22}.

\begin{lemma}\label{thm:radial-sandwich}
Let Assumption~\ref{ass:standing} hold, $r_0>0$, and let $A,B$ be as above.
\begin{enumerate}[label=(\roman*),leftmargin=2em]
\item If $A(r)\ge\mathcal A(y)$ and $B(r)\ge\mathcal B(y)$ for all $y\in D$ with $\tfrac12|y|^2=r\ge r_0$, and $v_{A,B}(\infty)=\infty$, then $\Qt$ does not explode at infinity: $\Qt_y(\tau_D<\infty,\ \sup_{t<\tau_D}|Y_t|=\infty)=0$ for all $y\in D$.
\item If $D$ is unbounded, $A(r)\le\mathcal A(y)$ and $B(r)\le\mathcal B(y)$ for all $y\in D$ with $\tfrac12|y|^2=r\ge r_0$, and $v_{A,B}(\infty)<\infty$, then $\Qt_y(\tau_D<\infty)>0$ for all $y\in D$, so $S$ is not a martingale on $[0,\infty)$ (and, by Proposition~\ref{prop:finite}, fails to be a martingale on every $[0,T]$ for some starting point).
\end{enumerate}
\end{lemma}

\begin{proof}
Let $f(r):=\int_{r_0}^r C(\rho)^{-1}\int_{r_0}^\rho C(s)A(s)^{-1}ds\,d\rho$ for $r\ge r_0$, $f:=0$ on $[0,r_0]$. Then $f\in C^{1,1}$, $f\ge0$, $f'\ge0$, $f''+Bf'=1/A$ a.e., and $f(\infty)=v_{A,B}(\infty)$. Put $V:=f(\tfrac12|y|^2)$; It\^o's formula applies (the Meyer--It\^o formula for the semimartingale $\tfrac12|Y|^2$ and the $C^{1,1}$ function $f$).

(i) On $\{r\ge r_0\}$, by \eqref{eq:radialgen}, $\Lt V=\tfrac12\mathcal A\big[1/A(r)-(B(r)-\mathcal B)f'\big]\le\tfrac12\mathcal A/A(r)\le\tfrac12$, and $\Lt V=0$ on $\{r<r_0\}$. Hence $\widetilde V:=V+\tfrac12$ satisfies $\Lt\widetilde V\le\widetilde V$ on $D$, so $e^{-t}\widetilde V(Y_t)$ is a nonnegative local supermartingale on $[0,\tau_D)$ and $\Qt_y(\sup_{s<t\wedge\tau_D}\widetilde V(Y_s)\ge k)\le e^t\widetilde V(y)/k$. Since $\widetilde V(y)\to\infty$ as $|y|\to\infty$, the path cannot be unbounded on $[0,t\wedge\tau_D)$; this holds for every $t$.

(ii) Now $f$ is bounded, $f(r_0)=0<f(r)$ for $r>r_0$, and on $D\setminus B$ with $B:=D\cap\{r\le r_0\}$ we get $\Lt V=\tfrac12\mathcal A\big[1/A(r)+(\mathcal B-B(r))f'\big]\ge\tfrac12\mathcal A/A(r)\ge\tfrac12\ge\lambda V$ with $\lambda:=1/(2f(\infty))$. Since $\sup_{\partial_DB}V=f(r_0)=0<V(y_0)$ for any $y_0\in D\setminus B$ (such points exist as $D$ is unbounded), Lemma~\ref{thm:lyap}(ii) applies.
\end{proof}

\begin{remark}[the Ikeda--Watanabe comparison]\label{rem:IW}
Lemma~\ref{thm:radial-sandwich} is a Lyapunov-function version of the comparison theorem of Ikeda--Watanabe \cite[Ch.~VI, Thm.~4.1]{IW89}, stated there for diffusions on the whole space: for $p\in C^2(\R^d)$ with range $I$, put $\mathfrak a:=\ip{a\nabla p}{\nabla p}$ and $\mathfrak b:=\Lt p/\mathfrak a$, let $\mathfrak a_\pm(z)$, $\mathfrak b_\pm(z)$ be the supremum and the infimum of $\mathfrak a$, $\mathfrak b$ over the level set $\{p=z\}$, assumed locally Lipschitz with $\mathfrak a_\pm>0$, and let $Z^{\pm\pm}$ be the scalar diffusions with generators $\mathfrak a_\pm(z)\big[\tfrac12\tfrac{d^2}{dz^2}+\mathfrak b_\pm(z)\tfrac{d}{dz}\big]$. Then $p(Y)$ can be coupled with the $Z^{\pm\pm}$ so that its running maximum lies between those of $Z^{--}$ and $Z^{++}$, and its running minimum between those of $Z^{+-}$ and $Z^{-+}$. With $p(y)=\tfrac12|y|^2$ one has $\mathfrak a=\mathcal A$, $\mathfrak b=\tfrac12\mathcal B$, and Lemma~\ref{thm:radial-sandwich} is the statement about $\sup I=\infty$. On a proper domain $D$ the killed projection $p(Y)$ is not a diffusion on $I$ and the conclusions must be read with the other exits included: the conclusion asserting exit becomes ``$Y$ leaves $D$, either through $\partial D$ or with $p(Y)\to\partial I$'', which is what Lemma~\ref{thm:radial-sandwich}(ii) asserts, while the conclusion excluding exit at the scalar endpoint only excludes $p(Y)\to\partial I$ before $Y$ leaves $D$ (Brownian motion killed on the unit disc, projected on the first coordinate, illustrates the point: the scalar comparison diffusion reaches $\pm1$, the projection does not). The Lyapunov proofs of Lemma~\ref{thm:radial-sandwich} are the safer basis in the killed case. The gap between the ``$+$'' and the ``$-$'' comparison diffusions is the gap between (i) and (ii) of the lemma; a monotone reparametrization $q=g(p)$ does not change the endpoint accessibility of the reference diffusions, so the choice among $p$, $p^2$, etc.\ matters only through the level-set family and the tightness of the envelopes, and other families (e.g.\ $\sum_i\lambda_iy_i^2$, or $\sum_i\log y_i$ on the orthant) can close the gap in specific models.
\end{remark}

\subsection{Linear growth, Osgood growth, and the radial case}
The two extreme cases of Lemma~\ref{thm:radial-sandwich} are the most useful in practice: a quadratic bound on $\tr a+2\ip{y}{\mut}$, which is the hypothesis \eqref{eq:lingrowth} of Theorem~\ref{thm:fichera}, excludes explosion, and a radial drift of Osgood growth that dominates the noise forces it. These are the multi-dimensional forms of the criteria of Lions--Musiela \cite[Thms.~2.4, 3.1]{LM07}.

\begin{proposition}\label{cor:osgood}
Let Assumption~\ref{ass:standing} hold.
\begin{enumerate}[label=(\roman*),leftmargin=2em]
\item If $\tr a(y)+2\ip{y}{\mut(y)}\le C(1+|y|^2)$ for all $y\in D$, then $\Qt$ does not explode at infinity.
\item Suppose $D$ is unbounded and that for $y\in D$ with $|y|\ge R$:
\begin{enumerate}[label=(\alph*),leftmargin=2em]
\item $c_1w(|y|)\le\ip{y}{a(y)y}\le C_1w(|y|)$ for a continuous $w:[R,\infty)\to(0,\infty)$ and constants $0<c_1\le C_1$ (the radial noise is comparable across spheres);
\item $\ip{y}{\mut(y)}\ge|y|\,\varphi(|y|)$ for a nondecreasing $\varphi\in C^1([R,\infty);(0,\infty))$ with
\begin{equation}\label{eq:osgood}
\int_R^\infty\frac{d\xi}{\varphi(\xi)}<\infty\qquad\text{(Osgood condition)};
\end{equation}
\item $\inf_{\xi\ge R}\ \xi^2\varphi(\xi)/w(\xi)>0$.
\end{enumerate}
Then $\Qt_y(\tau_D<\infty)>0$ for all $y\in D$, and $S$ is not a martingale on $[0,\infty)$.
\end{enumerate}
\end{proposition}

\begin{proof}
(i) With $V(y)=1+|y|^2$ one has $\Lt V=\tr a+2\ip{y}{\mut}\le CV$, and the supermartingale argument of the proof of Lemma~\ref{thm:radial-sandwich}(i) applies.

(ii) Put $r_0:=R^2/2$, $A(r):=c_1w(\sqrt{2r})$ and $B(r):=2\sqrt{2r}\,\varphi(\sqrt{2r})/(C_1w(\sqrt{2r}))$ for $r\ge r_0$. On $\{\tfrac12|y|^2=r\}$ we have $\mathcal A(y)\ge A(r)$ and $\mathcal B(y)\ge2\ip{y}{\mut}/\ip{y}{ay}\ge B(r)$, so by Lemma~\ref{thm:radial-sandwich}(ii) it suffices to show $v_{A,B}(\infty)<\infty$. Substituting $s=\xi^2/2$, $\rho=\eta^2/2$, $u=\zeta^2/2$ in \eqref{eq:feller} and writing $\hat g(\zeta):=\zeta B(\zeta^2/2)=2\zeta^2\varphi(\zeta)/(C_1w(\zeta))$, $M(\xi):=\xi\varphi(\xi)$, $K:=C_1/c_1$, we get $\xi/A(\xi^2/2)=K\hat g(\xi)/(2M(\xi))$ and, by Fubini,
\[
v_{A,B}(\infty)=K\int_R^\infty\eta\int_R^\eta\frac{\hat g(\xi)}{2M(\xi)}\,e^{-\int_\xi^\eta\hat g}\,d\xi\,d\eta .
\]
Since $\hat g(\xi)e^{-\int_\xi^\eta\hat g}=\partial_\xi e^{-\int_\xi^\eta\hat g}$, an integration by parts gives
\[
\int_R^\eta\frac{\hat g(\xi)}{2M(\xi)}e^{-\int_\xi^\eta\hat g}d\xi
=\frac1{2M(\eta)}-\frac{e^{-\int_R^\eta\hat g}}{2M(R)}+\int_R^\eta e^{-\int_\xi^\eta\hat g}\frac{M'(\xi)}{2M(\xi)^2}d\xi
\le\frac1{2M(\eta)}+\int_R^\eta e^{-\int_\xi^\eta\hat g}\frac{M'(\xi)}{2M(\xi)^2}d\xi ,
\]
where $M'\ge0$ because $\varphi$ is nondecreasing. By (c), $\hat g\ge g_0>0$ on $[R,\infty)$, hence $\hat E(\xi):=\int_\xi^\infty\eta e^{-\int_\xi^\eta\hat g}d\eta\le\xi/g_0+1/g_0^2$, and another application of Fubini yields
\[
v_{A,B}(\infty)\le K\int_R^\infty\frac{d\eta}{2\varphi(\eta)}+K\int_R^\infty\frac{M'(\xi)}{2M(\xi)^2}\Big(\frac{\xi}{g_0}+\frac1{g_0^2}\Big)d\xi<\infty ,
\]
because $\int_R^\infty\xi M'/M^2\,d\xi=R/M(R)+\int_R^\infty d\xi/(\xi\varphi(\xi))<\infty$ and $\int_R^\infty M'/M^2\,d\xi=1/M(R)$ by \eqref{eq:osgood}.
\end{proof}

When $\mathcal A$ and $\mathcal B$ are constant on each sphere outside a ball, the upper and lower comparison functions coincide, and the test is exact.

\begin{proposition}\label{thm:radial-iff}
Let Assumption~\ref{ass:standing} hold, let $D$ be unbounded, and suppose that $\mathcal A$ and $\mathcal B$ are radial outside a ball: $\mathcal A(y)=A(\tfrac12|y|^2)$, $\mathcal B(y)=B(\tfrac12|y|^2)$ for $y\in D$, $|y|^2\ge2r_0$, with $A,B$ continuous. Then
\[
\Qt\ \text{explodes at infinity with positive probability}\ \Longrightarrow\ v_{A,B}(\infty)<\infty\ \Longrightarrow\ \Qt_y(\tau_D<\infty)>0\ \ \forall y\in D ,
\]
and $v_{A,B}(\infty)=\infty$ implies no explosion at infinity. In particular, if $D=\R^d$, or more generally if $\Qt$ cannot leave $D$ through finite points of $\partial D$ (by Section~\ref{sec:boundary} and Proposition~\ref{prop:sing}), then
\[
S\ \text{is a true martingale}\iff v_{A,B}(\infty)=\infty .
\]
\end{proposition}

\begin{proof}
Both implications are Lemma~\ref{thm:radial-sandwich} with equalities in place of inequalities.
\end{proof}

The unboundedness of $D$ cannot be dropped: on a bounded domain the exterior hypotheses are vacuous for every $A,B$, while the conclusion may fail.

The test only sees explosion at infinity. Exits through finite boundary points, a volatility factor hitting zero or a coefficient singularity, are invisible to it and are the subject of Sections~\ref{sec:boundary} and \ref{sec:sing}.

\subsection{Comparison using a function with compact sublevel sets}\label{sec:adapted-infinity}

Remark~\ref{rem:IW} permits comparison along functions other than $|y|^2/2$. The following exact case allows a time change whose rate is bounded above and away from zero and may vary along each level set.

\begin{proposition}\label{prop:adapted-infinity}
Let Assumption~\ref{ass:standing} hold with $D=\R^d$ and locally H\"older coefficients $a,\mut$. Suppose there are an unbounded function $P\in C^2(\R^d;[0,\infty))$ with compact sublevel sets, a number $r_0>0$, continuous functions $A:[r_0,\infty)\to(0,\infty)$ and $B:[r_0,\infty)\to\R$, and a function $q$ satisfying $0<q_-\le q(y)\le q_+<\infty$ on $\{P\ge r_0\}$, such that
\begin{equation}\label{eq:adapted-coordinate}
\ip{a\nabla P}{\nabla P}=q(y)A(P(y)),
\qquad 2\Lt P=q(y)A(P(y))B(P(y))
\end{equation}
on $\{P\ge r_0\}$. Then, for every $T>0$ and every starting point,
\[
S\text{ is a true martingale on }[0,T]
\quad\Longleftrightarrow\quad v_{A,B}(\infty)=\infty,
\]
where $v_{A,B}$ is given by \eqref{eq:feller}.
\end{proposition}

\begin{proof}
Let $f(r)=v_{A,B}(r)$ be the integral in \eqref{eq:feller} with outer upper limit $r$. Then $f(r_0)=f'(r_0)=0$, $f'>0$ on $(r_0,\infty)$, and $f''+Bf'=1/A$. For $V=f\circ P$, the chain rule gives
\[
\Lt V=\tfrac12q(y)A(P)\,[f''(P)+B(P)f'(P)]=\tfrac12q(y)
\qquad(P>r_0).
\]
If $f(\infty)=\infty$, extend $f$ as a $C^2$ function below $r_0$ and add a constant so that $V\ge1$. Then $V(y)\to\infty$ as $|y|\to\infty$; its generator is bounded above outside a compact set, and Lemma~\ref{thm:lyap}(i) excludes explosion. If $f(\infty)<\infty$, then $V$ is bounded and positive on $\{P>r_0\}$, zero on its boundary, and
\[
\Lt V\ge q_-/2\ge\frac{q_-}{2f(\infty)}V.
\]
Lemma~\ref{thm:lyap}(ii), with the compact set $\{P\le r_0\}$, gives explosion with positive probability. Lemma~\ref{thm:main} and Proposition~\ref{prop:finite}(ii) give the stated finite-horizon equivalence.
\end{proof}

The generator identity \eqref{eq:adapted-coordinate} identifies $P(Y)$ outside a compact set with a scalar diffusion under a time change whose rate is bounded above and below. Finiteness of the explosion time is preserved by this time change.

A rate that varies with the scalar coordinate can be incorporated into the variance coefficient of the reference diffusion. The following corollary makes this dependence explicit. Related connections between accumulated time and explosion are developed in \cite[Thm.~2.2 and Sec.~3.2]{Cri21}.

\begin{corollary}\label{cor:variable-clock}
Let Assumption~\ref{ass:standing} hold with $D=\R^d$ and locally H\"older coefficients $a,\mut$. Let $P\in C^2(\R^d;[0,\infty))$ be unbounded with compact sublevel sets, let $r_0>0$, and let $A,g:[r_0,\infty)\to(0,\infty)$ and $B:[r_0,\infty)\to\R$ be continuous. Suppose that on $\{P\ge r_0\}$,
\[
\ip{a\nabla P}{\nabla P}=q(y)A(P(y)),
\qquad 2\Lt P=q(y)A(P(y))B(P(y)),
\]
and that constants $0<q_-\le q_+<\infty$ satisfy
\begin{equation}\label{eq:variable-clock}
q_-g(P(y))\le q(y)\le q_+g(P(y)).
\end{equation}
Then, for every $T>0$ and every starting point,
\[
S\text{ is a true martingale on }[0,T]
\quad\Longleftrightarrow\quad v_{gA,B}(\infty)=\infty,
\]
where, with $C(r)=\exp(\int_{r_0}^r B(u)\,du)$,
\[
v_{gA,B}(\infty)
=\int_{r_0}^\infty\frac1{C(r)}
  \int_{r_0}^r\frac{C(s)}{g(s)A(s)}\,ds\,dr.
\]
\end{corollary}

\begin{proof}
Set $\widehat A(r)=g(r)A(r)$ and $\widehat q(y)=q(y)/g(P(y))$ on $\{P\ge r_0\}$. The function $\widehat A$ is continuous and strictly positive, and \eqref{eq:variable-clock} gives $q_-\le\widehat q\le q_+$. Moreover,
\[
\ip{a\nabla P}{\nabla P}=\widehat q(y)\widehat A(P(y)),
\qquad 2\Lt P=\widehat q(y)\widehat A(P(y))B(P(y)).
\]
All hypotheses of Proposition~\ref{prop:adapted-infinity} hold with $\widehat A,\widehat q$ in place of $A,q$. That proposition yields the equivalence for every starting point and every finite horizon. Substituting $\widehat A=gA$ in \eqref{eq:feller} gives the displayed integral.
\end{proof}

The factor $g$ changes the speed density of the scalar reference diffusion from $2C/A$ to $2C/(gA)$, while its scale density remains $1/C$. Thus growth or decay of the rate is included in the Feller integral; the bounds in \eqref{eq:variable-clock} control its remaining variation along each level set of $P$.

\subsection{Comparison with existing martingality criteria}\label{sec:LMDP}

Criens--Glau \cite[Cor.~5.2]{CG18} give a necessary-and-sufficient martingality criterion in a radial multidimensional setting. In the factor-driven case $\theta=\sigma^{\top}\beta$, their Standing Assumption~5.2 requires the three functions
\[
\ip{y}{a(y)y},\qquad
\tfrac12\tr a(y)+\ip{y}{\mut(y)},\qquad
\ip{\beta(y)}{a(y)\beta(y)}=|\theta(y)|^2
\]
to depend only on $r=|y|^2/2$ throughout the state space $\R^d$. Under the condition that the auxiliary radial diffusion does not hit zero, they apply one-dimensional tests to its lifetime and to the integral of the quadratic-variation rate. Their original diffusion may be explosive; integrability at an accessible endpoint is therefore part of their criterion.

Proposition~\ref{thm:radial-iff} uses the same scalar reduction for the shifted squared radius. Under our standing assumption that the original diffusion is conservative, the question reduces to non-exit. The functions $\mathcal A,\mathcal B$ need be radial only outside a ball, and the criterion applies on an unbounded domain once exits through finite boundary points have been excluded. It imposes no radiality condition on $|\theta|^2$. This formulation separates the test at infinity from the local boundary conditions and permits the quadratic-variation rate of the exponential to vary along spheres. Proposition~\ref{prop:adapted-infinity} applies the same one-dimensional test through a function with compact sublevel sets and a common time-change rate bounded above and away from zero, as used for the anisotropic models of Example~\ref{ex:anisotropic}.

Criens--Glau \cite[Prop.~5.5]{CG18} also give comparison tests in a partially radial setting for global absolute continuity with the specified terminal density, and hence for uniform integrability. Criens \cite[Thm.~2.2, Cor.~2.8]{Cri21} characterizes uniform integrability by almost-sure explosion of a time-changed diffusion and gives Khasminskii-type tests in \cite[Sec.~3.1]{Cri21}. That time change uses the reciprocal of the quadratic-variation rate, assumed locally bounded away from zero. In Proposition~\ref{prop:adapted-infinity}, the two-sided bounds on $q$ instead preserve finiteness of the lifetime under a change of clock, giving a criterion for martingality on every finite horizon. Corollary~\ref{cor:variable-clock} incorporates a rate depending on the scalar coordinate into the Feller integral.

The sufficient conditions of \cite{LM07} and \cite{DP22} can now be placed within the exact criteria. In \cite{LM07} the volatility is one-dimensional: $dF_t=\sigma_tF_tdW_t$ with $d\sigma_t=\mu(\sigma_t)dZ_t+b(\sigma_t)dt$ on $[0,\infty)$, $\mu$ Lipschitz (their Section~2.5), or $dF_t=\sigma_t^\delta F_tdW_t$ with $d\sigma_t=\alpha\sigma_t^\gamma dZ_t+b(\sigma_t)dt$ (their Section~3), $\rho$ being the correlation of $W$ and $Z$. In the present notation $Y=\sigma$, the coefficient in $dF/F$ is $\sigma$ resp.\ $\sigma^\delta$, and the drift under the changed measure is $\mut(\xi)=\rho\mu(\xi)\xi+b(\xi)$, resp.\ $\rho\alpha\xi^{\gamma+\delta}+b(\xi)$. Their Theorem~2.4 (resp.\ Theorem~3.1, which further assumes $\rho>0$ and $\gamma + \delta >1$) states that $F$ is a martingale if $\limsup_{\xi\to\infty}\mut(\xi)/\xi<\infty$, and that $\E[F_t]<F_0$ for all $t>0$ if $\liminf_{\xi\to\infty}\mut(\xi)/\varphi(\xi)>0$ for some positive increasing $\varphi$ with $\int^\infty d\xi/\varphi<\infty$, the natural choice being $\varphi(\xi)=\xi^2$, resp.\ $\xi^{\gamma+\delta}$.

In \cite{DP22} the state process is $X=(M,v)$, where $M=(M^1,\dots,M^d)$ and $v$ solve
\begin{equation}\label{eq:DPmodel}
dM^j=M^j\sigma_{jj}(M,v)\,dB^j,\qquad dv=\bar\sigma(M,v)\,dZ+b(M,v)\,dt ,
\end{equation}
with $(B,Z)$ a correlated Brownian motion with correlation matrix $\Sigma$, and the question is which components $M^j$ are true martingales. This is the present setting with $Y=X$ on $D=(0,\infty)^d\times D_v$, $W$ a standard Brownian motion with $(B,Z)=\Sigma^{1/2}W$, and, for the component $j$, $\theta^{(j)}:=\sigma_{jj}(Y)\,(\Sigma^{1/2})^{\top}e_j$, so that $M^j=\cE(\int\theta^{(j)}(Y)\cdot dW)$ and $\sigma\theta^{(j)}=\beta^j$ is precisely the drift correction in \cite[(21)--(22)]{DP22}. Their Theorem~2 and Corollary~1, obtained via F\"ollmer's measure and Parthasarathy's extension theorem, are Lemma~\ref{thm:main}; their Theorems~3--4 are the two parts of Lemma~\ref{thm:lyap}; and their Theorems~6, 8 and~10 are Lemma~\ref{thm:radial-sandwich} applied to $\Lt^{(j)}$ for each $j$ (they write $A(\rho^2/2)$, $B(\rho^2/2)$ with $\rho=|y|$, which is the convention $r=\tfrac12|y|^2$ of \eqref{eq:AB}). The following proposition records these sufficient tests and their relation to the scalar and radial equivalences.

\begin{proposition}\label{prop:LMDP}
\begin{enumerate}[label=(\roman*),leftmargin=2em]
\item \emph{Scalar volatility} \cite[Thms.~2.4, 3.1]{LM07}. Let $d=1$, $D=(0,\infty)$, $dY=\mu(Y)dt+\sigma(Y)dZ$ with $\mu,\sigma$ locally Lipschitz, $\sigma>0$ on $(0,\infty)$ and $|\mu(\xi)|+|\sigma(\xi)|\le C(1+\xi)$, and assume that $0$ is not reached under $\Prob$. Let $S=\cE(\int Y^\delta\,dW)$ with $d\langle W,Z\rangle_t=\rho\,dt$ and $\delta>0$, so that $\mut(\xi)=\mu(\xi)+\rho\sigma(\xi)\xi^\delta$. Then:
\begin{enumerate}[label=(\alph*),leftmargin=2em]
\item $0$ is not reached under $\Qt$, and $S$ is a true martingale iff $+\infty$ is inaccessible for the scalar diffusion with drift $\mut$ and diffusion coefficient $\sigma$, i.e.\ iff the integral in Feller's test at $+\infty$ is infinite; this is Proposition~\ref{thm:radial-iff} with $d=1$, and the criterion of \cite{DS02,MU12,BKX12}.
\item If $\limsup_{\xi\to\infty}\mut(\xi)/\xi<\infty$, then $S$ is a true martingale.
\item If $\liminf_{\xi\to\infty}\mut(\xi)/\varphi(\xi)>0$ for a nondecreasing $\varphi\in C^1$ with $\int^\infty d\xi/\varphi<\infty$ and $\sigma(\xi)^2\le C\varphi(\xi)$ for large $\xi$, then $S$ is not a martingale on $[0,\infty)$; if moreover $\sigma,\mut$ are locally H\"older continuous, $\E[S_T]<1$ for every $T>0$.
\end{enumerate}
\item \emph{Vector models} \cite[Thms.~6, 8, 10]{DP22}. In the setting \eqref{eq:DPmodel}, assume Assumption~\ref{ass:standing} for $Y$ and for the generator $\Lt^{(j)}$ with drift $\mu+\sigma\theta^{(j)}$, and let $\mathcal A^{(j)},\mathcal B^{(j)}$ be the functions \eqref{eq:AB} of $\Lt^{(j)}$. Then:
\begin{enumerate}[label=(\alph*),leftmargin=2em]
\item If there are $A\ge\mathcal A^{(j)}$, $B\ge\mathcal B^{(j)}$ on spheres with $v_{A,B}(\infty)=\infty$, then $Y$ does not explode at infinity under $\Qt^{(j)}$; if moreover exits through $\partial D$ are excluded (Theorem~\ref{thm:fichera}, Proposition~\ref{prop:inherit}, Proposition~\ref{prop:sing}), $M^j$ is a true martingale.
\item If $D$ is unbounded and there are $A\le\mathcal A^{(j)}$, $B\le\mathcal B^{(j)}$ on spheres with $v_{A,B}(\infty)<\infty$, then $M^j$ is not a martingale on $[0,\infty)$; under local H\"older continuity, $\E[M^j_T]<M^j_0$ for every $T>0$ and every starting point.
\item If $D$ is unbounded, exits through $\partial D$ are excluded and $\mathcal A^{(j)},\mathcal B^{(j)}$ are radial outside a ball, then $M^j$ is a true martingale iff $v_{A,B}(\infty)=\infty$: the two alternatives of \cite[Thm.~10]{DP22} are then exhaustive and that theorem becomes an equivalence.
\end{enumerate}
In all cases what matters is the exit of the whole vector $Y$ under $\Qt^{(j)}$, through any component or through $\partial D$, not the explosion of the component $M^j$ alone (cf.\ \cite[Rem.~5]{DP22}).
\end{enumerate}
\end{proposition}

Statement (i) extends \cite{LM07} to diffusion $Y$ with at most linear growth coefficients. Statements (i)(b) and (i)(c) are the conditions of \cite[Thm.~2.4]{LM07} (for $\delta=1$) and \cite[Thm.~3.1]{LM07} (for $\sigma=\alpha\xi^\gamma$, $\gamma\le1$): they are the linear-growth and the Osgood regimes of the integral in Feller's test in (i)(a). The side condition $\sigma^2\le C\varphi$ in (i)(c) holds for $\varphi(\xi)=\xi^2$, the choice made in \cite{LM07}, and in the power model when $\delta\ge\gamma$; without it, (i)(c) still follows from (i)(a) or from the comparison argument of \cite{LM07}.

\begin{proof}
(i)(a) The coefficient $\xi^\delta$ in $dS/S$ is bounded near $0$, so Proposition~\ref{prop:inherit} transfers the non-attainability of $0$ from $\Prob$ to $\Qt$. In $d=1$ the sets $\{y\in D:\tfrac12y^2=r\}$ are single points, so $\mathcal A(\xi)=\sigma(\xi)^2\xi^2$ and $\mathcal B(\xi)=(\sigma(\xi)^2+2\xi\mut(\xi))/(\sigma(\xi)^2\xi^2)$ are trivially radial, $D$ is unbounded, and Proposition~\ref{thm:radial-iff} gives: $S$ is a true martingale iff $v_{A,B}(\infty)=\infty$. The substitution $r=\xi^2/2$ (as in the proof of Proposition~\ref{cor:osgood}) turns $v_{A,B}(\infty)$ into the integral in Feller's test $\int^\infty s'(\eta)\int^\eta m(\xi)\,d\xi\,d\eta$ of the scalar diffusion, with $s'=\exp(-\int2\mut/\sigma^2)$ and $m=2/(s'\sigma^2)$. (b) $\tr a+2\ip{y}{\mut}=\sigma(\xi)^2+2\xi\mut(\xi)\le C(1+\xi^2)$, so Proposition~\ref{cor:osgood}(i) excludes explosion at infinity, and (a) concludes. (c) Apply Proposition~\ref{cor:osgood}(ii) with $w(\xi)=\sigma(\xi)^2\xi^2$, $c_1=C_1=1$, and $c\varphi$ in place of $\varphi$, where $\mut\ge c\varphi$ for large $\xi$: condition (c) there reads $\xi^2c\varphi(\xi)/w(\xi)=c\varphi/\sigma^2\ge c/C$. Hence $\Qt_y(\tau_D<\infty)>0$ for all $y$, $S$ is not a martingale on $[0,\infty)$ by Lemma~\ref{thm:main}, and Proposition~\ref{prop:finite}(ii) gives the finite-horizon statement.

(ii)(a) is Lemma~\ref{thm:radial-sandwich}(i) followed by Lemma~\ref{thm:main}; (b) is Lemma~\ref{thm:radial-sandwich}(ii), Lemma~\ref{thm:main}(ii) and Proposition~\ref{prop:finite}; (c) is Proposition~\ref{thm:radial-iff}. The last assertion is Lemma~\ref{thm:main}: $\E[M^j_T]=M^j_0\,\Qt^{(j)}_y(\tau_D>T)$, where $\tau_D$ is the exit time of the vector $Y$.
\end{proof}

For the log-normal model $dF=\sigma F\,dW$, $d\sigma=\alpha\sigma\,dZ$, part (i) recovers the basic dichotomy of \cite{LM07}: $F$ is a martingale iff $\rho\le0$. The examples of \cite{DP22} call for one comment on the comparison functions used there.

\begin{remark}[on {\cite[Examples~1--2]{DP22}}]\label{rem:DPex}
In the strict-local-martingale case of \cite[Example~1]{DP22},
\[
f=\frac{1}{2x_1^2}\Big(\frac{|x|}{2}\Big)^{2+\epsilon},\qquad
\sigma=\frac{1}{2x_2}\Big(\frac{|x|}{2}\Big)^{2+\epsilon},\qquad
b=\frac{|x|^{2+\epsilon}}{x_2},\qquad \rho=1,\quad \epsilon<1.
\]
The proposed comparison functions are $A(r)=r^{2+\epsilon}$ and $B(r)=1/r$ on $[3/2,\infty)$. For any $r_0\ge3/2$, these give
\[
v_{A,B}(\infty)
=\int_{r_0}^\infty\frac1r\int_{r_0}^r s^{-1-\epsilon}\,ds\,dr
=\infty,
\]
since the inner integral is bounded below by a positive constant for all sufficiently large $r$. Thus the explosion condition of their Theorem~8 is not met by these functions.

The coefficients give a sharper lower bound. Put $r=|x|^2/2$ and $g=(|x|/2)^{2+\epsilon}$. Under the auxiliary law $\Prob^1$ of \cite{DP22}, direct calculation gives
\begin{align*}
\ip{x}{ax}&=g^2=\Big(\frac r2\Big)^{2+\epsilon},\\
\frac{\tr a+2\ip{x}{\mut}}{\ip{x}{ax}}
&=\frac{5}{4x_1^2}+\frac1{4x_2^2}
+\frac{2^{4+3\epsilon/2}}{r^{1+\epsilon/2}}
\ \ge\ \frac{\beta_0}{r},
\qquad \beta_0:=\frac{3+\sqrt5}{4}>1.
\end{align*}
The lower bound follows from
$(5/x_1^2+1/x_2^2)(x_1^2+x_2^2)\ge(\sqrt5+1)^2$.
With $A(r)=(r/2)^{2+\epsilon}$ and $B(r)=\beta_0/r$, Tonelli's theorem gives the exact expression
\[
v_{A,B}(\infty)
=\int_{r_0}^\infty\frac{s^{\beta_0}}{A(s)}
  \int_s^\infty r^{-\beta_0}\,dr\,ds
=\frac{2^{2+\epsilon}}{\beta_0-1}\int_{r_0}^\infty s^{-1-\epsilon}\,ds.
\]
It is finite if and only if $\epsilon>0$. Hence the corrected lower comparison yields the stated strict-local-martingale conclusion for $\epsilon\in(0,1)$. For $\epsilon\le0$ it is inconclusive; the upper envelope is infinite near the axes. The same divergence for $B(r)=1/r$ occurs in the comparison functions proposed for explosion in \cite[Example~2]{DP22}. In Example~1, the choices $\rho=\pm1$ give a rank-one diffusion matrix. For $\rho=-1$ in its non-explosion case, $\ip{x}{ax}\equiv0$ and the squared radius has no martingale part. The coefficients are singular on the coordinate axes, so boundary behaviour must be checked separately from the radial comparison; the diffusion matrix also falls outside Assumption~\ref{ass:standing}(ii).
\end{remark}

\section{Singular boundary points}\label{sec:sing}

At intersections of boundary faces the smooth-boundary criterion does not apply. On the orthant, Proposition~\ref{prop:sing} states the classical sufficient non-attainment test obtained by projection onto the coordinates that vanish, following Ramasubramanian \cite{Ram88}. This test supplies the condition at nonsmooth boundary pieces in Theorem~\ref{thm:SV}(iv).

Proposition~\ref{thm:weighted-stratum} treats a separate class with a boundary piece of codimension at least two and an inverse-distance integrand modulated by correlated factors. For this class, both directions of the martingality criterion are proved, including critical equality. The scalar threshold is the Bessel threshold; the necessity argument also controls the singular drift induced in the tangential factors. Example~\ref{ex:weighted-stratum} shows concretely that a weighted distance gives an exact answer when the Euclidean projection test is inconclusive.

Let $D=(0,\infty)^d$. A boundary stratum consists of points with a specified set of zero coordinates and all remaining coordinates positive. The open faces have one zero coordinate; intersections of two or more faces form the nonsmooth part of the boundary. The estimates below apply to closed intersections of faces. Fix one such intersection, which after relabelling is $\Gamma_k:=\{y\in\bar D:y_1=\dots=y_{d-k}=0\}$, $0\le k\le d-2$, and let $\pi(y):=(y_1,\dots,y_{d-k})$. Following \cite{Ram88}, set
\[
A_\pi(y):=\sum_{i,j\le d-k}a_{ij}(y)\frac{y_iy_j}{|\pi(y)|^2},\qquad B_\pi(y):=\sum_{i\le d-k}a_{ii}(y),\qquad C_\pi(y):=2\sum_{i\le d-k}y_i\mut_i(y),
\]
\[
\beta(r):=\inf\Big\{\frac{B_\pi(y)-A_\pi(y)+C_\pi(y)}{A_\pi(y)}:\ y\in D,\ |\pi(y)|=r\Big\}.
\]
$A_\pi$ and $(B_\pi-A_\pi+C_\pi)/(2|\pi|)$ are the diffusion and drift rates of the distance $|\pi(Y)|$ to $\Gamma_k$, and $\beta$ plays the role of ``dimension minus one'' of a Bessel process: for $a=I$ and $\mut=0$ one gets $\beta\equiv d-k-1$.

\begin{proposition}\label{prop:sing}
Let Assumption~\ref{ass:standing} hold, let $r_0>0$, assume that $\beta$ is measurable and locally bounded on $(0,r_0]$, and that
\begin{equation}\label{eq:betacond}
\int_0^{r_0}\exp\Big(\int_u^{r_0}\frac{\beta(s)}{s}\,ds\Big)du=\infty\qquad(\text{e.g.\ }\beta(s)\ge1\text{ for all small }s).
\end{equation}
Then $\Qt_y\big(\tau_D<\infty,\ |\pi(Y_t)|\to0\text{ as }t\uparrow\tau_D\big)=0$ for all $y\in D$; in particular the diffusion cannot exit through $\Gamma_k$.
\end{proposition}

\begin{proof}
For $0<c<r_0$ let $F(r)=F(r;c):=\int_c^r\exp(-\int_c^u\beta(s)s^{-1}ds)du$, $c\le r\le r_0$. Then $F\in C^{1,1}([c,r_0])$, $F'\ge0$, $F(c)=0$, and $F''+(\beta/r)F'=0$ a.e. By the Meyer--It\^o formula for the semimartingale $|\pi(Y)|$ (whose quadratic variation rate is $A_\pi$ and whose drift rate is $(B_\pi-A_\pi+C_\pi)/(2|\pi|)$, by \eqref{eq:radialgen}-type calculus with $p(y)=|\pi(y)|$), on $\{c<|\pi(Y)|<r_0\}$,
\begin{align*}
d F(|\pi(Y_t)|)&=\Big(\tfrac12A_\pi F''+\frac{F'}{2|\pi|}\big(B_\pi-A_\pi+C_\pi\big)\Big)(Y_t)\,dt+dM_t\\
&\ge\ \tfrac12A_\pi\Big[F''+\frac{\beta(|\pi|)}{|\pi|}F'\Big](Y_t)\,dt+dM_t=dM_t ,
\end{align*}
so $F(|\pi(Y)|;c)$ is a bounded local submartingale up to $\zeta:=\tau_D\wedge\tau^{r_0}\wedge\tau_c$, where $\tau_c:=\inf\{t:|\pi(Y_t)|\le c\}$ and $\tau^{r_0}:=\inf\{t:|\pi(Y_t)|\ge r_0\}$. For $y$ with $c<|\pi(y)|<r_0$, optional sampling at $\zeta\wedge t$ and $t\to\infty$ give, since $F$ is increasing, $F=0$ at $c$ and $F\le F(r_0;c)$,
\[
F(|\pi(y)|;c)\ \le\ F(r_0;c)\,\Qt_y\big(\{\tau_D\wedge\tau^{r_0}<\tau_c\}\cup\{\zeta=\infty\}\big).
\]
With $g(u):=\exp(\int_u^{r_0}\beta(s)s^{-1}ds)$ one has
\[
\frac{F(r;c)}{F(r_0;c)}=\frac{\int_c^rg(u)\,du}{\int_c^{r_0}g(u)\,du}=1-\frac{\int_r^{r_0}g(u)\,du}{\int_c^{r_0}g(u)\,du}\ \longrightarrow\ 1\quad(c\downarrow0)
\]
exactly under \eqref{eq:betacond}. The events $\{\tau_D\wedge\tau^{r_0}<\tau_c\}\cup\{\zeta=\infty\}$ increase, as $c\downarrow0$, to $\{\tau_D\wedge\tau^{r_0}<\tau_0\}\cup\{\tau_0=\infty\}$ with $\tau_0:=\lim_c\tau_c$, so this set has full $\Qt_y$-measure: a.s., either $|\pi(Y)|$ never approaches $0$, or $Y$ leaves $D$ or the cylinder $\{|\pi|<r_0\}$ before doing so. On the event $\{\tau_D<\infty,\ |\pi(Y_t)|\to0\text{ as }t\uparrow\tau_D\}$ there is a rational $q<\tau_D$ with $|\pi(Y_t)|<r_0$ on $[q,\tau_D)$, and for the process restarted at time $q$ one has $\tau_0=\tau_D-q\le\tau_D\wedge\tau^{r_0}$; by the Markov property this event is therefore null.
\end{proof}

\begin{remark}
For $k=0$ (the origin) the infimum in $\beta$ is over a compact set and \eqref{eq:betacond} is often easy to check; e.g.\ for $a=I$, $\mut=0$ one recovers that Brownian motion does not hit the origin in $d\ge2$. For $k\ge1$ (axes in $d\ge3$) the infimum is over a cylinder and may be $\le0$ even in natural models (the three-dimensional volatility-stabilized model of Example~\ref{ex:VS} is such a case); this test alone does not give a necessary and sufficient condition for general orthant intersections. Proposition~\ref{thm:weighted-stratum} gives an exact weighted-distance criterion for a specified class of boundary pieces of positive dimension, and Remark~\ref{rem:VS-clock} excludes intersections at the first exit in Example~\ref{ex:VS}. In many stochastic volatility models the strata are excluded for structural reasons instead (Example~\ref{ex:twofactor}), and for affine models the proposition applies with $\beta(r)\ge1-Cr$ (Proposition~\ref{prop:affine}).
\end{remark}

\begin{example}[a two-dimensional diffusion on the orthant]\label{ex:VS}
Consider on $D=(0,\infty)^2$ the diffusion
\begin{equation}\label{eq:VS}
dY^1_t=Y^1_tdt+\sqrt{Y^1_t(Y^1_t+Y^2_t)}\,dW^1_t,\qquad dY^2_t=Y^2_tdt+\sqrt{Y^2_t(Y^1_t+Y^2_t)}\,dW^2_t ,
\end{equation}
the auxiliary diffusion of the volatility-stabilized model of Fernholz--Karatzas \cite[Sec.~12.1]{FK10}, taken here only as a test case for the criteria of Sections~\ref{sec:boundary} and \ref{sec:sing}. Here $a=\operatorname{diag}\big(y_1(y_1+y_2),\,y_2(y_1+y_2)\big)$ and $\mut(y)=y$. For the origin ($k=0$),
\[
A_\pi=\frac{(y_1+y_2)(y_1^3+y_2^3)}{y_1^2+y_2^2},\qquad B_\pi=(y_1+y_2)^2,\qquad C_\pi=2(y_1^2+y_2^2),
\]
and a direct computation gives
\[
\frac{B_\pi-A_\pi+C_\pi}{A_\pi}=1+\frac{y_1^4+6y_1^2y_2^2+y_2^4}{(y_1+y_2)^2(y_1^2-y_1y_2+y_2^2)} ,
\]
a function homogeneous of degree $0$ whose second term extends continuously to the closed quarter circle and is positive there (its minimum is $\approx0.96$, attained near the axes; it equals $1$ on the axes and $2$ on the diagonal), so that the infimum over the open quarter circle is controlled by the minimum over its compact closure. Hence $\beta\equiv\beta_0\ge1.96>1$, \eqref{eq:betacond} holds, and by Proposition~\ref{prop:sing} the process cannot reach the origin before reaching one of the open faces. On the faces, the Fichera drift is attractive: on $\{y_1=0\}$, with $\nu=-e_1$,
\[
F(0,y_2)=-\Big(\mut_1-\tfrac12\partial_1a_{11}\Big)\Big|_{y_1=0}=-\Big(0-\tfrac{y_2}2\Big)=\frac{y_2}{2}>0 ,
\]
so $\{y_1=0\}\subset\Sig_2$ (and likewise $\{y_2=0\}$), and Proposition~\ref{thm:attain}(i) shows that the process leaves $D$ with positive probability, necessarily through a face. In fact it leaves a.s.: the ratio $\nu_t:=Y^1_t/(Y^1_t+Y^2_t)$ satisfies, by It\^o's formula, $d\nu_t=\sqrt{\nu_t(1-\nu_t)}\,d\widetilde W_t$ for a Brownian motion $\widetilde W$ (all drift terms cancel), a driftless Wright--Fisher diffusion on $(0,1)$ for which both endpoints are attainable by Feller's test; being a bounded martingale with $\int_0^\infty\nu_t(1-\nu_t)dt<\infty$, it reaches $\{0,1\}$ in finite time a.s., i.e.\ one of the coordinates of $Y$ hits $0$ a.s. In dimension three the stratum $\{y_1=y_2=0\}$ ($k=1$) has $\beta\le0$, because the infimum over the cylinder $\{y_1^2+y_2^2=r^2\}$ is approached as $y_3\to\infty$, and Proposition~\ref{prop:sing} gives no information there. Remark~\ref{rem:VS-clock} resolves this case through the classical squared-Bessel representation: in every dimension the first exit occurs almost surely through a single open face. Thus the inconclusive projection bound is a limitation of that test, rather than an unresolved boundary classification for this particular model.
\end{example}

\begin{remark}[first exit in Example~\ref{ex:VS}]\label{rem:VS-clock}
The failure of the projection test at the axes of the three-dimensional example can also be resolved structurally. In any dimension $d\ge2$, consider
\[
dY_i=Y_i\,dt+\sqrt{Y_i\sum_jY_j}\,dW_i
\]
with positive initial coordinates, independent drivers and absorption of each coordinate at $0$. Set $X_i(t)=e^{-t}Y_i(t)$. A common time change of independent squared Bessel processes, as in \cite[Eq.~(4)]{Pal11}, gives
\[
X_i(t)=R_i(A_t),\qquad
A_t=\frac14\int_0^t\sum_jX_j(s)\,ds,
\]
where the $R_i$ are independent squared Bessel processes of dimension zero. Their absorption times $T_i$ are independent and have continuous distributions
$\Prob(T_i\le u)=\exp(-Y_i(0)/(2u))$, $u>0$. Thus $T_*:=\min_iT_i$ has a unique minimizing index almost surely. Since $\sum_iR_i$ stays positive on the compact interval $[0,T_*]$, the inverse time change
\[
t(T_*)=4\int_0^{T_*}\frac{du}{\sum_iR_i(u)}
\]
is finite. Hence the first exit from the open orthant occurs almost surely through an open face: exactly one coordinate vanishes. Intersections of two or more faces are not possible locations of the first exit in this model, including the axes at which the projection test is inconclusive.
\end{remark}
\subsection{Weighted distance to a boundary piece}\label{sec:weighted-stratum}

The next class has a noncompact boundary piece of codimension $m\ge2$. Its normal distance is exactly reducible to a Bessel comparison, even with correlated tangential factors.

\begin{proposition}\label{thm:weighted-stratum}
Let $m\ge2$, $k\ge1$, and let $(W,V)$ be a Brownian motion in $\R^m\times\R^k$ with covariance
\[
\begin{pmatrix}I_m&C\\C^\top&I_k\end{pmatrix}\succ0.
\]
Let $U_t=u+W_t$, $u\ne0$, and $dZ_t=-\Lambda Z_t\,dt+H\,dV_t$, where $\Lambda$ and $H$ are positive diagonal matrices. For a bounded, locally H\"older function $c:\R^k\to[0,\infty)$, define
\[
S=\cE\left(-\int c(Z_t)\frac{U_t}{|U_t|^2}\cdot dW_t\right).
\]
For every $T>0$ and every $(u,z)\in(\R^m\setminus\{0\})\times\R^k$,
\begin{equation}\label{eq:stratum-iff}
S\text{ is a true martingale on }[0,T]
\quad\Longleftrightarrow\quad \sup_{z\in\R^k}c(z)\le\frac{m-2}{2}.
\end{equation}
If the inequality fails, $\E_{u,z}[S_T]<1$ for every $T>0$ and every starting point.
\end{proposition}

\begin{proof}
Brownian motion in dimension $m\ge2$ does not hit the origin, so the original diffusion remains in its domain. The shifted system is
\[
dU_t=-c(Z_t)\frac{U_t}{|U_t|^2}\,dt+d\widehat W_t,
\qquad
dZ_t=\left(-\Lambda Z_t-HC^\top c(Z_t)\frac{U_t}{|U_t|^2}\right)dt
+H\,d\widehat V_t.
\]
For $R=|U|$, It\^o's formula gives
\[
dR_t=d\beta_t+\frac{m-1-2c(Z_t)}{2R_t}\,dt,
\]
where $\beta$ is a Brownian motion up to the lifetime. If $\bar c:=\sup c\le(m-2)/2$, stop $R$ in $(1/n,n)$. The stopped integrand is bounded, and under the corresponding probability measure $R$ lies between Bessel processes of dimensions $m-2\bar c\ge2$ and $m$. Their exit probabilities tend to zero, proving martingality by the same stopped-density argument as above.

If the threshold fails, choose a tangential box $K$ on which $c\ge c_*>(m-2)/2$, and choose $0<p<\min\{1,2c_*-m+2\}$. Then
\[
\Lt R^p=\frac p2R^{p-2}(p+m-2-2c(z))\le-\delta R^{p-2}
\]
on $0<R<\varepsilon$, $z\in K$. The weight $R^{p-2}$ dominates $1+R^{-1}$, which bounds both the normal and the tangential drift on this cylinder. The stopped It\^o and maximal estimates in the proof of Theorem~\ref{thm:coupled-log}, with $f(R)=R^p$, show that for small initial $R$ the process reaches $R=0$ in finite time with positive probability before a tangential exit. All coordinates have finite limits by the occupation bound. Lemma~\ref{thm:main} and Proposition~\ref{prop:finite}(ii) complete the proof.
\end{proof}

The scalar calculation is the classical Bessel boundary test. A change of coordinates makes its usefulness beyond Euclidean projection explicit. Let $M:\R^k\to\mathrm{GL}_m(\R)$ be smooth and set $X=M(Z)U$. The joint diffusion $(X,Z)$ is defined by It\^o's formula for this change of coordinates, and
\begin{equation}\label{eq:weighted-distance}
r(x,z)=|M(z)^{-1}x|,
\qquad
r(x,z)^2=x^\top M(z)^{-\top}M(z)^{-1}x.
\end{equation}
The integrand against $W$ is $-c(z)M(z)^{-1}x/r(x,z)^2$, so the same stochastic exponential satisfies \eqref{eq:stratum-iff}. This supplies weights depending on the tangential factors. All derivative and cross-variation terms in the coefficients of $X$ are included by the coordinate transformation. Since $M$ and $M^{-1}$ are bounded on compact subsets of $\R^k$, the boundary piece is $\{x=0\}$ and conservativeness of the original system is preserved.

\begin{example}[Euclidean projection can miss nonattainment]\label{ex:weighted-stratum}
Take $m=3$, $k=1$, $M=\operatorname{diag}(2,1,1)$ and
\[
c(z)=\frac15+\frac15\frac{z^2}{1+z^2}.
\]
Then $\sup c=2/5<1/2$, and Proposition~\ref{thm:weighted-stratum} gives a true martingale for every admissible $C$. In the $x$ coordinates the normal covariance is $\operatorname{diag}(4,1,1)$ and the shifted normal drift is $-c(z)x/|M^{-1}x|^2$. Along the first normal axis, the quantities used by Euclidean projection are
\[
A_\pi=4,\qquad B_\pi=6,\qquad C_\pi=-8c(z),\qquad
\frac{B_\pi-A_\pi+C_\pi}{A_\pi}=\frac12-2c(z).
\]
Thus the infimum over every Euclidean cylinder is at most $-3/10$. The integral in \eqref{eq:betacond} is then finite, so that sufficient test is inconclusive, while the weighted distance \eqref{eq:weighted-distance} gives the exact answer. In this domain, $(\R^3\setminus\{0\})\times\R$, the excluded boundary piece has dimension one and codimension three.
\end{example}

\section{Applications to multi-factor stochastic volatility models}\label{sec:SV}

We now return to the model $dS_t=S_t\,\theta(Y_t)\cdot dW_t$, $dY_t=\mu(Y_t)dt+\sigma(Y_t)dW_t$, and collect the criteria in the form in which they are applied. Theorem~\ref{thm:SV} is the summary of Sections~\ref{sec:exact}--\ref{sec:sing}; the examples that follow are the point of the exercise. The power-volatility examples are decided at infinity. The positive square-root-factor models have coefficients $\theta$ singular at one face and closed-form parameter thresholds, including correlated modulation. Section~\ref{sec:absorbed-model} treats a density that vanishes when an absorbed factor reaches zero.

\begin{theorem}\label{thm:SV}
Let Assumption~\ref{ass:standing} hold.
\begin{enumerate}[label=(\roman*),leftmargin=2em]
\item $S$ is a true martingale on $[0,T]$ iff the diffusion with generator $\Lt$ (drift $\mu+\sigma\theta$) started at $y$ stays in $D$ up to time $T$, a.s.; on $[0,\infty)$ iff it never leaves $D$, for one or equivalently all $y$. If, additionally, $a$ and $\tilde{\mu}$ are locally H\"{o}lder continuous, then $S$ is a true martingale on $[0,\infty)$ iff $\Lt u=\lambda u$ has no bounded positive solution on $D$.
\item Suppose $\partial D\in C^2$, $a\in C^1(\bar D)$, $\mu,\mut\in C(\bar D)$, and that \eqref{eq:lingrowth} holds or, more generally, that Lemma~\ref{thm:radial-sandwich}(i) applies. Then $S$ is a true martingale if $\ip{\mut-\tfrac12\dv a}{\nu}<0$ on $\partial D$, and only if $\ip{\mut-\tfrac12\dv a}{\nu}\le0$ on $\partial D$; the two coincide (an ``iff'') when (T) holds at the characteristic points with vanishing Fichera drift near which $\theta$ is unbounded.
\item Suppose $D=\R^d$ (or $D$ unbounded with exits through $\partial D$ excluded by (ii) and Proposition~\ref{prop:sing}) and $\mathcal A,\mathcal B$ in \eqref{eq:AB} are radial outside a ball. Then $S$ is a true martingale iff $v_{A,B}(\infty)=\infty$.
\item Suppose $D=(0,\infty)^d$, $a\in C^1(\bar D)$ and $\mu,\mut\in C(\bar D)$ (so that the coefficients are bounded on bounded subsets of $\bar D$, corners included, and every bounded path with $\tau_D<\infty$ has a limit in $\partial D$), \eqref{eq:lingrowth} holds, and every stratum satisfies \eqref{eq:betacond} (or is excluded by a structural argument). Then $S$ is a true martingale iff $\ip{\mut-\tfrac12\dv a}{\nu}\le0$ on every face, provided (T) holds at the face points with vanishing Fichera drift near which $\theta$ is unbounded.
\end{enumerate}
\end{theorem}

\begin{proof}
(i) is Lemmas~\ref{thm:main} and \ref{thm:khas}; (ii) is Theorem~\ref{thm:fichera}; (iii) is Proposition~\ref{thm:radial-iff}; (iv) combines Theorem~\ref{thm:fichera} (whose covering argument applies to the faces, the $C^2$ part of the boundary) with Proposition~\ref{prop:sing} for the strata: under \eqref{eq:lingrowth} the path is bounded on $\{\tau_D<\infty\}$, the limit point $Y_{\tau_D-}$ exists because the coefficients are bounded on bounded subsets of $\bar D$, and it lies either on a face or on a stratum; both events have probability zero. For the necessity direction use Proposition~\ref{thm:attain}(i) on a face.
\end{proof}

The first two examples are decided at infinity; the subsequent models are decided at a finite boundary.

\begin{example}[isotropic multi-factor model with power volatility]\label{ex:iso}
Let $D=\R^d$, $d\ge2$, $Z$ a $d$-dimensional Brownian motion, $\kappa,\eta>0$, $p>0$, $\rho\in[-1,1]$, and
\[
dY_t=-\kappa Y_tdt+\eta\,dZ_t,\qquad dS_t=S_t|Y_t|^p\,dB_t,\qquad dB_t=\rho\Big\langle\frac{Y_t}{|Y_t|},dZ_t\Big\rangle+\sqrt{1-\rho^2}\,dB^\perp_t ,
\]
with $B^\perp$ independent of $Z$ ($B$ is a Brownian motion by L\'evy's characterisation; $Y$ never hits $0$ since $d\ge2$). Here $W=(Z,B^\perp)$, $\sigma=[\eta I_d\ 0]$, $\theta(y)=(\rho|y|^{p-1}y,\sqrt{1-\rho^2}|y|^p)$ with $\theta(0):=0$, so $\sigma\theta=\rho\eta|y|^{p-1}y$ and
\[
\mut(y)=\big(-\kappa+\rho\eta|y|^{p-1}\big)y,\quad \mathcal A(y)=\eta^2|y|^2=2\eta^2r,\quad \mathcal B(y)=\frac{d\eta^2-4\kappa r+2\rho\eta(2r)^{(p+1)/2}}{2\eta^2r}=:B(r).
\]
Both functions depend only on the radius, so Theorem~\ref{thm:SV}(iii) applies and
\[
S\ \text{is a true martingale}\iff\rho\le0\ \text{ or }\ p\le1 .
\]
Indeed, if $\rho>0$ and $p>1$ then $B(u)\ge cu^{(p-1)/2}$ for $u\ge r_1$, so $\int_s^\rho B\ge cs^{(p-1)/2}(\rho-s)$ and
\[
v_{A,B}(\infty)=\int_{r_1}^\infty\frac{1}{2\eta^2s}\int_s^\infty e^{-\int_s^\rho B(u)du}\,d\rho\,ds\le\int_{r_1}^\infty\frac{ds}{2\eta^2c\,s^{(p+1)/2}}<\infty .
\]
If $\rho\le0$ or $p\le1$, then $B\le\bar b$ on $[r_0,\infty)$ for some constant $\bar b\ge1$, and the same double integral is $\ge\int_{r_0}^\infty\frac{ds}{2\eta^2\bar b\,s}=\infty$. For $d=1$ and $p=1$ this is the Stein--Stein/Sch\"obel--Zhu model. The dichotomy ``$\rho>0$ and superlinear volatility $\Rightarrow$ strict local martingale'' is the additive-noise, multi-factor counterpart of \cite[Thm.~3.1]{LM07} (whose condition for the model $dF=\sigma^\delta F\,dW$, $d\sigma=\alpha\sigma^\gamma dZ+b(\sigma)dt$ reads ``$\rho>0$ and $\gamma+\delta>1$''; here $\gamma=0$, $\delta=p$), and of \cite{Jou04,AP07}, all of which concern one volatility factor. The multi-factor model and its characterization are new; for $p>1$ the model is not affine and the Riccati approach does not apply.
\end{example}

\begin{example}[anisotropy and tangential rotation]\label{ex:anisotropic}
Let $d\ge2$, let $M$ be invertible, and write $z=M^{-1}y$. Let $q:\R^d\to(0,\infty)$ be smooth and bounded above and away from zero, and let $\Omega:\R^d\to\R^{d\times d}$ be smooth with $\Omega(z)^\top=-\Omega(z)$. With $W$ a standard $d$-dimensional Brownian motion, take
\[
dY=q(Y)M[-\kappa z+\Omega(z)z]dt+\eta\sqrt{q(Y)}M\,dW,
\qquad \kappa,\eta>0.
\]
For an independent scalar Brownian motion $B^\perp$, parameters $\gamma>0$ and $\rho\in[-1,1]$, define the coefficient of $(dW,dB^\perp)$ in $dS/S$ by
\[
\theta(y)=\sqrt{q(y)}\big(\rho|z|^{\gamma-1}z,\ \sqrt{1-\rho^2}|z|^\gamma\big),
\qquad\theta(0)=0.
\]
Set $P(y)=|M^{-1}y|^2/2$. Since $z^\top\Omega(z)z=0$, the original generator satisfies $LP=q(y)(d\eta^2/2-2\kappa P)$, proving conservativeness. Under the shifted generator,
\begin{align*}
\ip{a\nabla P}{\nabla P}&=2\eta^2q(y)P,\\
2\Lt P&=q(y)\big[d\eta^2-4\kappa P+2\rho\eta(2P)^{(\gamma+1)/2}\big].
\end{align*}
Thus \eqref{eq:adapted-coordinate} holds with
\[
A(r)=2\eta^2r,\qquad
B(r)=\frac{d\eta^2-4\kappa r+2\rho\eta(2r)^{(\gamma+1)/2}}{2\eta^2r}.
\]
These are the scalar functions of Example~\ref{ex:iso}. Consequently
\begin{equation}\label{eq:anisotropic-iff}
S\text{ is a true martingale on every finite horizon}
\quad\Longleftrightarrow\quad \rho\le0\ \text{or}\ \gamma\le1.
\end{equation}
For instance, $M=\operatorname{diag}(2,1,\ldots,1)$, $q(y)=1+\tfrac12\sin y_1$ and $\Omega(z)=|z|^2J$, with a nonzero constant skew-symmetric matrix $J$, give anisotropic covariance, a time-change rate that varies along the ellipsoids, and a cubic tangential drift. Under the stated hypotheses, the martingality threshold \eqref{eq:anisotropic-iff} is independent of $M$, $q$ and $\Omega$.
\end{example}
\begin{example}[the $3/2$ model recovered from the Fichera drift]\label{ex:32}
Let $d=1$, $D=(0,\infty)$, $dY=(\alpha-\kappa Y)dt+\eta\sqrt Y\,dZ$ with $2\alpha\ge\eta^2$ (so that Assumption~\ref{ass:standing}(iii) holds), and let the asset have volatility $Y^{-1/2}$, i.e.\ $\theta(y)=y^{-1/2}(\rho,\sqrt{1-\rho^2})$ with $W=(Z,B^\perp)$. Then $V:=1/Y$ is the $3/2$ process $dV=V\big(\kappa-(\alpha-\eta^2)V\big)dt-\eta V^{3/2}dZ$ and $S$ has volatility $\sqrt V$; this is the $3/2$ model of \cite{CS07}. Here $\sigma\theta=\rho\eta$ is constant, $\mut(y)=\alpha+\rho\eta-\kappa y$, and at the boundary point $0$ (outward normal $\nu=-1$)
\[
\Ft(0)=-\big(\alpha+\rho\eta-\tfrac12\eta^2\big),\qquad \frac{a(y)}{y}-2\mut(y)=\eta^2-2\alpha-2\rho\eta+2\kappa y ,
\]
so (T) holds on $\{\Ft=0\}$ and the linear growth of $\mut$ gives non-explosion at infinity. Theorem~\ref{thm:SV}(ii) yields
\[
S\ \text{is a true martingale}\iff\alpha+\rho\eta\ge\tfrac12\eta^2 ,
\]
which is the known criterion for the $3/2$ model (in the $V$-parametrisation $dV=\kappa'V(\theta'-V)dt+\epsilon V^{3/2}dZ'$ with correlation $\rho'$ between $Z'$ and $B$: $\kappa'-\rho'\epsilon\ge-\epsilon^2/2$), obtained there from Feller's test at $V=\infty$; the boundary derivation is the starting point for the coupled models below. Note that $\theta$ is unbounded at $0$, in accordance with Corollary~\ref{cor:inherit}.
\end{example}

\begin{example}[a two-factor inverse-CIR model]\label{ex:twofactor}
Let $D=(0,\infty)^2$, $Z^1,Z^2,B^\perp$ independent Brownian motions,
\[
dY^i_t=(\alpha_i-\kappa_iY^i_t)dt+\eta_i\sqrt{Y^i_t}\,dZ^i_t\ \ (i=1,2),\qquad dS_t=S_t\,(Y^1_t)^{-1/2}g(Y^2_t)\,dB_t,
\]
where $B=\rho Z^1+\sqrt{1-\rho^2}B^\perp$, $2\alpha_i\ge\eta_i^2$ and $g:[0,\infty)\to(0,\infty)$ continuous and bounded, $\bar g:=\sup g$. (The variance $g(Y^2)^2/Y^1$ is a $3/2$-type variance modulated by a second factor.) Then $\sigma\theta=(\rho\eta_1g(y_2),0)$,
\[
\mut(y)=\big(\alpha_1+\rho\eta_1g(y_2)-\kappa_1y_1,\ \alpha_2-\kappa_2y_2\big)\in C(\bar D),
\]
and $\tr a+2\ip{y}{\mut}$ has linear growth. On the face $\{y_1=0\}$, $\nu=-e_1$ and
\[
\Ft(0,y_2)=-\big(\alpha_1-\tfrac12\eta_1^2+\rho\eta_1g(y_2)\big),\qquad \frac{a_{11}}{y_1}-2\mut_1=-2\big(\alpha_1-\tfrac12\eta_1^2+\rho\eta_1g(y_2)\big)+2\kappa_1y_1 ,
\]
so $\Ft\le0$ on this face iff $\alpha_1-\tfrac12\eta_1^2\ge\eta_1\rho^-\bar g$ ($\rho^-:=\max(-\rho,0)$), and then (T) holds on the characteristic boundary where the Fichera drift vanishes with $C=2\kappa_1$. On the face $\{y_2=0\}$, $\Ft=-(\alpha_2-\tfrac12\eta_2^2)\le0$ and (T) holds. Under $\Qt$ the second factor is an autonomous CIR process with $2\alpha_2\ge\eta_2^2$, so it never hits $0$ and neither the face $\{y_2=0\}$ nor the origin can be a limit point of $Y$ (alternatively, the computation in the proof of Proposition~\ref{prop:affine} gives $\beta(r)\ge1-Cr$ at the origin, using $2\alpha_1+2\rho\eta_1g\ge\eta_1^2$ and $2\alpha_2\ge\eta_2^2$, so Proposition~\ref{prop:sing} applies). Theorem~\ref{thm:SV}(iv) gives
\[
S\ \text{is a true martingale}\iff 2\alpha_1-\eta_1^2\ \ge\ 2\eta_1\,\rho^-\sup g .
\]
In particular $S$ is always a martingale when $\rho\ge0$, while for $\rho<0$ (positive correlation between the asset and its volatility) the second factor enters only through $\sup g$. For $g\equiv1$ this reduces to Example~\ref{ex:32}. For nonconstant $g$ and $\rho\ne0$, the first coordinate of the shifted system is not an autonomous scalar diffusion and its drift is non-affine. The inequality above therefore extends the direct scalar and affine criteria to this coupled model. The correlated extension and the comparison with other multidimensional criteria are given in Section~\ref{sec:correlated}.
\end{example}

\begin{example}[$4/2$-type models]\label{ex:42}
Grasselli's $4/2$ model \cite{Gra17} has $dS_t=S_t\big(a\sqrt{V_t}+b/\sqrt{V_t}\,\big)dB_t$, $dV_t=(\alpha-\kappa V_t)dt+\eta\sqrt{V_t}\,dZ_t$, $d\langle B,Z\rangle_t=\rho\,dt$, with $a,b\ge0$ ($b=0$ is Heston, $a=0$ is the $3/2$ model in the parametrisation of Example~\ref{ex:32}). The coefficient $\theta(v)=(a\sqrt v+b/\sqrt v)(\rho,\sqrt{1-\rho^2})$ is unbounded at $0$, the model is not affine (the term $\tfrac12|\theta|^2$ contains $b^2/(2v)$), but the Girsanov shift $\sigma\theta=\rho\eta(av+b)$ is affine, so
\[
\mut(v)=\alpha+\rho\eta b-(\kappa-\rho\eta a)v,\qquad \Ft(0)=-\big(\alpha+\rho\eta b-\tfrac12\eta^2\big),
\]
(T) holds on $\{\Ft=0\}$ as in Example~\ref{ex:32}, and $\mut$ has linear growth. Under $2\alpha\ge\eta^2$, Theorem~\ref{thm:SV}(ii) gives: $S$ is a true martingale iff $\alpha+\rho\eta b\ge\tfrac12\eta^2$, the condition obtained in \cite{Gra17} from the Bessel-process representation of $V$; the parameter $a$ (the Heston part) plays no role, in accordance with Corollary~\ref{cor:inherit}. A multi-factor version with a modulating factor,
\[
dS_t=S_t\Big(a\sqrt{Y^1_t}+\frac{b\,g(Y^2_t)}{\sqrt{Y^1_t}}\Big)dB_t,\qquad dY^i_t=(\alpha_i-\kappa_iY^i_t)dt+\eta_i\sqrt{Y^i_t}\,dZ^i_t\ \ (i=1,2),
\]
with $B=\rho Z^1+\sqrt{1-\rho^2}B^\perp$, $Z^1,Z^2,B^\perp$ independent, $2\alpha_i\ge\eta_i^2$ and $g$ continuous, bounded and positive, is treated exactly as Example~\ref{ex:twofactor}: $\mut_1(y)=\alpha_1+\rho\eta_1bg(y_2)-(\kappa_1-\rho\eta_1a)y_1$, the second factor is autonomous under $\Qt$, and
\[
S\ \text{is a true martingale}\iff 2\alpha_1-\eta_1^2\ \ge\ 2\eta_1\,(\rho b)^-\,\sup g .
\]
The one-factor criterion recovers \cite{Gra17}; its shifted factor remains CIR. For the modulated model, a nonconstant $g$ and $\rho b\ne0$ make the shifted factor system non-affine, and the boundary criterion gives the displayed exact inequality. Theorem~\ref{thm:correlated} extends this conclusion to correlated factors.
\end{example}

\subsection{Correlated factors and boundary singularities}\label{sec:correlated}

The next theorem extends Examples~\ref{ex:twofactor} and \ref{ex:42} to correlated square-root factors and a modulation depending on all the other factors. Correlation presents two obstacles to applying Theorem~\ref{thm:fichera}: the mixed covariance coefficients are not $C^1$ at the faces, and the shifted drifts of the other factors can be singular at $\{y_1=0\}$. We use scalar comparison for sufficiency and a local power barrier for necessity.

\begin{theorem}\label{thm:correlated}
Let $d\ge2$ and let $(B,Z^1,\ldots,Z^d)$ be a Brownian motion with correlation matrix
\[
\begin{pmatrix}1&\boldsymbol\rho^{\top}\\ \boldsymbol\rho&R\end{pmatrix}\succeq0,
\qquad R\succ0,\qquad R_{ii}=1,
\]
so that $d\langle Z^i,Z^j\rangle_t=R_{ij}\,dt$ and $d\langle B,Z^i\rangle_t=\rho_i\,dt$. Suppose
\begin{equation}\label{eq:correlated-CIR}
dY^i_t=(\alpha_i-\kappa_iY^i_t)dt+\eta_i\sqrt{Y^i_t}\,dZ^i_t,
\qquad \alpha_i,\kappa_i,\eta_i>0,\quad 2\alpha_i\ge\eta_i^2,
\end{equation}
with $Y_0=y\in(0,\infty)^d$. Let $u,v\ge0$ and let $g:(0,\infty)^{d-1}\to(0,\infty)$ be bounded and locally H\"older continuous. Write $\bar g:=\sup g$ and set
\begin{equation}\label{eq:correlated-loading}
\ell(y):=u\sqrt{y_1}+\frac{v\,g(y_2,\ldots,y_d)}{\sqrt{y_1}},
\qquad S:=\cE\Big(\int_0^\cdot\ell(Y_t)\,dB_t\Big).
\end{equation}
For every $T>0$ and every $y\in(0,\infty)^d$,
\begin{equation}\label{eq:correlated-iff}
S\text{ is a true martingale on }[0,T]
\quad\Longleftrightarrow\quad
2\alpha_1-\eta_1^2\ \ge\ 2\eta_1(\rho_1v)^-\bar g,
\end{equation}
where $q^\pm:=\max\{\pm q,0\}$. If the inequality fails, $\E_y[S_T]<1$ for every $T>0$ and every $y\in(0,\infty)^d$.
\end{theorem}

The factor correlations $R_{ij}$ and the asset correlations $\rho_i$, $i\ge2$, may take any values compatible with the displayed matrix conditions. They change the law and the shifted dynamics, while the threshold in \eqref{eq:correlated-iff} depends only on $\rho_1$. The equality case and a supremum $\bar g$ that is not attained are included. The result concerns a coefficient $\ell$ singular at one face, with bounded modulation by the remaining factors.

\begin{proof}
Each original factor is a strictly positive, non-explosive CIR process. On $D=(0,\infty)^d$ the covariance matrix
\[
a_{ij}(y)=R_{ij}\eta_i\eta_j\sqrt{y_iy_j}
\]
is smooth and locally uniformly positive definite. The shifted coefficients are locally H\"older continuous, and the shifted diffusion is well posed up to its lifetime (localize to compact subsets of $D$ and use Girsanov). Thus Lemma~\ref{thm:main} and Proposition~\ref{prop:finite} apply. Put $x=y_1$, $z=(y_2,\ldots,y_d)$ and
\[
c(z):=\alpha_1+\rho_1\eta_1v g(z),\qquad k:=\kappa_1-\rho_1\eta_1u.
\]
The shifted drift is
\begin{align}
\mut_1(x,z)&=c(z)-kx,\label{eq:correlated-shift}\\
\mut_i(x,z)&=\alpha_i-\kappa_i y_i+
\rho_i\eta_i\left(u\sqrt{xy_i}+v g(z)\sqrt{y_i/x}\right),
\qquad i\ge2.\nonumber
\end{align}

\emph{Sufficiency.} Suppose the inequality in \eqref{eq:correlated-iff} holds, and define
\[
c_-:=\alpha_1-\eta_1(\rho_1v)^-\bar g\ge\tfrac12\eta_1^2,
\qquad c_+:=\alpha_1+\eta_1(\rho_1v)^+\bar g.
\]
Then $c_-\le c(z)\le c_+$. Fix $T>0$ and let
\[
\tau_n:=\inf\{t\ge0:Y^1_t\notin(1/n,n)\},\qquad
d\Q^{n,T}:=S_{T\wedge\tau_n}\,d\Prob.
\]
For sufficiently large $n$, the initial point lies in $(1/n,n)$. Before $\tau_n$ the coefficient $\ell(Y)$ is bounded by $\sqrt n(u+v\bar g)$, so Novikov's condition makes $\Q^{n,T}$ a probability measure. Under this measure, up to $\tau_n\wedge T$, the first factor has drift $c(Y^2,\ldots,Y^d)-kY^1$. Scalar comparison with the same Brownian driver gives
\[
X^-_t\le Y^1_t\le X^+_t\qquad(t\le\tau_n\wedge T),
\]
where $dX^\pm=(c_\pm-kX^\pm)dt+\eta_1\sqrt{X^\pm}\,d\widehat Z^1$ and $X^\pm_0=y_1$. The comparison follows from the usual Tanaka argument for the square-root coefficient; see \cite[Ch.~VI]{IW89}. Both comparison processes are strictly positive and non-explosive, including when $k\le0$, by the scalar CIR boundary test and linear growth. Their laws on $[0,T]$ are independent of $n$, hence
\[
\Q^{n,T}(\tau_n\le T)
\le \Q^{n,T}\Big(\inf_{t\le T}X^-_t\le1/n\Big)
+\Q^{n,T}\Big(\sup_{t\le T}X^+_t\ge n\Big)\longrightarrow0.
\]
Under $\Prob$, $\tau_n\uparrow\infty$. Consequently
\[
\E_y[S_T]=\lim_n\E_y[S_T\mathbf1_{\{\tau_n>T\}}]
=\lim_n\Q^{n,T}(\tau_n>T)=1.
\]

\emph{Necessity.} If the inequality in \eqref{eq:correlated-iff} fails, then $\rho_1v<0$ and there is an interior point $z^0$ with $c(z^0)<\eta_1^2/2$. Choose a closed box $K\subset(0,\infty)^{d-1}$ with $z^0\in\operatorname{int}K$ and
\[
\sup_{z\in K}c(z)<\tfrac12\eta_1^2.
\]
There exist $p\in(0,\tfrac12)$, $\varepsilon\in(0,1)$ and $\delta>0$ such that, on $U=(0,\varepsilon)\times\operatorname{int}K$,
\begin{equation}\label{eq:correlated-barrier}
\Lt x^p=p x^{p-1}\big(c(z)-kx-\tfrac12\eta_1^2(1-p)\big)
\le-\delta x^{p-1}.
\end{equation}
Start the shifted diffusion at $(x,z^0)\in U$. For $n>1/x$, let $\sigma_n$ be its exit time from $(1/n,\varepsilon)\times\operatorname{int}K$. Stopped It\^o's formula, first at $t\wedge\sigma_n$ and then with $t\to\infty$, gives
\begin{equation}\label{eq:correlated-occupation}
\E^{\Qt}_{(x,z^0)}[(Y^1_{\sigma_n})^p]
+\delta\,\E^{\Qt}_{(x,z^0)}\int_0^{\sigma_n}(Y^1_s)^{p-1}ds
\le x^p.
\end{equation}
In particular $\E^{\Qt}[\sigma_n]\le Cx^p$, with $C$ independent of $n$ and $x$. Write $Z'=(Y^2,\ldots,Y^d)$ for the tangential coordinates. On $U$, \eqref{eq:correlated-shift} gives
\[
|\mut_{2:d}(x,z)|\le C(1+x^{-1/2}),\qquad
\tr a_{2:d,2:d}(x,z)\le C.
\]
Because $p<1/2$ and $x<\varepsilon<1$, $1+x^{-1/2}\le2x^{p-1}$. Thus \eqref{eq:correlated-occupation} bounds the expected total variation of the tangential drift by $Cx^p$. If $M'$ is its martingale part, Doob's inequality also gives
\[
\E^{\Qt}\sup_{t\le\sigma_n}|M'_t|^2\le C\E^{\Qt}[\sigma_n]\le Cx^p.
\]
The positive distance from $z^0$ to $\partial K$ now yields
\[
\Qt_{(x,z^0)}(Z'_{\sigma_n}\in\partial K)\le Cx^p,
\qquad
\Qt_{(x,z^0)}(Y^1_{\sigma_n}=\varepsilon)\le(x/\varepsilon)^p.
\]
As $n\to\infty$, $\sigma_n\uparrow\sigma$ with $\E^{\Qt}[\sigma]\le Cx^p$. The occupation bound makes the tangential drift integrable up to $\sigma$; the bounded covariance gives a finite limit for its martingale part. The first coordinate also has a limit, since its drift and covariance are bounded on $U$. Therefore, except on an event of probability at most $Cx^p+(x/\varepsilon)^p$, the diffusion reaches $\{x=0\}$ in finite time with its other coordinates in $K$. This probability is positive for sufficiently small $x$. Lemma~\ref{thm:main} gives loss of mass, and Proposition~\ref{prop:finite}(ii) extends it to every starting point and every positive horizon.
\end{proof}

\begin{example}[a correlated model with an unattained supremum]\label{ex:correlated}
Take $d=2$, $\eta_1=\eta_2=\kappa_1=\kappa_2=1$, $\alpha_2=1$, $\alpha_1\ge1/2$, and
\[
d\langle Z^1,Z^2\rangle_t=\tfrac35dt,\qquad
B=-\tfrac12 Z^1+\tfrac{\sqrt3}{2}B^\perp,\qquad
g(z)=1+\frac{z}{1+z},\quad v=1,
\]
where $B^\perp$ is independent of $(Z^1,Z^2)$. Then $\rho_1=-1/2$, $\rho_2=-3/10$ and $\sup g=2$. For every $u\ge0$, Theorem~\ref{thm:correlated} gives
\[
S\text{ is a true martingale on every finite horizon}
\quad\Longleftrightarrow\quad \alpha_1\ge\tfrac32.
\]
Here $u=0$ gives inverse-CIR volatility and $u>0$ gives a modulated $4/2$ volatility coefficient. The second shifted drift is
\[
\mut_2(x,z)=1-z-\tfrac3{10}\left(u\sqrt{xz}+g(z)\sqrt{z/x}\right),
\]
which is singular as $x\downarrow0$. At $\alpha_1=3/2$ the infimum of $c(z)=\alpha_1-g(z)/2$ is approached only as $z\to\infty$, yet the equality case is a martingale. For $1/2\le\alpha_1<3/2$, a finite interval of sufficiently large $z$ has $c(z)<1/2$ and produces strict local martingality on every positive horizon. This tests sharpness without requiring the supremum of the modulation to be attained.
\end{example}

\paragraph{Comparison with existing criteria.}
Examples~\ref{ex:twofactor}, \ref{ex:42} and Theorem~\ref{thm:correlated} give exact conditions when $\theta$ is singular at the boundary and modulated by other factors. For nonconstant $g$ and $\rho_1v\ne0$, the first shifted coordinate is not an autonomous scalar diffusion. For non-diagonal $R$, the covariance matrix of the full factor system is also non-affine. Thus direct scalar Feller tests and affine Riccati criteria do not give \eqref{eq:correlated-iff}. The one-factor $4/2$ criterion is already in \cite{Gra17}; in the multifactor model of Gnoatto--Grasselli--Platen \cite[Sec.~3.2, eqs.~(3.7)--(3.9)]{GGP22}, the density separates across independent CIR factors and each shifted factor remains CIR. The modulation and factor correlations above remove this separation.

The radial tests of Dandapani--Protter \cite[Thms.~6, 8, 10]{DP22} give sufficient conditions for explosion or non-explosion at infinity. In Example~\ref{ex:twofactor}, the shifted drift has linear growth on both sides of the martingality threshold: the distinction is exit through a finite boundary point. Theorem~\ref{thm:correlated} decides this exit even with the singular tangential drift in \eqref{eq:correlated-shift}. Its proof uses classical comparison and stopping arguments, with \eqref{eq:correlated-occupation} controlling the additional singularity. This is the specific extension in applicability supplied by the correlated examples.

\begin{remark}[bounded integrands]\label{rem:correlated}
If $\theta$ is locally bounded on the closed orthant, Proposition~\ref{prop:inherit} transfers non-attainment of every finite boundary piece from $\Prob$ to $\Qt$, including intersections of faces. For correlated square-root factors and $\theta$ satisfying $|\theta(y)|\le C\sqrt{1+|y|}$, the shifted drift has linear growth, so explosion at infinity is excluded as well. Theorem~\ref{thm:correlated} treats the singular coefficient $\ell$ in \eqref{eq:correlated-loading} by a different localization: stopping the first coordinate makes the entire integrand bounded, uniformly in the other coordinates.
\end{remark}

\subsection{An absorbed factor with correlated modulation}\label{sec:absorbed-model}

For an absorbed factor, a singular coefficient $\theta$ can force the density to vanish at the boundary. The following model has an exact coefficient criterion with correlated modulation. The boundary continuation is part of its definition, as discussed in Section~\ref{sec:extended-state}.

\begin{theorem}\label{thm:absorbed}
Let $(W^0,W^1,\ldots,W^k)$, $k\ge1$, be a Brownian motion whose constant correlation matrix is positive definite. Set $\rho_i=d\langle W^0,W^i\rangle_t/dt$. For $x>0$, let
\[
X_t=x+W^0_{t\wedge\tau_0},\qquad
\tau_0:=\inf\{t\ge0:x+W^0_t=0\},\qquad
dZ^i_t=-\lambda_iZ^i_t\,dt+\eta_i\,dW^i_t,
\]
where $\lambda_i,\eta_i>0$. Let $c:\R^k\to(0,\infty)$ be bounded and locally H\"older continuous, and define
\begin{equation}\label{eq:absorbed-density}
S_t=\begin{cases}
\displaystyle\exp\left(\int_0^t\frac{c(Z_s)}{X_s}\,dW^0_s
-\frac12\int_0^t\frac{c(Z_s)^2}{X_s^2}\,ds\right),&t<\tau_0,\\[4pt]
0,&t\ge\tau_0.
\end{cases}
\end{equation}
Then $S$ is a continuous nonnegative local martingale, and for every $T>0$ and $(x,z)\in(0,\infty)\times\R^k$,
\begin{equation}\label{eq:absorbed-iff}
S\text{ is a true martingale on }[0,T]
\quad\Longleftrightarrow\quad \inf_{z\in\R^k}c(z)\ge\tfrac12.
\end{equation}
If the inequality fails, $\E_{x,z}[S_T]<1$ for every $T>0$ and every starting point.
\end{theorem}

\begin{proof}
First we justify the continuation by zero. It\^o's formula for $\log X$ shows that $\int_0^{\tau_0}X_s^{-2}ds=\infty$ on $\{\tau_0<\infty\}$: otherwise both terms in the formula would have finite limits. The path of $Z$ up to a finite $\tau_0$ is compact, and $c$ has a positive minimum on this path. Thus
\[
A_t:=\int_0^t\frac{c(Z_s)^2}{X_s^2}\,ds\longrightarrow\infty
\qquad(t\uparrow\tau_0<\infty).
\]
The time-change representation of a continuous local martingale gives $S_t=\exp(\beta_{A_t}-A_t/2)$ before $\tau_0$, for a Brownian motion $\beta$. Since $\beta_u/u\to0$, $S_t\to0$, and
\[
\int_0^{\tau_0}S_t^2\frac{c(Z_t)^2}{X_t^2}\,dt
=\int_0^\infty e^{2\beta_u-u}\,du<\infty
\qquad\text{on }\{\tau_0<\infty\}.
\]
Consequently \eqref{eq:absorbed-density} solves
$dS_t=S_t c(Z_t)X_t^{-1}\mathbf1_{\{t<\tau_0\}}dW^0_t$
on the whole time axis, with a locally square-integrable integrand. This proves the local martingale assertion.

Before $X$ reaches $0$, the shifted system is
\begin{equation}\label{eq:absorbed-shift}
dX_t=\frac{c(Z_t)}{X_t}\,dt+d\widehat W^0_t,
\qquad
dZ^i_t=\left(-\lambda_iZ^i_t+\frac{\eta_i\rho_i c(Z_t)}{X_t}\right)dt
+\eta_i\,d\widehat W^i_t.
\end{equation}
It is locally well posed on $D=(0,\infty)\times\R^k$ by stopped Girsanov. Exhaust $D$ by relatively compact sets and stop the density at their exit times. Under $\Prob$ these times increase to $\tau_0$; since $S_T=0$ on $\{\tau_0\le T\}$, monotone convergence in the stopped identity gives
\begin{equation}\label{eq:absorbed-mass}
\E_{x,z}[S_T]=\Qt_{x,z}(\tau_D>T),
\end{equation}
where the shifted process in \eqref{eq:absorbed-shift} is killed on leaving $D$.

If $\underline c:=\inf c\ge1/2$, stop only $X$ in $(1/n,n)$. The integrand $c(Z)/X$ is bounded on this stopped interval because $c$ is bounded. Under the resulting probability measures, scalar comparison places $X$ between Bessel processes of dimensions $1+2\underline c$ and $1+2\sup c$, both at least $2$. Their probabilities of leaving $(1/n,n)$ before $T$ tend to zero. The stopped-density argument of Theorem~\ref{thm:correlated}, with $S_T=0$ after $\tau_0$, therefore gives $\E_{x,z}[S_T]=1$.

Conversely, if $\inf c<1/2$, choose a closed box $K\subset\R^k$, an interior point $z^0$ and $c_*<1/2$ such that $c\le c_*$ on $K$. Take $0<p<\min\{1,1-2c_*\}$. For some $\delta>0$,
\[
\Lt x^p=p x^{p-2}\left(c(z)-\frac{1-p}{2}\right)
\le-\delta x^{p-2}\qquad\text{on }(0,\varepsilon)\times K.
\]
If $\sigma_n$ is the first exit from $(1/n,\varepsilon)\times\operatorname{int}K$, stopped It\^o gives
\[
\E^{\Qt}[X_{\sigma_n}^p]
+\delta\E^{\Qt}\int_0^{\sigma_n}X_s^{p-2}ds\le x^p.
\]
Here $1+x^{-1}\le Cx^{p-2}$ for $0<x<\varepsilon$. The estimate therefore bounds $\E^{\Qt}[\sigma_n]$ and the expected total variation of the tangential drift by $Cx^p$. Its martingale part has bounded quadratic-variation rate. Doob's inequality and the drift bound show that the probability of a tangential displacement reaching $\dist(z^0,\partial K)$ is at most $Cx^p$. The probability of exit at $x=\varepsilon$ is at most $(x/\varepsilon)^p$. For small $x$, with positive probability the limit $\sigma=\lim_n\sigma_n$ is finite, $Z$ remains in $\operatorname{int}K$, and $X\to0$. All coordinates have limits there, since the same occupation bound makes their drift integrable and their quadratic variation finite. This is an exit of the shifted system. Identity \eqref{eq:absorbed-mass} proves strict local martingality for some starting points and times. Local ellipticity and H\"older continuity on $D$ allow the strong-maximum-principle argument of Proposition~\ref{prop:finite}(ii), applied to the killed system, to give the assertion for every $T>0$ and every starting point.
\end{proof}

For $c\equiv1$, the density is $S_t=X_t/x$ and the shifted first factor is a three-dimensional Bessel process, recovering the classical change of measure from absorbed Brownian motion. For nonconstant $c$, the shifted first coordinate is not autonomous and the correlated factors acquire the singular drift in \eqref{eq:absorbed-shift}. The occupation estimate controls this drift in the necessity argument.

For example, with $k=1$, $\rho_1=1/2$ and $c(z)=a+b/(1+z^2)$, $a,b>0$, the exact condition is $a\ge1/2$. At $a=1/2$ the infimum is not attained. For $a=2/5$ and $b=1/5$ the density is strict even when $Z_0=0$, where $c(Z_0)=3/5$: a favourable initial value does not replace the global coefficient condition.

\section{Consistency with the affine theory}\label{sec:affine}

For affine models the martingale property of exponentially affine local martingales is characterized through generalized Riccati equations \cite{KMK10,MMS11,KRM15}. This section checks that the criteria of Sections~\ref{sec:boundary}--\ref{sec:sing} give the same answer and identifies the ingredients on both sides. The upshot is that for affine \emph{diffusions} the answer is always ``true martingale'', that the Fichera drift at the faces is literally the Feller condition and is not changed by the change of measure (Corollary~\ref{cor:inherit} in explicit form), and that on the cone of positive semidefinite matrices the Fichera drift reproduces the exact Wishart boundary condition, provided the curvature term of Lemma~\ref{lem:T} is kept.

\subsection{Canonical affine diffusions}
Let $d=m+n$, $I=\{1,\dots,m\}$, $J=\{m+1,\dots,d\}$, $D=(0,\infty)^m\times\R^n$, and let $(a_0,\alpha_1,\dots,\alpha_m,b^0,B)$ be admissible in the sense of \cite[Def.~2.6]{DFS03} (diffusion case): $a(y)=a_0+\sum_{i\in I}y_i\alpha_i$, $b(y)=b^0+By$, with $a_0,\alpha_i\in\mathbb S^d_+$, $(a_0)_{kl}=0$ unless $k,l\in J$, $(\alpha_i)_{kl}=0$ whenever $k\in I\setminus\{i\}$ or $l\in I\setminus\{i\}$, $b^0_I\in\R^m_+$, $B_{ij}\ge0$ for $i\ne j$ in $I$, and $B_{IJ}=0$. We assume $(a_0)_{JJ}\succ0$ and $(\alpha_i)_{ii}>0$, so that $a\succ0$ on $D$, and the boundary non-attainment condition
\begin{equation}\label{eq:feller-affine}
b^0_i\ \ge\ \tfrac12(\alpha_i)_{ii}\qquad(i\in I),
\end{equation}
which by Proposition~\ref{thm:attain} and Proposition~\ref{prop:sing} (see the proof below) is exactly the condition for $Y$ not to leave $D$ under $\Prob$. For $\Lambda\in\R^d$ the exponentially affine local martingale is
\[
M^\Lambda_t:=\cE\Big(\int_0^{\cdot}\Lambda^{\top}\sigma(Y_s)\,dW_s\Big)_t=\exp\Big(\Lambda\cdot(Y_t-Y_0)-\int_0^t\big(F(\Lambda)+R(\Lambda)\cdot Y_s\big)ds\Big),
\]
where $F(u)=\tfrac12u^{\top}a_0u+u\cdot b^0$, $R_i(u)=\tfrac12u^{\top}\alpha_iu+(B^{\top}u)_i$ for $i\in I$ and $R_k(u)=(B^{\top}u)_k$ for $k\in J$; here $\theta=\sigma^{\top}\Lambda$ and $\sigma\theta=a\Lambda$.

\begin{proposition}\label{prop:affine}
Under the above assumptions $M^\Lambda$ is a true martingale for every $\Lambda\in\R^d$.
\end{proposition}

\begin{proof}
\emph{Riccati route.} The drift under the changed measure is $\mut(y)=b(y)+a(y)\Lambda=\tilde b^0+\tilde By$ with $\tilde b^0=b^0+a_0\Lambda$ and $\tilde B=B+[\alpha_1\Lambda\,|\cdots|\,\alpha_m\Lambda\,|\,0]$. These parameters are admissible: $\tilde b^0_I=b^0_I$ since $(a_0)_{I\cdot}=0$; $\tilde B_{ji}=B_{ji}$ for $j\in I\setminus\{i\}$ since $(\alpha_i)_{j\cdot}=0$; $\tilde B_{IJ}=0$. Hence the drift-changed process is an affine diffusion, and it is conservative because the associated Riccati function is a polynomial, hence locally Lipschitz \cite[Sec.~9]{DFS03}. By \cite[Thm.~3.2]{KMK10} (sufficiency) and \cite{MMS11} (necessity), $M^\Lambda$ is a true martingale. Equivalently, by the affine transform formula with discounting, $\E_y[M^\Lambda_T]=\exp(\Phi(T)+(\Psi(T)-\Lambda)\cdot y)$ where $\Psi'=R(\Psi)-R(\Lambda)$, $\Psi(0)=\Lambda$, $\Phi'=F(\Psi)-F(\Lambda)$, $\Phi(0)=0$; the unique solution is $\Psi\equiv\Lambda$, $\Phi\equiv0$.

\emph{Route through Sections~\ref{sec:boundary}--\ref{sec:sing}.} (i) \emph{Infinity:} $\tr a+2\ip{y}{\mut}$ is a polynomial of degree at most two in $y$, so \eqref{eq:lingrowth} holds and Proposition~\ref{cor:osgood}(i) gives no explosion. (ii) \emph{Faces:} for $i\in I$, $a_{ii}(y)=(\alpha_i)_{ii}y_i$, $a_{ij}=0$ for $j\in I\setminus\{i\}$, $a_{ik}=y_i(\alpha_i)_{ik}$ for $k\in J$, and $(a\Lambda)_i=y_i(\alpha_i\Lambda)_i$; so on the face $\{y_i=0\}$, with $\nu=-e_i$, $\ip{a\nu}{\nu}=0$, the Girsanov shift vanishes (this is Corollary~\ref{cor:inherit}: $\theta=\sigma^{\top}\Lambda$ is bounded near the faces), and
\[
\Ft=F=-\Big(b^0_i+\sum_{j\in I\setminus\{i\}}B_{ij}y_j-\tfrac12(\alpha_i)_{ii}\Big)\ \le\ -\big(b^0_i-\tfrac12(\alpha_i)_{ii}\big)\ \le0
\]
by \eqref{eq:feller-affine} and $B_{ij}\ge0$; conversely $\Ft>0$ somewhere on the face if \eqref{eq:feller-affine} fails, so \eqref{eq:feller-affine} is exactly the non-exit condition for $\Prob$ (Proposition~\ref{thm:attain}(i)). Condition (T) holds on the characteristic boundary where the Fichera drift vanishes: $a_{ii}/y_i-2\mut_i=-2(b^0_i-\tfrac12(\alpha_i)_{ii})-2\sum_{j\in I\setminus\{i\}}B_{ij}y_j-2\tilde B_{ii}y_i\le2|\tilde B_{ii}|\,y_i$. Proposition~\ref{thm:attain}(ii) applies to all faces. (iii) \emph{Intersections of faces:} for $I_0\subset I$ with $|I_0|\ge2$ and $\pi$ the projection onto the $I_0$-coordinates, $A_\pi=\sum_{i\in I_0}(\alpha_i)_{ii}y_i^3/|\pi(y)|^2$, $B_\pi=\sum_{i\in I_0}(\alpha_i)_{ii}y_i$ and $C_\pi=2\sum_{i\in I_0}y_i\mut_i(y)\ge2\sum_{i\in I_0}y_i(b^0_i-|\tilde B_{ii}|y_i)$, because $\tilde B_{ij}\ge0$ for $j\in I\setminus\{i\}$ and $\tilde B_{iJ}=0$ (the coordinates outside $I_0$ only increase $C_\pi$ and do not enter $A_\pi,B_\pi$). Writing $\pi(y)=r\omega$ with $\omega$ on the unit sphere of $\R^{I_0}_+$,
\[
\frac{B_\pi-A_\pi+C_\pi}{A_\pi}\ \ge\ \frac{\sum_{i\in I_0}c_i\omega_i(1-\omega_i^2)+2\sum_{i\in I_0}\omega_ib^0_i}{\sum_{i\in I_0}c_i\omega_i^3}-Cr ,\qquad c_i:=(\alpha_i)_{ii},
\]
where the factor $r$ common to $A_\pi$, $B_\pi$ and the leading part of $C_\pi$ has cancelled (the ratio is of order one, not $1/r$), and $C$ depends only on the $|\tilde B_{ii}|$ because $D(\omega):=\sum_{i\in I_0}c_i\omega_i^3$ is bounded away from zero on the closed positive unit sphere. By \eqref{eq:feller-affine}, $2\sum_i\omega_ib^0_i\ge\sum_ic_i\omega_i$, so the numerator is at least $\sum_ic_i\omega_i(2-\omega_i^2)\ge\sum_ic_i\omega_i^3=D(\omega)$, and $\beta(r)\ge1-Cr$ for small $r$. Hence $\int_u^{r_0}\beta(s)s^{-1}ds\ge\log(r_0/u)-Cr_0$, \eqref{eq:betacond} holds, and Proposition~\ref{prop:sing} shows that no stratum is an exit set. (The bound is sharp: for two independent CIR factors at the critical Feller parameter, $c_i=1$, $b^0_i=\tfrac12$ and no linear drift, the infimum over the quarter circle is exactly $1$, attained on the axes.) By Lemma~\ref{thm:main}, $M^\Lambda$ is a true martingale.
\end{proof}

\begin{remark}
(a) If \eqref{eq:feller-affine} fails for some $i$, the process must be treated on the closed state space with the prescribed boundary behavior, which could be reflection on the face $\{y_i =0\}$ or absorption. Remark~\ref{rem:reflect} applies with the face included in the state space, and the conclusion persists; this matches the affine theory. (b) The check shows where the two approaches store the same information: conservativeness of the drift-changed affine process (global existence for the Riccati equation) corresponds to the quadratic bound \eqref{eq:lingrowth}; admissibility of the changed parameters ($\tilde b^0_I=b^0_I$) corresponds to the invariance of the Fichera drift at the faces.
\end{remark}

\subsection{Wishart processes}
Let $\mathbb S_n$ be the space of symmetric $n\times n$ matrices with the Frobenius inner product, $D=\mathbb S^{++}_n$ the open cone of positive definite matrices, and consider the Wishart process \cite{Bru91,CFMT11}
\begin{equation}\label{eq:wishart}
d\Sigma_t=\big(\Omega\Omega^{\top}+M\Sigma_t+\Sigma_tM^{\top}\big)dt+\sqrt{\Sigma_t}\,dW_t\,Q+Q^{\top}dW_t^{\top}\sqrt{\Sigma_t},
\end{equation}
with $W$ an $n\times n$ matrix Brownian motion and $Q$ invertible. A solution exists on the closed cone iff $\Omega\Omega^{\top}\succeq(n-1)Q^{\top}Q$, and it never reaches the boundary of the cone iff $\Omega\Omega^{\top}\succeq(n+1)Q^{\top}Q$ \cite{Bru91,MPS11,CFMT11}; both conditions are independent of $M$. The smooth part of $\partial D$ consists of the matrices of rank $n-1$; near such a point the distance to the boundary is the smallest eigenvalue, $\rho(\Sigma)=\lambda_{\min}(\Sigma)$, which is real-analytic where it is simple.

\begin{proposition}\label{prop:wishart}
Let $\Sigma^0\in\partial D$ have rank $n-1$ and let $v$ be a unit vector spanning its kernel. Then $\ip{a\nu}{\nu}(\Sigma^0)=0$ and the Fichera drift of \eqref{eq:wishart} at $\Sigma^0$ (outward normal $\nu=-vv^{\top}$) is
\begin{equation}\label{eq:wishart-fichera}
F(\Sigma^0)=-\Big(v^{\top}\Omega\Omega^{\top}v-(n+1)\,v^{\top}Q^{\top}Qv\Big).
\end{equation}
In particular $F\le0$ on the smooth part of $\partial D$ iff $\Omega\Omega^{\top}\succeq(n+1)Q^{\top}Q$, $F<0$ iff $\succ$, and $F$ does not depend on $M$.
\end{proposition}

\begin{proof}
Let $\lambda_1\ge\dots\ge\lambda_n$ be the eigenvalues of $\Sigma\in D$ near $\Sigma^0$ with orthonormal eigenvectors $v_1,\dots,v_n$, so that $\lambda_n=\rho$ is simple and $v_n\to v$. Perturbation theory for a simple eigenvalue gives $\nabla\rho=v_nv_n^{\top}$ and $\nabla^2\rho[E,E]=2\sum_{j<n}(v_n^{\top}Ev_j)^2/(\lambda_n-\lambda_j)$. Writing $\ip{x}{dW\,y}$ for $x^{\top}dW\,y$, the martingale part of $v_n^{\top}d\Sigma\,v_j$ is $\sqrt{\lambda_n}\ip{v_n}{dW\,Qv_j}+\sqrt{\lambda_j}\ip{v_j}{dW\,Qv_n}$, whose quadratic variation rate is $\lambda_n|Qv_j|^2+\lambda_j|Qv_n|^2$ (the cross term vanishes because $v_n\perp v_j$). Hence
\[
h=\ip{a\nabla\rho}{\nabla\rho}=4\lambda_n|Qv_n|^2,\qquad L\rho=v_n^{\top}\big(\Omega\Omega^{\top}+M\Sigma+\Sigma M^{\top}\big)v_n+\sum_{j<n}\frac{\lambda_n|Qv_j|^2+\lambda_j|Qv_n|^2}{\lambda_n-\lambda_j} ,
\]
and as $\Sigma\to\Sigma^0$: $h/\rho\to4|Qv|^2$, $v_n^{\top}(M\Sigma+\Sigma M^{\top})v_n=2\lambda_nv_n^{\top}Mv_n\to0$, and the sum tends to $-(n-1)|Qv|^2$. By \eqref{eq:limitT}, $2F(\Sigma^0)=\lim(h/\rho-2L\rho)=4|Qv|^2-2v^{\top}\Omega\Omega^{\top}v+2(n-1)|Qv|^2$, which is \eqref{eq:wishart-fichera}. The first assertion follows from $h(\Sigma^0)=0$.
\end{proof}

\begin{remark}
(a) The sum in $L\rho$ is $\tfrac12\tr(a\nabla^2\rho)$, the curvature term of Lemma~\ref{lem:T}. Ignoring it, i.e.\ treating $v^{\top}\Sigma v$ as a CIR process with $v$ frozen, gives $v^{\top}\Omega\Omega^{\top}v-2|Qv|^2$ and the threshold $\Omega\Omega^{\top}\succeq2Q^{\top}Q$, which is wrong for $n\ge2$; the correct threshold $(n+1)$ is the squared-Bessel dimension $2$ for the smallest eigenvalue in Bru's eigenvalue dynamics \cite{Bru91}. (b) In the Wishart stochastic volatility model of Da Fonseca--Grasselli--Tebaldi \cite{DFGT08}, $dS_t/S_t=\tr(\sqrt{\Sigma_t}\,dZ_t)$ with $Z=WR^{\top}+B(I-RR^{\top})^{1/2}$, the coefficient of the stochastic exponential is bounded on bounded subsets of $\bar D$, hence bounded near every boundary point (the boundary of the cone is itself unbounded, but Proposition~\ref{prop:inherit} and Corollary~\ref{cor:inherit} are local), so the boundary classification is not affected by the measure change, and the Girsanov shift replaces $M$ by $\tilde M=M+Q^{\top}R^{\top}$, which is again linear. With $\Omega\Omega^{\top}\succeq(n+1)Q^{\top}Q$ the process does not reach the boundary under $\Prob$, the drift-changed process does not reach it either (Proposition~\ref{prop:wishart}, whose right side is $M$-independent) and does not explode (linear growth), so $S$ is a true martingale for all $M,Q,R$; this is the statement that the drift-changed Wishart process is a conservative affine process on $\mathbb S^+_n$ \cite{CFMT11,KMK10}. The lower-rank strata of the boundary are all contained in $\{\lambda_{\min}=0\}$ and are excluded by the same global results.
\end{remark}

\section{Conclusion}\label{sec:open}

For the stochastic exponential $S$, defined by $dS_t=S_t\theta(Y_t)\cdot dW_t$ and $S_0=1$, necessary and sufficient coefficient conditions for martingality can be obtained by combining local boundary analysis with a test for explosion at infinity. Here $Y$ is the factor diffusion, $W$ its $m$-dimensional driving Brownian motion, and $\theta:D\to\R^m$ specifies the integrand $\theta(Y_t)$ in the stochastic exponential. The Fichera criterion provides such an equivalence under the stated hypotheses at boundary points with zero Fichera drift, while bounded $\theta$ preserves the boundary behaviour of the original diffusion. This boundary analysis uses the full multidimensional generator and the local boundary geometry; the distance to the boundary need not satisfy an autonomous scalar stochastic differential equation. Boundary singularities in $\theta$ with correlated modulation require control of the tangential drift as well as the normal motion. This control yields the exact thresholds for the square-root models of Theorem~\ref{thm:correlated}, the logarithmic family of Theorem~\ref{thm:coupled-log}, and the inverse-distance family of Proposition~\ref{thm:weighted-stratum}. Each threshold includes equality, even when the relevant extremum of the modulation is not attained. In the mean-reverting logarithmic model, every nontrivial exponential converges to zero almost surely, including in the finite-horizon martingale regime; uniform integrability holds exactly when $c\equiv0$ (Corollary~\ref{cor:coupled-log-ui}).

The choice of coordinate also determines the scope of scalar comparison. Weighted distance can establish nonattainment when Euclidean projection is inconclusive, and a scalar function with compact sublevel sets gives an exact Feller test at infinity, including time-change rates comparable to a positive function of that coordinate (Corollary~\ref{cor:variable-clock}). The resulting classes include anisotropic diffusions with nonlinear motion along the level sets. For factors that can reach the boundary under the original measure, the criterion depends additionally on the prescribed boundary behaviour and the density at absorption. Reflection preserves mass at boundary visits under the hypotheses of Proposition~\ref{prop:closed-state}; Theorem~\ref{thm:absorbed} gives an exact coefficient condition when a singular integrand modulated by correlated factors sends the density to zero.

The remaining boundary questions concern critical points where condition (T) fails and the normal and tangential coefficients do not have the structures treated here, as well as intersections of boundary faces not decided by the available non-attainment tests. At infinity, the unresolved cases are those in which neither radial comparison nor comparison through a function with compact sublevel sets decides explosion. Accessible boundaries at which some density mass survives introduce the further term in \eqref{eq:retained-mass}. The boundary questions require control of the joint normal and tangential dynamics.

\appendix
\section{Regularity of the loss of mass}\label{app:reg}

This appendix contains the regularity lemma used in Proposition~\ref{prop:finite}(ii), Remark~\ref{rem:pde} and Appendix~\ref{sec:pde}. Recall that the exhaustion $(D_n)$ of $D$ is ours to choose and that we take $\partial D_n\in C^2$.

\begin{lemma}\label{lem:regularity}
Let $a,\mut$ be locally H\"older continuous on $D$ and $\partial D_n\in C^2$, $\langle a(y) v, v \rangle>0$, for $y\in D$, $v\neq 0$. For $w_0\in C_b(D)$ the functions
\[
w(t,y):=\E_y\big[w_0(Y_t);\,t<\tau_D\big],\qquad m(t,y):=\Qt_y(\tau_D\le t),\qquad (t,y)\in[0,T]\times D,
\]
belong to $C([0,T]\times D)\cap C^{1,2}((0,T]\times D)$ and satisfy $\partial_tw=\Lt w$, $\partial_tm=\Lt m$ there, with $w(0,\cdot)=w_0$ and $m(0,\cdot)=0$.
\end{lemma}

\begin{proof}
By linearity we may take $w_0\ge0$. Fix $n$. On the bounded $C^2$ domain $D_n$ the operator $\Lt$ is uniformly elliptic with H\"older coefficients, so the first initial--boundary value problem $\partial_tv=\Lt v$ in $(0,T]\times D_n$, $v(0,\cdot)=w_0$, $v=0$ on $(0,T]\times\partial D_n$, has a bounded classical solution $v_n$, continuous on $[0,T]\times\bar D_n$ away from the corner $\{0\}\times\partial D_n$ \cite[Ch.~3]{Fri64} (approximate $w_0$ by functions vanishing near $\partial D_n$ and use the maximum principle together with the interior estimates below). Every point of $\partial D_n$ is regular \cite[Thm.~2.3.3]{Pin95}, and It\^o's formula applied to $s\mapsto v_n(t-s,Y_s)$ on $[0,t\wedge\tau_n)$ gives
\[
v_n(t,y)=\E_y\big[w_0(Y_t);\,t<\tau_n\big]+\E_y\big[v_n(t-\tau_n,Y_{\tau_n});\,\tau_n\le t\big]=\E_y\big[w_0(Y_t);\,t<\tau_n\big]=:w_n(t,y),
\]
because $v_n$ vanishes on the lateral boundary (for an approximant $w_0^k$ vanishing near $\partial D_n$ the solution is continuous on $[0,T]\times\bar D_n$; the limit $k\to\infty$ is monotone on the probabilistic side and controlled by the interior estimates on the analytic side). As $n\to\infty$, $w_n\uparrow w$ by monotone convergence, the $w_n$ are bounded by $\|w_0\|_\infty$, and by the interior Schauder estimates \cite[Ch.~3, Sec.~5]{Fri64} their $C^{1+\gamma/2,2+\gamma}$-norms are bounded on every compact subcylinder of $(0,T]\times D$; hence a subsequence converges in $C^{1,2}_{\mathrm{loc}}$, the limit $w$ is in $C^{1,2}((0,T]\times D)$ and solves the equation. Continuity at $t=0$: $|w(t,y)-w_0(y)|\le\E_y|w_0(Y_t)-w_0(y)|+\|w_0\|_\infty\Qt_y(\tau_D\le t)$, and both terms tend to $0$ as $t\downarrow0$, locally uniformly in $y$, by the continuity of the paths and $\tau_D\ge\tau_n$. The same argument applies to $m_n(t,y):=\Qt_y(\tau_n\le t)=1-\E_y[1;t<\tau_n]$, which decreases to $m$.
\end{proof}

\begin{remark}[on probabilistic regularity]\label{rem:regularity}
Fernholz--Karatzas \cite[Sec.~9]{FK10} carry out the same localization and need every point of the boundary of the auxiliary set to be probabilistically regular; this is the delicate point of that argument when the auxiliary sets are not chosen freely. Here the exhaustion is ours to choose, the operator is uniformly elliptic inside $D$, and the $C^2$ regularity of $\partial D_n$ settles the matter. Nothing of the sort is required at $\partial D$ itself, where the operator degenerates: the points of $\widetilde\Sig_1\cup\Sig_1^0$ are not attainable, when (T) holds for points of $\Sig_1^0$, and at the points of $\Sig_2\cup\Sig_3$ the function $U$ of Lemma~\ref{lem:exit}, with $U(x^0)=0$ and $\Lt U\le-\gamma$, is exactly a barrier in the sense of the classical regularity theory, now placed at a characteristic boundary point. With merely continuous coefficients (the Bass--Perkins-type conditions of \cite{FK10}) classical solutions must be replaced by $W^{2,1,p}_{\mathrm{loc}}$ solutions and It\^o's formula by the It\^o--Krylov formula; the uniqueness statements of Appendix~\ref{sec:pde} remain true in that class with the same proofs.
\end{remark}

\section{Uniqueness for the associated linear PDEs}\label{sec:pde}

For models satisfying Assumption~\ref{ass:standing}, the correspondence of Remark~\ref{rem:pde} turns the martingale criteria into uniqueness criteria for the Cauchy problem of a linear parabolic equation with degenerate boundary and unbounded coefficients, in the class of solutions with prescribed growth. The reflected and absorbed models of Sections~\ref{sec:extended-state} and \ref{sec:absorbed-model} have different boundary specifications. Throughout this appendix we assume, in addition to Assumption~\ref{ass:standing}(ii), that the coefficients are locally H\"older continuous on $D$, so that Lemma~\ref{lem:regularity} applies, and $\E_y$ denotes expectation under $\Qt_y$ unless said otherwise.

\subsection{The Cauchy problem}
Consider, for $w_0\in C_b(D)$,
\begin{equation}\label{eq:CP}
\partial_tw=\Lt w\ \text{ in }(0,T]\times D,\qquad w(0,\cdot)=w_0\ \text{ on }D,
\end{equation}
and call $w$ a bounded solution if $w\in C([0,T]\times D)\cap C^{1,2}((0,T]\times D)$ is bounded and satisfies \eqref{eq:CP}.

\begin{lemma}\label{prop:CP}
The following are equivalent:
\begin{enumerate}[label=(\roman*),leftmargin=2em]
\item $\Qt_y(\tau_D\le T)=0$ for all $y\in D$;
\item for every $w_0\in C_b(D)$, \eqref{eq:CP} has at most one bounded solution;
\item the only bounded solution of \eqref{eq:CP} with $w_0=0$ is $w\equiv0$.
\end{enumerate}
Moreover, $m(t,y):=\Qt_y(\tau_D\le t)$ is a bounded solution of \eqref{eq:CP} with $w_0=0$, nontrivial iff (i) fails; for $w_0\ge0$, $w_{\min}(t,y):=\E_y[w_0(Y_t);t<\tau_D]$ is the minimal nonnegative bounded solution; and under (i) the unique bounded solution is $w(t,y)=\E_y[w_0(Y_t)]$.
\end{lemma}

\begin{proof}
(ii)$\Rightarrow$(iii) is trivial and (iii)$\Rightarrow$(i) follows from the fact that $m$ is a bounded solution with zero initial data (Lemma~\ref{lem:regularity}). (i)$\Rightarrow$(ii): let $w$ be a bounded solution and fix $(t,y)$, $0<\varepsilon<t$. By It\^o's formula, $s\mapsto w(t-s,Y_s)$ is a bounded martingale on $[0,(t-\varepsilon)\wedge\tau_n]$, so $w(t,y)=\E_y[w(t-(t-\varepsilon)\wedge\tau_n,\,Y_{(t-\varepsilon)\wedge\tau_n})]$. Since $\tau_n\uparrow\tau_D>T$ a.s.\ by (i), bounded convergence gives $w(t,y)=\E_y[w(\varepsilon,Y_{t-\varepsilon})]$, and $\varepsilon\downarrow0$ gives $w(t,y)=\E_y[w_0(Y_t)]$ by the continuity of $w$ on $[0,T]\times D$ and of the paths. The minimality of $w_{\min}$ follows from the same argument stopped at $\tau_n$: a nonnegative solution $w$ satisfies $w(t,y)\ge\E_y[w_0(Y_t);t<\tau_n]$ for every $n$.
\end{proof}

The criteria of Sections~\ref{sec:boundary}--\ref{sec:sing} now read as follows.

\begin{proposition}\label{thm:pde-unique}
Let $\partial D\in C^2$, $a\in C^1(\bar D)$, $\mut\in C(\bar D)$, and let $\Ft=\ip{\mut-\tfrac12\dv a}{\nu}$ be the Fichera drift of $\Lt$ on $\partial D$.
\begin{enumerate}[label=(\roman*),leftmargin=2em]
\item If $\ip{a\nu}{\nu}>0$ or $\Ft>0$ at some point of $\partial D$, then \eqref{eq:CP} with $w_0=0$ has a nontrivial bounded nonnegative solution: uniqueness in the bounded class fails, for every $T>0$.
\item If $\ip{a\nu}{\nu}=0$ and $\Ft\le0$ on $\partial D$, (T) holds on $\{\Ft=0\}$, and $\Qt$ does not explode at infinity (e.g.\ $\tr a+2\ip{y}{\mut}\le C(1+|y|^2)$, or Lemma~\ref{thm:radial-sandwich}(i)), then \eqref{eq:CP} has exactly one bounded solution for every $w_0\in C_b(D)$, namely $\E_y[w_0(Y_t)]$.
\item If $D=\R^d$ and $\mathcal A,\mathcal B$ in \eqref{eq:AB} are radial outside a ball, then uniqueness in the bounded class holds iff $v_{A,B}(\infty)=\infty$.
\item On $D=(0,\infty)^d$, uniqueness in the bounded class holds if the faces satisfy the conditions in (ii) and every stratum satisfies \eqref{eq:betacond}; it fails if some face contains a point with $\Ft>0$.
\end{enumerate}
\end{proposition}

\begin{proof}
Lemma~\ref{prop:CP} together with Proposition~\ref{thm:attain}, Proposition~\ref{thm:radial-iff}, Proposition~\ref{prop:sing} and Proposition~\ref{prop:finite}(i) (which shows that (i) of Lemma~\ref{prop:CP} for one $T$ is the same as for all $T$).
\end{proof}

\begin{remark}
(a) For bounded domains with smooth coefficients this is Fichera's theory \cite{Fic56,OR73}: boundary data are prescribed on $\Sig_2\cup\Sig_3$. The probabilistic route gives the same dichotomy for unbounded domains with the condition at infinity built in, for merely continuous drifts, and identifies the additional solutions explicitly ($m$ and $w_{\min}$). (b) The elliptic form is a Liouville statement: under the hypotheses of (ii) the only bounded solution of $\Lt u=\lambda u$ in $D$, $\lambda>0$, is $u\equiv0$, while under (i) $u(y)=\E_y[e^{-\lambda\tau_D}]$ is a nontrivial bounded positive solution (Lemma~\ref{thm:khas}). (c) On $D=\R^d$ the statement ``uniqueness of bounded solutions of the heat equation $\iff$ non-explosion'' is the stochastic completeness theory of Azencott \cite{Aze74} and Grigor'yan \cite{Gri99}; Khasminskii's radial test and Proposition~\ref{cor:osgood} are the Euclidean, variable-coefficient counterparts of Grigor'yan's volume-growth tests, and Proposition~\ref{thm:radial-iff} is their sharp form when $\mathcal A$ and $\mathcal B$ depend only on the radius. The uniqueness of \emph{nonnegative} solutions (Widder-type classes) is a different question, see Ishige--Murata \cite{IM01}.
\end{remark}

\subsection{Growth classes and valuation equations}
Now let $S$ be the asset of Section~\ref{sec:SV} and consider the joint state $(S,Y)$ on $G:=(0,\infty)\times D$, whose generator is
\begin{equation}\label{eq:pricing-op}
\mathcal L u=\tfrac12s^2|\theta(y)|^2\partial_{ss}u+s\sum_{i}(\sigma\theta)_i(y)\,\partial_{sy_i}u+\tfrac12\sum_{i,j}a_{ij}(y)\partial_{y_iy_j}u+\sum_i\mu_i(y)\partial_{y_i}u ,
\end{equation}
and the valuation equation in forward time,
\begin{equation}\label{eq:VE}
\partial_\tau u=\mathcal Lu\ \text{ in }(0,T]\times G,\qquad u(0,\cdot)=\varphi\ \text{ on }G,
\end{equation}
whose probabilistic solution is the price $u(\tau,s,y)=\E[\varphi(S_\tau,Y_\tau)\,|\,S_0=s,Y_0=y]$. Two growth classes are natural: $\mathcal C_b$, the bounded solutions, and $\mathcal C_{\mathrm{lin}}:=\{u:\sup|u(\tau,s,y)|/(1+s)<\infty\}$, which contains the call payoff and the forward $\varphi(s,y)=s$. The functions $1$ and $s$ are $\mathcal L$-harmonic, and $h(s,y):=1+s$ is the natural weight for $\mathcal C_{\mathrm{lin}}$: the Doob transform $\mathcal L^hw:=h^{-1}\mathcal L(hw)=\mathcal Lw+\ip{a^{\mathcal L}\nabla\log h}{\nabla w}$ is, when $S$ is a true martingale on $[0,T]$, the generator of $(S,Y)$ under the probability measure $\Prob^h$ with density $(1+S_T)/(1+s)$, the mixture $\tfrac1{1+s}\Prob+\tfrac{s}{1+s}\Qt$ on $\mathcal F_T$; here $a^{\mathcal L}$ is the diffusion matrix in \eqref{eq:pricing-op}, and the extra drift in the $y$-directions is $\tfrac{s}{1+s}\sigma\theta$. (When $S$ is a strict local martingale the density has mass less than one and $\mathcal L^h$ is the generator of a killed process; we only use $\Prob^h$ in the martingale case.)

\begin{theorem}\label{thm:VE}
Let Assumption~\ref{ass:standing} hold with locally H\"older coefficients.
\begin{enumerate}[label=(\roman*),leftmargin=2em]
\item (\emph{bounded class}) For every $\varphi\in C_b(G)$, \eqref{eq:VE} has at most one solution in $\mathcal C_b$, and if a solution exists it is the price.
\item (\emph{linear-growth class}) If $\Qt_y(\tau_D\le T)=0$ for all $y\in D$, i.e.\ if $S$ is a true martingale on $[0,T]$ from every starting point, then for every continuous $\varphi$ with $|\varphi|\le C(1+s)$, \eqref{eq:VE} has at most one solution in $\mathcal C_{\mathrm{lin}}$. If $\Qt_y(\tau_D\le T)>0$ for some $y$, then \eqref{eq:VE} with $\varphi(s,y)=s$ has the two distinct solutions
\[
u_1(\tau,s,y)=s\qquad\text{and}\qquad u_2(\tau,s,y)=s\,\E_y[S_\tau]=s\big(1-m(\tau,y)\big)
\]
in $\mathcal C_{\mathrm{lin}}$, the second being the price of the forward.
\item Consequently, under the regularity hypotheses of Theorem~\ref{thm:fichera} and non-explosion at infinity under $\Qt$: uniqueness in $\mathcal C_{\mathrm{lin}}$ holds if $\Ft=\langle \mu - \frac12 \text{div} a, \nu \rangle+\ip{\sigma\theta}{\nu}<0$ on $\partial D$ (or $\le0$ and $\Lt$ satisfies condition (T) at every point of $\{\Ft =0\}$ near which $\theta$ is unbounded), and fails if $\Ft>0$ somewhere on $\partial D$.
\end{enumerate}
\end{theorem}

Classical existence requires regularity of the joint operator $\mathcal L$, which the factor ellipticity does not provide (for $Y$ a Brownian motion, $\theta=0$ and $\varphi(s,y)=\min\{|s-1|,1\}$ the price is $\varphi$ itself, which is not $C^2$ in $s$).

\begin{proof}
(i) is the uniqueness part of Lemma~\ref{prop:CP} for the $(S,Y)$-diffusion under $\Prob$, which never leaves $G$ ($S>0$ is finite and $Y$ stays in $D$); It\^o's formula along the paths identifies any bounded classical solution with the price. (ii) For the first assertion put $w:=u/h$; then $w$ is a bounded solution of $\partial_\tau w=\mathcal L^hw$, and since $S$ is a true martingale on $[0,T]$, $\Prob^h$ is a probability measure under which the $(S,Y)$-diffusion does not leave $G$ by time $T$ (neither $\Prob$ nor $\Qt$ do), so $w$ is unique by the argument of Lemma~\ref{prop:CP}. (Equivalently, one may avoid $\Prob^h$ altogether: $s\mapsto u(\tau-s,S_s,Y_s)$ is a local martingale bounded by $C(1+S_s)$, hence a true martingale by the uniform integrability of the martingale $S$ on $[0,T]$.) For the second, $u_1$ solves \eqref{eq:VE} since $s$ is $\mathcal L$-harmonic; $u_2$ is the price of the claim $S_T$ by the Markov property and the multiplicative structure of $S$, and it is a classical solution because $m$ is (Lemma~\ref{prop:CP} for $\Lt$); $u_2\ne u_1$ since $m\not\equiv0$. (iii) is Theorem~\ref{thm:fichera} together with (ii) and Proposition~\ref{prop:finite}(i).
\end{proof}

\begin{example}
(a) For the $4/2$ model of Example~\ref{ex:42} the valuation operator on $G=(0,\infty)^2$ is
\[
\mathcal L=\tfrac12s^2\Big(a\sqrt y+\frac{b}{\sqrt y}\Big)^2\partial_{ss}+s\rho\eta(ay+b)\,\partial_{sy}+\tfrac12\eta^2y\,\partial_{yy}+(\alpha-\kappa y)\partial_y ,
\]
with a coefficient of $\partial_{ss}$ that is singular on $\{y=0\}$. Under $2\alpha\ge\eta^2$ the Cauchy problem has at most one solution in $\mathcal C_b$, and uniqueness holds in $\mathcal C_{\mathrm{lin}}$ iff $\alpha+\rho\eta b\ge\tfrac12\eta^2$; for $b=0$ (Heston) the two classes behave alike, for $a=0$ this is the $3/2$ valuation equation. (b) For the two-factor models of Examples~\ref{ex:twofactor} and \ref{ex:42}, uniqueness in $\mathcal C_{\mathrm{lin}}$ holds iff $2\alpha_1-\eta_1^2\ge2\eta_1(\rho b)^-\sup g$ (with $b=1$ for inverse-CIR volatility). For the correlated $d$-factor model of Theorem~\ref{thm:correlated}, the corresponding criterion is \eqref{eq:correlated-iff}, by Theorem~\ref{thm:VE}(ii). (c) For the Wishart model of Section~\ref{sec:affine}, on the open cone $G=(0,\infty)\times\mathbb S^{++}_n$ (the model on the closed cone, with the boundary behaviour of the process included, is a different problem) with the pricing operator of \cite{DFGT08}: uniqueness holds in $\mathcal C_b$ and in $\mathcal C_{\mathrm{lin}}$ iff $\Omega\Omega^{\top}\succeq(n+1)Q^{\top}Q$ (Proposition~\ref{prop:wishart}; $\theta$ is locally bounded near the boundary, so the two classes coincide). The uniqueness statements in (a) for $b>0$ and in (b), (c) are new; for one-factor models with $\theta$ bounded near finite boundary points, such as the Heston case $b=0$, uniqueness in $\mathcal C_{\mathrm{lin}}$ is in \cite{ET10}.
\end{example}

\begin{proof}[Proof of the converse in (c)]
Suppose $\Omega\Omega^{\top}\not\succeq(n+1)Q^{\top}Q$. Choose a unit vector $v$ such that, at the rank-$(n-1)$ boundary point $\Sigma^0=I-vv^{\top}$, Proposition~\ref{prop:wishart} gives
\[
F(\Sigma^0)=v^{\top}\big[(n+1)Q^{\top}Q-\Omega\Omega^{\top}\big]v>0.
\]
Let $\Prob_{\Sigma}$ be the law of the original Wishart factor in \eqref{eq:wishart}, started at $\Sigma\in D=\mathbb S_n^{++}$ and stopped at its first boundary hit $\tau_{\partial}$. The drift and covariance coefficients are smooth up to the smooth characteristic boundary patch containing $\Sigma^0$, so Lemma~\ref{thm:local}(i) applies. For a sufficiently small neighbourhood $W=D\cap B(\Sigma^0,\varepsilon)$, Lemma~\ref{lem:exit} gives $\E_{\Sigma}[\tau_W]\to0$ and $\Prob_{\Sigma}(\Sigma_{\tau_W}\in\partial D)\to1$ as $\Sigma\to\Sigma^0$ from $D$, where the expectation is under $\Prob_{\Sigma}$. Hence, for every $T>0$,
\[
\Prob_{\Sigma}(\tau_{\partial}\le T)
\ge \Prob_{\Sigma}(\Sigma_{\tau_W}\in\partial D)
   -\frac{\E_{\Sigma}[\tau_W]}{T}>0
\]
for $\Sigma$ sufficiently close to $\Sigma^0$.

Write $L_{\Sigma}$ for the factor generator. Its covariance quadratic form is
\[
\langle a(\Sigma)H,H\rangle
=4\tr(\Sigma H Q^{\top}QH)>0
\qquad(\Sigma\in D,\ H\in\mathbb S_n\setminus\{0\}),
\]
since $Q$ is invertible. Thus $L_{\Sigma}$ has smooth coefficients and is locally uniformly elliptic in $D$. Lemma~\ref{lem:regularity}, applied to this stopped factor under $\Prob_{\Sigma}$, shows that
\[
m(t,\Sigma):=\Prob_{\Sigma}(\tau_{\partial}\le t)
\]
is a bounded classical solution of $\partial_t m=L_{\Sigma}m$ with $m(0,\cdot)=0$, and the preceding estimate makes it nontrivial on every positive time horizon. For functions independent of $s$, the pricing operator $\mathcal L$ reduces to $L_{\Sigma}$. Consequently $u_0(t,s,\Sigma)=0$ and $u_1(t,s,\Sigma)=m(t,\Sigma)$ are distinct solutions of \eqref{eq:VE} with zero initial data on $G$ and no prescribed boundary values. Since $0\le u_1\le1$, both belong to $\mathcal C_b$ and $\mathcal C_{\mathrm{lin}}$. This proves necessity in both classes.
\end{proof}

\begin{remark}
(a) In one dimension, Theorem~\ref{thm:VE}(ii)--(iii) is the characterization of uniqueness for the Black--Scholes and stochastic-volatility valuation equations of Ekstr\"om--Tysk \cite{ET09,ET10}, Bayraktar--Xing \cite{BX10} and \cite{BKX12}, where the criterion is Feller's test; in stochastic portfolio theory, Fernholz--Karatzas \cite[Thm.~1]{FK10} characterize the optimal arbitrage $U=1-m$ as the smallest nonnegative supersolution of the Cauchy problem with initial value $1$, so that arbitrage ($U<1$) is precisely the failure of uniqueness, and \cite[Prop.~3]{FK10} decides it through the Fichera drifts on the faces of the orthant; Cox--Hobson \cite{CH05} is the financial origin of the ``bubble'' interpretation of the missing mass. The uniqueness mechanism itself (Lemma~\ref{thm:main} plus a Doob transform) is dimension-independent and not new; what is specific to the multi-dimensional case is the description of the uniqueness class through the shifted Fichera drift, i.e.\ the coefficient conditions of Sections~\ref{sec:boundary}--\ref{sec:sing}; the analytic uniqueness theory for Heston-type degenerate operators in weighted H\"older spaces \cite{FP15,EM13} concerns regularity in weighted H\"older spaces rather than uniqueness in the class of functions with linear growth. (b) Power growth classes $|u|\le C(1+s^p)$, $p>1$, correspond to the moments $\E[S_T^p]$; the weight $s^p$ is not $\mathcal L$-harmonic, the Doob transform acquires the potential $\tfrac12p(p-1)|\theta|^2$. Finite horizon moments and their explosions have been studied by \cite{AP07, KR11}. A separate analysis is needed for uniqueness in the power growth class.
\end{remark}

\section{Proofs of the preliminary results}\label{app:prelim-proofs}

This appendix proves the results stated in Section~\ref{sec:exact}. The stopping times for the open-state-space model are the common coordinate functionals defined in \eqref{eq:canonical-exit}; the law is indicated by the probability or expectation symbol.

\begin{proof}[Proof of Lemma~\ref{thm:main}]
Fix $n$ and $T$. On $D_n$ the coefficients are bounded, so Novikov's condition holds for the stopped exponential $S^{\tau_n}$, which is a true martingale. Define $\Prob^n:=S_{T\wedge\tau_n}\cdot\Prob$ on $\mathcal F_T$. By Girsanov, $W-\int_0^{\cdot\wedge\tau_n}\theta(Y_s)ds$ is a $\Prob^n$-Brownian motion, hence under $\Prob^n$ the process $Y^{\tau_n}$ solves, up to $\tau_n$, the SDE with drift $\mut$ and diffusion $\sigma$. Uniqueness in law for the $\Lt$-martingale problem stopped at the exit from $D_n$ (Assumption~\ref{ass:standing}(ii) and \cite[Thm.~1.13.1]{Pin95}) gives $\Prob^n(\tau_n>T)=\Qt_y(\tau_n>T)$. Since $\{\tau_n>T\}\in\mathcal F_{T\wedge\tau_n}$,
\[
\E_y[S_T\mathbf 1_{\{\tau_n>T\}}]=\E_y[S_{T\wedge\tau_n}\mathbf 1_{\{\tau_n>T\}}]=\Prob^n(\tau_n>T)=\Qt_y(\tau_n>T).
\]
Let $n\to\infty$: on the left, $\tau_n\uparrow\tau_D=\infty$ $\Prob$-a.s.\ by Assumption~\ref{ass:standing}(iii) and monotone convergence gives $\E_y[S_T]$; on the right $\{\tau_n>T\}\uparrow\{\tau_D>T\}$ because $\tau_n\uparrow\tau_D$. This proves \eqref{eq:mass}; (i) and (ii) follow because $t\mapsto\E_y[S_t]$ is non-increasing.
\end{proof}

\begin{proof}[Proof of Proposition~\ref{prop:finite}]
(i) By Lemma~\ref{thm:main}, the hypothesis gives $\Qt_y(\tau_D>T)=1$ for every $y\in D$. The Markov property of the factor coordinate under $\Qt$ yields
\[
\Qt_y(\tau_D>2T)
=\E_y^{\Qt}\big[\mathbf1_{\{\tau_D>T\}}\Qt_{Y_T}(\tau_D>T)\big]
=1.
\]
Iteration gives $\Qt_y(\tau_D>kT)=1$ for every integer $k\ge1$, hence $\Qt_y(\tau_D=\infty)=1$. Lemma~\ref{thm:main} concludes.

(ii) By Lemma~\ref{lem:regularity} in Appendix~\ref{app:reg}, $m(t,y):=\Qt_y(\tau_D\le t)=1-\E_y[S_t]$ is a nonnegative classical solution of $\partial_tm=\Lt m$ on $(0,\infty)\times D$. Suppose $m(t_0,y_0)=0$ for some $t_0>0$, $y_0\in D$. Since $\Lt$ is uniformly elliptic with continuous coefficients on compact subsets of $D$ and $(0,t_0]\times D$ is connected, the strong maximum principle for parabolic equations \cite[Ch.~2, Thm.~5]{Fri64}, applied to $-m$, which attains its maximum $0$ at the interior point $(t_0,y_0)$, gives $m\equiv0$ on $(0,t_0]\times D$. The Markov property at time $t_0$ then gives, for every $y\in D$,
\[
\Qt_y(\tau_D\le2t_0)\le\Qt_y(\tau_D\le t_0)+\E_y\big[\Qt_{Y_{t_0}}(\tau_D\le t_0);\,t_0<\tau_D\big]=0 ,
\]
and inductively $m\equiv0$ on $(0,\infty)\times D$: $S$ is a martingale on $[0,\infty)$ from every starting point. Otherwise $m>0$ on $(0,\infty)\times D$, i.e.\ $\E_y[S_T]<1$ for every $T>0$ and $y\in D$.
\end{proof}

With $\partial D_n\in C^2$, every boundary point of $D_n$ is regular for the killed diffusion by the exterior cone criterion \cite[Thm.~2.3.3]{Pin95}. Under the local H\"older hypothesis of Lemma~\ref{thm:khas}, the Dirichlet problem $\Lt u=\lambda u$ in $D_n$, $u=1$ on $\partial D_n$, has a classical solution $u\in C^2(D_n)\cap C(\bar D_n)$ \cite[Thm.~6.13]{GT01}. It\^o's formula identifies it with $u_n(x)=\E_x^{\Qt}[e^{-\lambda\tau_n}]$. The boundary regularity here concerns the auxiliary domains $D_n$.

\begin{proof}[Proof of Lemma~\ref{thm:khas}]
By Lemma~\ref{thm:main} we must show: $\Qt_y(\tau_D<\infty)>0$ for some $y$ iff $\Lt u=\lambda u$ has a bounded positive solution. If $\Qt_y(\tau_D<\infty)>0$, then $u(y):=\E^{\Qt}_y[e^{-\lambda\tau_D}]$ is bounded, positive (by \cite[Thm.~1.15.1]{Pin95}) and solves $\Lt u=\lambda u$. Conversely let $0<u_0\le1$ solve $\Lt u_0=\lambda u_0$ and let $u_n(y):=\E^{\Qt}_y[e^{-\lambda\tau_n}]$, which solves $\Lt u_n=\lambda u_n$ in $D_n$ with $u_n=1$ on $\partial D_n$ (see above). The maximum principle applied to $u_n-u_0$ (a negative minimum at $y_0\in D_n$ would give $0\le\Lt(u_n-u_0)(y_0)=\lambda(u_n-u_0)(y_0)<0$) yields $u_n\ge u_0$ on $D_n$; letting $n\to\infty$, $\E_y^{\Qt}[e^{-\lambda\tau_D}]\ge u_0(y)>0$, so $\Qt_y(\tau_D<\infty)>0$.
\end{proof}

\begin{proof}[Proof of Lemma~\ref{thm:lyap}]
(i) Let $\zeta_k$ be the exit time from $\{V<k\}$ (a relatively compact open subset of $D$), so $\zeta_k\uparrow\tau_D$. By It\^o's formula $e^{-\lambda(t\wedge\zeta_k)}V(Y_{t\wedge\zeta_k})$ is a supermartingale, whence $k\,e^{-\lambda t}\Qt_y(\zeta_k\le t)\le\E_y[e^{-\lambda(t\wedge\zeta_k)}V(Y_{t\wedge\zeta_k})]\le V(y)$ and $\Qt_y(\tau_D\le t)\le\lim_k\Qt_y(\zeta_k\le t)=0$. If $\Lt V\le\lambda V$ only outside $B\subset\subset D$, replace $V$ by $V+c$ and $\lambda$ by $\lambda'$ so large that $\Lt V\le\lambda'(V+c)$ on $B$ as well (the coefficients and $V$ are bounded on $B$).

(ii) Let $\varsigma$ be the hitting time of $B$ ($\varsigma=\infty$ if $B=\emptyset$) and suppose $\Qt_{y_0}(\tau_D=\infty)=1$. Then $e^{-\lambda(t\wedge\varsigma)}V(Y_{t\wedge\varsigma})$ is a bounded submartingale, and since discounting can only decrease a positive value and increase a negative one, with $V^+:=\max\{V,0\}$,
\begin{align*}
V(y_0)&\le\E_{y_0}\big[e^{-\lambda(t\wedge\varsigma)}V(Y_{t\wedge\varsigma})\big]\le e^{-\lambda t}\sup V^++\E_{y_0}\big[e^{-\lambda\varsigma}V^+(Y_\varsigma);\varsigma\le t\big]\\
&\le e^{-\lambda t}\sup V^++\max\Big\{0,\sup_{\partial_DB}V\Big\}\xrightarrow[t\to\infty]{}\max\Big\{0,\sup_{\partial_DB}V\Big\} ,
\end{align*}
contradicting $V(y_0)>\max\{0,\sup_{\partial_DB}V\}$. The last assertion is \cite[Thm.~1.15.1]{Pin95}. (Without the positivity requirement the statement is false: for Brownian motion on $\R$, $B=(-\infty,0]$, $\lambda=1$ and $V(x)=-1-e^{-x}$ one has $LV\ge V$, $V(0)<V(1)$, and no explosion.)
\end{proof}

\begin{proof}[Proof of Proposition~\ref{prop:closed-state}]
On the common canonical space for the closed-state-space model, let $\tau_n$ be the first exit of the factor coordinate from $E\cap\{|y|<n\}$. The same path functional is used under both laws. The stopped integrand is bounded, so $S_{T\wedge\tau_n}$ defines a probability measure. Girsanov changes the drift to $\mu+\sigma\theta$ and leaves the finite-variation boundary term, including its support and direction constraints, unchanged. Uniqueness in law therefore gives
\[
\E_y[S_T\mathbf1_{\{\tau_n>T\}}]=\Qt_y(\tau_n>T).
\]
Since the original process is conservative, $\tau_n\uparrow\infty$; passing to the limit proves \eqref{eq:closed-mass}. The same proof applies to a compact exhaustion when $E$ has several boundary faces.
\end{proof}

\clearpage
\section{Overview of the martingale criteria}\label{app:hierarchy}

\begin{figure}[!ht]
\centering
\resizebox{!}{0.80\textheight}{%
\begin{tikzpicture}[
  font=\footnotesize,
  box/.style={draw=black, text=black, fill=white, line width=0.6pt, rounded corners=2pt, align=flush center, inner sep=4pt, text width=#1},
  exact/.style={box=#1, double, double distance=0.8pt},
  cond/.style={box=#1, solid, line width=0.8pt},
  suff/.style={box=#1, dashed, line width=0.8pt},
  tool/.style={box=#1, densely dotted, line width=0.8pt},
  ex/.style={box=#1, fill=white, font=\scriptsize},
  head/.style={font=\footnotesize\bfseries},
  arr/.style={-{Latex[length=2mm]}, thick},
  darr/.style={-{Latex[length=2mm]}, thick, black},
  ]
\node[head] (title) at (0,0) {Is $S=\mathcal E\big(\int\theta(Y)\cdot dW\big)$ a true martingale?\quad ($dY=\mu(Y)dt+\sigma(Y)dW$ on $D$, $\ \tilde\mu:=\mu+\sigma\theta$)};
\node[exact=4.9cm] (A) at (-5.6,-1.7) {\textbf{Lemma \ref{thm:main}} (probabilistic)\\ $S$ true martingale $\iff$ the $\tilde L$-diffusion never leaves $D$\\[1pt] $\mathbb E_y[S_T]=\tilde{\mathbb Q}_y(\tau_D>T)$};
\node[exact=4.9cm] (B) at (0,-1.7) {\textbf{Lemma \ref{thm:khas}} (Khasminskii)\\ $S$ true martingale $\iff$ $\tilde Lu=\lambda u$ has no bounded positive solution in $D$ ($\lambda>0$)};
\node[exact=4.9cm] (C) at (5.6,-1.7) {\textbf{Rem.\ \ref{rem:pde}, App.\ \ref{sec:pde}} (PDE)\\ $S$ true martingale $\iff$ the Cauchy problem $\partial_tm=\tilde Lm$, $m(0,\cdot)=0$, has only the zero bounded solution};
\draw[arr] (title) -- (A.north); \draw[arr] (title) -- (B.north); \draw[arr] (title) -- (C.north);
\draw[<->,thick] (A) -- (B); \draw[<->,thick] (B) -- (C);
\node[tool=14cm] (T) at (0,-4.0) {\textbf{Lemma \ref{thm:lyap} + Lemma \ref{thm:main}:} Lyapunov tests for martingality\\
$V\ge0$, $\tilde LV\le\lambda V$, $V\to\infty$ at $\partial D$ and infinity $\Rightarrow$ no exit $\Rightarrow$ $S$ true martingale.\\
A bounded subsolution satisfying the exit condition in Lemma~\ref{thm:lyap}(ii)\\
$\Rightarrow$ positive probability of exit $\Rightarrow$ $S$ strict local martingale.\\
\textbf{The coefficient tests below examine the possible exits separately.}};
\draw[arr] (B.south) -- (T.north);
\node[head] (h1) at (-5.6,-5.65) {smooth boundary points $\partial D^{(2)}$};
\node[head] (h2) at (0,-5.65) {infinity};
\node[head] (h3) at (5.6,-5.65) {singular strata (corners, axes)};
\coordinate (j1) at ($(T.south)!0.5!(h1.north)$);
\draw[arr] ([xshift=1.2cm]T.south west) |- (j1 -| h1.north) -- (h1.north);
\draw[arr] (T.south) -- (h2.north);
\coordinate (j3) at ($(T.south)!0.5!(h3.north)$);
\draw[arr] ([xshift=-1.2cm]T.south east) |- (j3 -| h3.north) -- (h3.north);
\node[cond=4.9cm, below=2mm of h1] (F) {Fichera drift of $\tilde L$, $\nu$ outward normal:\\ $\tilde F=\langle\tilde\mu-\tfrac12\operatorname{div}a,\ \nu\rangle$\\[2pt]
\textbf{Prop.\ \ref{thm:attain}, Thm.\ \ref{thm:fichera}:} $\tilde F\le0$ is necessary for $S$ to be a true martingale.\\
$\tilde F\le0$, with (T) at boundary points with zero Fichera drift where $\theta$ is unbounded, excludes smooth-boundary exits.\\[2pt]
\textbf{Prop.\ \ref{prop:inherit}:} $\theta$ bounded near a boundary patch $\Rightarrow$ no exit through that patch, hence no loss of mass there.\\
\textbf{Cor.\ \ref{cor:inherit}:} $\tilde F=F\le0$ holds automatically on such smooth patches.};
\node[cond=4.9cm, below=2mm of h2] (K) {Khasminskii's radial test, $r=\tfrac12|y|^2$:\\ $\mathcal A=\langle y,ay\rangle$, $\mathcal B=(\operatorname{tr}a+2\langle y,\tilde\mu\rangle)/\mathcal A$\\[2pt]
\textbf{Lemma \ref{thm:radial-sandwich}:} on large spheres, $A_-\le\mathcal A\le A_+$, $B_-\le\mathcal B\le B_+$.\\
$v_{A_+,B_+}(\infty)=\infty$ $\Rightarrow$ no explosion at infinity;\\
$v_{A_-,B_-}(\infty)<\infty$ $\Rightarrow$ $S$ strict local martingale.\\[2pt]
\textbf{Prop.\ \ref{thm:radial-iff}:} if exits through finite boundary points are excluded and $\mathcal A=A(r)$, $\mathcal B=B(r)$,\\
$S$ true martingale\\ $\iff v_{A,B}(\infty)=\infty$.\\[2pt]
\textbf{Prop.\ \ref{cor:osgood}:} linear growth bound $\Rightarrow$ no explosion at infinity;\\
Osgood drift with the stated noise bounds $\Rightarrow$ $S$ strict local martingale.};
\node[exact=4.9cm, below=3mm of K] (G) {\textbf{Prop.\ \ref{prop:adapted-infinity}, Cor.\ \ref{cor:variable-clock}:}\\
On $D=\mathbb R^d$, a $C^2$ scalar function $P$ with compact sublevel sets satisfies, outside a compact set,\\
$\langle a\nabla P,\nabla P\rangle=qA(P)$,\\
$2\tilde LP=qA(P)B(P)$.\\
If $0<q_-\le q/g(P)\le q_+<\infty$, then\\
$S$ true martingale\\ $\iff v_{gA,B}(\infty)=\infty$.};
\draw[darr] (K) -- (G);
\node[suff=4.9cm, below=2mm of h3] (P) {Projection onto the stratum (Ramasubramanian):\\ $\beta(r)=\inf_{|\pi(y)|=r}\frac{B_\pi-A_\pi+C_\pi}{A_\pi}$\\[2pt]
\textbf{Prop.\ \ref{prop:sing}:} $\beta(s)\ge1$ near $0$ $\Rightarrow$ no exit through this boundary piece $\Rightarrow$ no loss of mass there.\\[2pt]
Together with non-exit at the rest of the boundary and at infinity, this proves that $S$ is a true martingale.\\[2pt]
Sufficient only: failure of this condition does not imply strict local martingality.};
\node[ex=4.9cm, below=3mm of F] (F2) {\textbf{recovers / yields:} $3/2$ and $4/2$ models (Ex.\ \ref{ex:32}, \ref{ex:42}); two-factor inverse-CIR (Ex.\ \ref{ex:twofactor}); faces of the volatility-stabilized diffusion (Ex.\ \ref{ex:VS}); CIR and Wishart boundary non-attainment conditions (Sec.\ \ref{sec:affine})};
\node[ex=4.9cm, below=3mm of G] (K2) {\textbf{recovers / yields:} Dandapani--Protter Thms 6, 8, 10; Lions--Musiela Thms 2.4, 3.1; isotropic power-volatility model (Ex.\ \ref{ex:iso}); anisotropic model (Ex.\ \ref{ex:anisotropic}); conservativeness of affine diffusions (Sec.\ \ref{sec:affine})};
\node[ex=4.9cm, below=3mm of P] (P2) {\textbf{recovers / yields:} origin of the volatility-stabilized diffusion (Ex.\ \ref{ex:VS}); boundary pieces of affine models (Sec.\ \ref{sec:affine})};
\draw[darr] (F) -- (F2); \draw[darr] (G) -- (K2); \draw[darr] (P) -- (P2);
\node[exact=4.9cm, below=5mm of F2] (U) {\textbf{Coupled logarithmic model}\\
\textbf{Thm.\ \ref{thm:coupled-log}:} $S$ true martingale on every finite horizon\\ $\iff\sup c\le\tfrac12$.\\[3pt]
\textbf{Cor.\ \ref{cor:coupled-log-ui}: infinite horizon}\\
$S$ uniformly integrable martingale on $[0,\infty)$ $\iff c\equiv0$.\\
$c\not\equiv0\Rightarrow S_t\to0$ $\mathbb P$-a.s.};
\coordinate (log-entry) at ([yshift=-7mm]U.north west);
\coordinate (log-branch) at ([xshift=-9mm]A.west);
\draw[arr] (A.west) -- (log-branch)
  -- node[midway, above, sloped, font=\scriptsize, align=center]
  {Scalar comparison and Feller test: non-exit\\
   Logarithmic barrier: positive-probability exit}
  (log-branch |- log-entry) -- (log-entry);
\node[tool=16cm, below=4mm of K2] (combine) {\textbf{Combining the tests:} excluding all exits through finite boundary points and explosion at infinity proves that $S$ is a true martingale. Any positive probability of exit gives $\mathbb E_y[S_T]<1$ for some finite $T$.};
\node[exact=16cm, below=3mm of combine] (D1) {$d=1$: Feller's test at both endpoints gives an exact martingality criterion (Mijatović--Urusov, Bayraktar--Kardaras--Xing)};
\node[below=3mm of D1, font=\scriptsize, align=center] {\tikz{\node[exact=2mm, inner sep=2pt]{};} if and only if\qquad
\tikz{\node[cond=2mm, inner sep=2pt]{};} coefficient tests\qquad
\tikz{\node[suff=2mm, inner sep=2pt]{};} sufficient only\qquad
\tikz{\node[tool=2mm, inner sep=2pt]{};} tool};
\end{tikzpicture}
}
\caption{Martingality on every finite horizon, under Assumption~\ref{ass:standing} and the hypotheses of the cited results. The tests at infinity concern unbounded domains.}
\label{fig:hierarchy}
\end{figure}
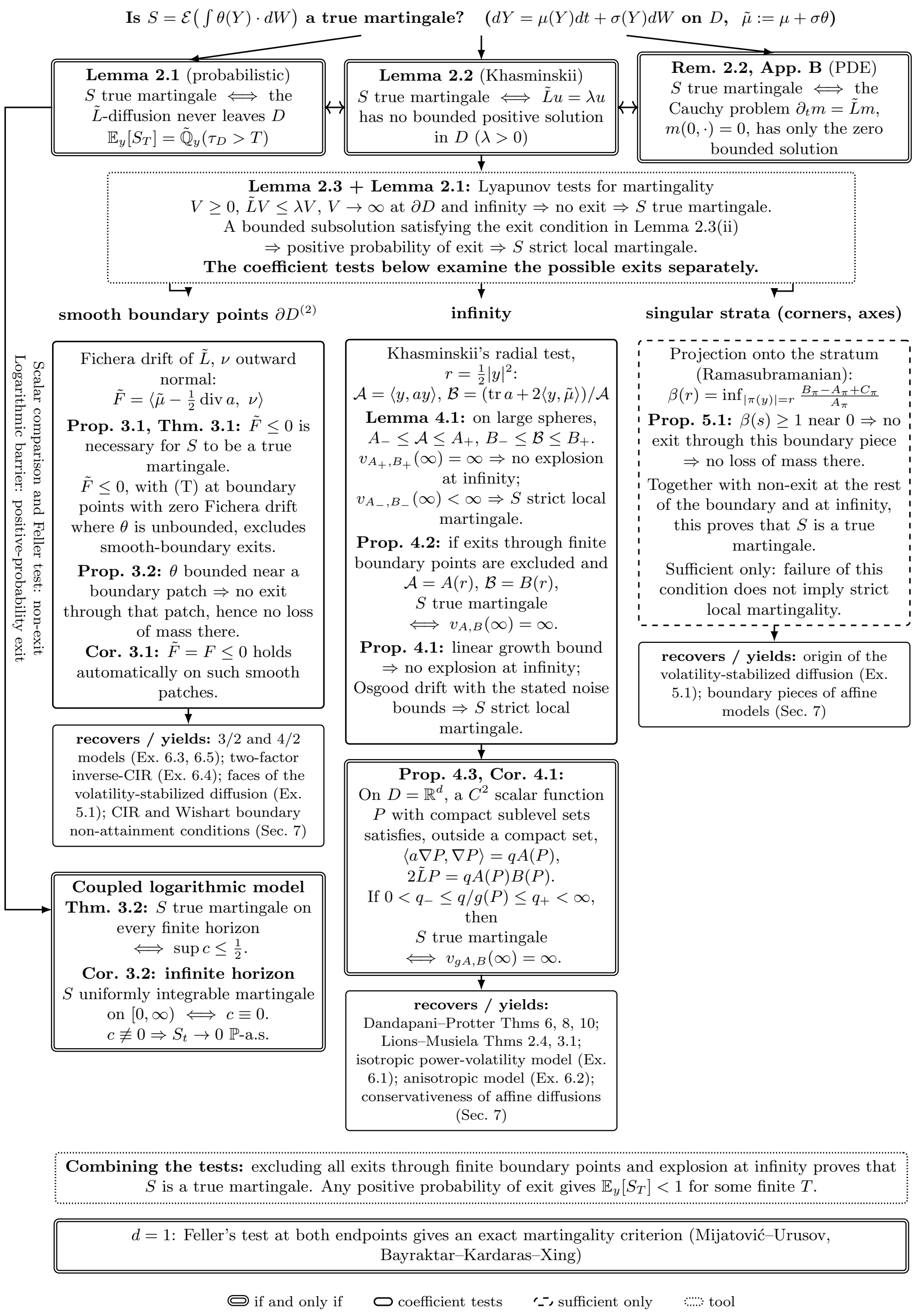
\clearpage

\bibliographystyle{amsplain}
\bibliography{v5-refs}
\end{document}